\documentclass[a4paper,mathserif,11pt]{article}
\usepackage{a4wide}

\usepackage[utf8]{inputenc}
\usepackage[T1]{fontenc} 
\usepackage{amsmath}
\usepackage{amsfonts}
\usepackage{amssymb}
\usepackage{amsthm}
\usepackage{mathabx}
\usepackage{color}
\usepackage{esint}
\usepackage{enumerate}
\usepackage{dsfont}
\usepackage{psfrag}
\usepackage{stmaryrd}
\usepackage{hyperref}
\usepackage{graphicx}
\usepackage{auto-pst-pdf}
\usepackage{graphics}
\usepackage{epstopdf}
\usepackage{epsfig}

\newtheorem{theorem}{Theorem}[section]
\newtheorem{remark}[theorem]{Remark}

\newcommand{\RR}{\mathbb{R}}
\newcommand{\NN}{\mathbb{N}}
\newcommand{\Aa}{\mathbb{A}}
\newcommand{\D}{\mathbb{D}}
\newcommand{\eps}{\varepsilon}
\newcommand{\dis}{\displaystyle}
\newcommand{\intO}{\int_{\Omega}}
\newcommand{\nab}{\boldsymbol{\nabla}}
\newcommand{\bn}{\boldsymbol{n}}

\newcommand{\bx}{\boldsymbol{x}}
\newcommand{\by}{\boldsymbol{y}}
\newcommand{\bs}{\boldsymbol{s}}
\newcommand{\bw}{\boldsymbol{w}}
\newcommand{\bq}{\boldsymbol{q}}
\newcommand{\bQ}{\boldsymbol{Q}}
\newcommand{\bE}{\boldsymbol{E}}
\newcommand{\bp}{\boldsymbol{p}}
\newcommand{\bX}{\boldsymbol{X}}
\newcommand{\bA}{\boldsymbol{A}}

\newcommand{\be}{\boldsymbol{e}}
\newcommand{\btau}{\boldsymbol{\tau}}
\newcommand{\M}{\mathbb{M}}
\newcommand{\mT}{\mathcal{T}}
\newcommand{\mR}{\mathcal{R}}

\newcommand{\mF}{\mathcal{F}}
\newcommand{\mH}{\mathcal{H}}
\newcommand{\Om}{\Omega}
\newcommand{\cQ}{\mathcal{Q}}
\newcommand{\cR}{\mathcal{R}}
\newcommand{\II}{\mathbb{I}}

\newcommand{\vertiii}[1]{{\left\vert\kern-0.25ex\left\vert\kern-0.25ex\left\vert #1 
    \right\vert\kern-0.25ex\right\vert\kern-0.25ex\right\vert}}

\begin{document}
\date{\today}
\author{Ludovic Chamoin$^{1,2}$ and Frédéric Legoll$^{3,2}$
\\ \\
{\footnotesize $^1$ LMT (ENS Cachan, CNRS, Universit\'e Paris-Saclay),}\\
{\footnotesize 61 avenue du Pr\'esident Wilson, 94235 Cachan, France}\\
{\footnotesize $^2$ Inria Paris, MATHERIALS project-team,}\\
{\footnotesize 2 rue Simone Iff, CS 42112, 75589 Paris Cedex 12, France}\\
{\footnotesize $^3$ Laboratoire Navier (ENPC, Paris-Est University),}\\
{\footnotesize 6 et 8 avenue Blaise Pascal, 77455 Marne-La-Vall\'ee Cedex 2, France}
}

\title{\textit{A posteriori} error estimation and adaptive strategy for the control of MsFEM computations}

\maketitle

\begin{abstract}
We introduce quantitative and robust tools to control the numerical accuracy in simulations performed using the Multiscale Finite Element Method (MsFEM). First, we propose a guaranteed and fully computable \textit{a posteriori} error estimate for the global error measured in the energy norm. It is based on dual analysis and the Constitutive Relation Error (CRE) concept, with recovery of equilibrated fluxes from the approximate MsFEM solution. Second, the estimate is split into several indicators, associated to the various MsFEM error sources, in order to drive an adaptive procedure. The overall strategy thus enables to automatically identify an appropriate trade-off between accuracy and computational cost in the MsFEM numerical simulations. Furthermore, the strategy is compatible with the offline/online paradigm of MsFEM. The performances of our approach are demonstrated on several numerical experiments.
\end{abstract}

\section{Introduction and objectives}

Developing multiscale numerical methods for mechanical problems with highly heterogeneous material structure is an increasingly active research field. Such methods are becoming a standard approach in Material Sciences and Computational Mechanics. As an alternative to solving a full fine-scale problem, with a usually prohibitive computational cost, multiscale modeling aims at linking the different scales for the accurate description of physical phenomena and/or the prediction of macroscopic properties (effective conductivity, elastic moduli, \dots). One of the goals is to capture the impact of the smaller scales on the larger scales, in order to observe the influence of the microscopic structure on the macroscopic behavior.

In the last 20 years, a lot of efforts have been put in the design of multiscale approaches for elliptic PDEs. They are inspired by the framework of the Finite Element Method (FEM) but take into account scale separation between macroscopic and microscopic features. A pioneering work using multiscale finite element basis functions is reported in~\cite{BAB94}. It can be considered as an extension of an earlier development, the Generalized Finite Element Method (GFEM)~\cite{BAB83,STR01,FIS04}. All these contributions share the idea of adapting the FE space to the particular fine-scale features of the problem by means of handbook functions and partition of unity. Other classical multiscale methods include the Multiscale Finite Element Method (MsFEM) developed in~\cite{HOU97,HOU99} and also based on modified basis functions obtained from fine-scale equations (we refer to~\cite{EFE09} for a review), the Variational Multiscale Method (VMS) developed in~\cite{HUG98}, the sparse FEM introduced in~\cite{HOA05}, multigrid methods~\cite{HAC85}, and multilevel finite elements methods such as the Heterogeneous Multiscale Method (HMM) proposed in~\cite{E03,ABD12} or the FE$^2$ method proposed in~\cite{FEY03}. We also mention the large literature on numerical homogenization from RVEs developed in the Structural Mechanics and Engineering communities to obtain macroscopic phenomenological constitutive laws~\cite{SUQ87}.

In this context, we focus on MsFEM which is a powerful tool to capture information at various scales without resorting to a full fine-scale computation~\cite{EFE09,LEB14}. It defines an approximate solution in a finite dimensional space, related to a macroscopic mesh and generated by basis functions which encode details of the fine-scale heterogeneities. MsFEM performs the computations in a two-stage procedure: (i) an \textit{offline} stage in which basis functions are computed solving local fine-scale problems; (ii) an \textit{online} stage in which an inexpensive Galerkin approximation problem is solved. We refer to Section~\ref{section:MSFEM} below for more details. As for all numerical methods, a crucial issue is to control the accuracy of the obtained approximate numerical solution. In particular, it is important to understand and master how errors generated in problems at different scales influence the overall error in the MsFEM solution. Furthermore, the MsFEM numerical model may be costly, since fine-scale \textit{offline} computations may not be possible over the whole structure, even with parallel computing, in practical situations (limited resources available). Consequently, a major problem is to identify an appropriate trade-off between accuracy and computational cost. The ultimate goal is then to design MsFEM algorithms such that: (i) a given precision is attained at the end of the simulation; (ii) the computational work needed is as small as possible (computing resources are optimized).

\medskip

For many years, numerous research works have been conducted to develop effective tools for the estimation of discretization error and the adaptation of discretization parameters, mainly in the monoscale framework of FEM. A review can be found in~\cite{VER96,AIN00,LAD04,CHA15}. In contrast to \textit{a priori} estimates that merely enable to evaluate the convergence order (with respect to discretization parameters) of a given numerical method, \textit{a posteriori} error estimates provide quantitative information on the error as well as criteria for mesh adaptation. A particular and robust \textit{a posteriori} error estimation tool, using the concept of Constitutive Relation Error (CRE), has been used in the Computational Mechanics community for more than thirty years~\cite{LAD04}. It is based on a dual analysis approach~\cite{FRA72,ODE74,FRA01} and thus requires the computation (by post-processing the primal approximate solution field) of flux fields that verify equilibrium in a strong sense~\cite{PLE11}. This enables to derive, in a convenient way, fully computable (i.e. without any unknown multiplicative constant) and guaranteed error bounds on both the global solution (error in energy norm) and given quantities of interest (goal-oriented error estimation) using an adjoint problem~\cite{LAD10}.

In the multiscale framework, and for MsFEM-like methods particularly, the state-of-the-art is very different. On the one hand, a large literature has addressed \textit{a priori} error estimation for the last fifteen years, enabling to evaluate the convergence of the error with respect to the macroscopic mesh size and the characteristic length of the microscopic heterogeneities~\cite{HOU99,EFE00,HOU04,ALL06,EFE09,HES14}. On the other hand, despite the fact that \textit{a posteriori} error control and adaptivity have become an important issue for reliable and efficient multiscale computations, very few tools on this topic are available (see~\cite{STR06,STR07,LAR07,NOL08,ABD09,ABD11a,LAR11,PAL17}) and there are many open questions: quantitative assessment of error propagation across scales, relevant adaptive strategies, \dots Currently, and as far as the authors know, there is only one \textit{a posteriori} error estimate approach available for MsFEM approximations. It was recently proposed in~\cite{HEN14} using the explicit residuals method. This latter approach, splitting the error estimate with respect to the different error sources, is effective to drive an adaptive algorithm. However, it does not enable to quantify the error level accurately (but only up to some unknown constants), which is a drawback in engineering activities where certification of the accuracy and stopping criteria for adaptive algorithms are required. We also mention the \textit{a posteriori} error indicator 
proposed in~\cite{CHU14} for a variant of MsFEM, namely the Generalized Multiscale Finite Element Method (GMsFEM).

\medskip

In this work, we develop a robust \textit{a posteriori} error estimate, as well as an associated adaptive strategy, for MsFEM computations. To achieve this goal, the CRE concept is extended to the multiscale framework. We first show how equilibrated fluxes can be recovered from an extension of the method given for FEM in~\cite{LAD04}, so that a guaranteed and fully computable error estimate on the overall error can be derived. This is, to the best authors knowledge, the first estimate of this kind for MsFEM. In order to minimize the computational cost associated with error estimation, as well as to conserve the MsFEM philosophy, we pay attention to perform most of the additional computations in the \textit{offline} step and with multi-query constraints. Conforming and non-conforming (oversampling) variants of MsFEM are considered in this framework.

Second, we detail how CRE can be effectively used to localize and assess the various MsFEM error sources, before driving a robust adaptive algorithm from local error maps and with respect to a given error tolerance. The adaptive strategy enables to choose automatically the various discretization parameters which are involved in MsFEM. In particular, it circumvents \textit{offline} fine-scale computations over all elements of the macro mesh. Local remeshing being not suitable in the MsFEM context (since it would lead to costly recomputations of basis functions), we use a variant ($h$-version) of the technique proposed in~\cite{ALL06}, which avoids to go back to the \textit{offline} step when adapting the macro mesh.

Third, the proposed techniques for robust error estimation and multiscale adaptivity are illustrated on several examples of the literature. Numerical results confirm the practical interest of these techniques, in terms of both certification of the numerical solution and computational savings, when performing multiscale simulations using MsFEM.

\medskip

The article is organized as follows. The model problem is described in Section~\ref{section:modelpb}, where we also recall some elements of homogenization theory. Basics on MsFEM are recalled in Section~\ref{section:MSFEM}. The \textit{a posteriori} error estimation method based on CRE is introduced in Section~\ref{section:CRE}. The adaptive strategy, associated to indicators related to various error sources, is detailed in Section~\ref{section:adaptive}. Numerical experiments illustrating the performances of the proposed approach are discussed in Section~\ref{section:results}. Conclusions and prospects are eventually drawn in Section~\ref{section:conclusions}. 

\section{Model problem and homogenized limit}\label{section:modelpb}

\subsection{Problem of interest}

We consider a fine-scale elliptic problem, typically associated with heat transport in composite materials or Darcy flow in porous media, that is characterized by strong material heterogeneities. The characteristic scale of the material coefficients, denoted $\eps$, is hence small. The associated boundary value problem reads
\begin{equation}\label{modelpb}
-\nab \cdot [\Aa^\eps \nab u^\eps]=f \ \ \text{in $\Omega$}, \qquad u^\eps=0 \ \ \text{on $\Gamma_D$}, \qquad (\Aa^\eps \nab u^\eps) \cdot \bn = g \ \ \text{on $\Gamma_N$},
\end{equation}
where $\Omega$ is an open bounded domain of $\RR^d$. Its boundary $\Gamma = \partial \Omega$ is split into a subset $\Gamma_D \subset \partial \Omega$ (with measure $|\Gamma_D|>0$) where homogeneous Dirichlet boundary conditions are applied, and a subset $\Gamma_N = \partial \Omega \setminus \Gamma_D \subset \partial \Omega$ where Neumann boundary conditions are applied. We assume that the body loading $f \in L^2(\Omega)$ and the traction loading $g \in L^2(\Gamma_N)$ are given and slowly varying (they do not depend on $\eps$). The quantity $\bq^\eps=\Aa^\eps \nab u^\eps$ is the flux associated to $u^\eps$.

The second-order diffusion tensor $\Aa^\eps \in (L^{\infty}(\Omega))^{d\times d}$ is a rapidly oscillating function, $\eps$ representing the typical (presumably small) size of the heterogeneities. We assume that $\Aa^\eps$ is a symmetric matrix (see Remark~\ref{rem:non-sym} below for a discussion of that assumption) and that it is uniformly elliptic and bounded, in the sense that there exist $\alpha>0$ and $\beta>0$ such that
\begin{equation}
  \label{eq:elliptic}
\forall \eps \ge 0, \quad \forall \boldsymbol{\xi} \in \RR^d, \quad
\alpha |\boldsymbol{\xi}|^2 \le (\Aa^\eps(\bx)\boldsymbol{\xi}) \cdot \boldsymbol{\xi} \le \beta |\boldsymbol{\xi}|^2 \quad \text{a.e. in $\Omega$}.
\end{equation}
The weak formulation of~\eqref{modelpb} consists in finding $u^\eps \in V$ such that
\begin{equation}\label{modelpbweak}
\forall v \in V, \quad B^\eps(u^\eps,v) = F(v),
\end{equation}
where 
\begin{equation}
  \label{eq:def_V}
V = \mH^1_{0,\Gamma_D}(\Omega) = \left\{ v \in \mH^1(\Omega); \quad v=0 \ \ \text{on $\Gamma_D$} \right\}
\end{equation}
and 
\begin{equation*}
B^\eps(u,v) = \intO \left[ \Aa^\eps \, \nab u \right] \cdot \nab v, \qquad F(v) = \intO f \, v + \int_{\Gamma_N} g \, v.
\end{equation*}
The bilinear form $B^\eps$ is symmetric, continuous and coercive on $V$. It hence defines an inner product and induces the energy norm $\vertiii{v} = \sqrt{B^\eps(v,v)}$ on $V$. In what follows, for any $k \in \NN$ and any $v \in \mH^k(\Omega)$, we set $\| v \|_k = \| v \|_{\mH^k(\Omega)}$. Note in particular that $\| v \|_0 = \| v \|_{L^2(\Omega)}$.

\begin{remark}
\label{rem:non-sym}
We refer to~\cite{VER96,ERN10}, as well as~\cite{CHA16} and the extensive bibliography therein, for \textit{a posteriori} error estimations (in a single-scale framework) in the case when the bilinear form $B^\eps$ in~\eqref{modelpbweak} is non-symmetric.
\end{remark}

Let $\mT_H$ be a partition of $\Omega$. We denote $\dis H_K=\sup_{M,N \in K}|MN|$ the diameter of each element $K$, and $\dis H=\max_{K \in \mT_H}H_K$. Introducing the P1-Lagrange finite dimensional space 
\begin{equation}
\label{eq:def_VH}
V_H^0 = \left\{
\begin{array}{c}
  v \in C^0(\overline{\Omega}), \ \ v \in V, \ \ \text{$v_{|K}$ is a polynomial function}
  \\
  \text{of degree up to 1 on any element $K$}
\end{array}
\right\},
\end{equation}
the Finite Element (FE) approximation of~\eqref{modelpb} reads:
\begin{equation*}
\text{Find $u_H \in V_H^0$ such that, for any $v \in V_H^0$}, \quad B^\eps(u_H,v)=F(v).
\end{equation*}
Note that $V_H^0 \subset V$, and hence the approximation is conforming. Assuming that $u^\eps \in \mH^2(\Omega)$, standard FEM \textit{a priori} error estimates yield that
\begin{equation*}
\| u^\eps - u_H \|_1 \le C H \|u^\eps \|_2 \le C \frac{H}{\eps}\|f\|_0,
\end{equation*} 
where $C$ is independent of $H$ and $\eps$. We hence see that obtaining an accurate approximation $u_H$ requires to choose $H \ll \eps$. The number of degrees of freedom is then very large. This leads to a usually prohibitive computational complexity. Non-convergence to the correct solution in the limit $\eps \to 0$ is due to the fact that the small scale information is averaged out incorrectly in the classical FEM approach.

\subsection{Homogenization strategy}\label{section:homogenization}

From an engineering perspective, it may be sufficient to only predict the {\em macroscopic} properties of the solution to~\eqref{modelpb}. Powerful homogenization methods can be used for this purpose, based on a scale separation assumption, i.e. considering the regime $\eps \ll L$ where $L$ is the characteristic size of the macroscopic variations of $u^\eps$ inside $\Om$. The idea is to approximate the behavior of the heterogeneous medium with an effective averaged behavior at a macroscopic scale. From a mathematical point of view, homogenization consists in identifying the limit of the operator $\dis \big( -\nab \cdot [\Aa^\eps \nab \cdot] \big)^{-1}$ when $\eps$ tends to zero.

For the considered problem, homogenization theory (see the textbooks~\cite{BEN78,SAN80,JIK94,MUR97,TAR10} and the review article~\cite{KAN09}) yields the following result, if no further assumptions are made on $\Aa^\eps$. Up to the extraction of a subsequence, $u^\eps$ converges to some $u^0$ (strongly in the $L^2$-norm, weakly in the $\mH^1$-norm), solution to the macroscale (homogenized) problem 
\begin{equation*}
-\nab \cdot [\Aa^0 \nab u^0] = f \ \ \text{in $\Omega$}, \qquad u^0=0 \ \ \text{on $\Gamma_D$}, \qquad (\Aa^0 \nab u^0) \cdot \bn = g \ \ \text{on $\Gamma_N$},
\end{equation*}
where the second-order tensor $\Aa^0$ is symmetric and satisfies the same bounds~\eqref{eq:elliptic} as $\Aa^\eps$. The tensor $\Aa^0$ depends on the considered subsequence, but is independent of the loads $f$ and $g$.  

When $\Aa^\eps$ has a periodic (or locally periodic) structure with period $\eps$, the whole sequence $u^\eps$ converges to $u^0$, and the homogenized tensor $\Aa^0$ can be characterized from the solutions to \textit{cell problems}, i.e. boundary value problems posed on a domain of size $\eps^d$ and involving the fine-scale tensor $\Aa^\eps$. Indeed, introducing the fast variable $\by=\bx/\eps$ and assuming that $\Aa^\eps(\bx)=\Aa(\bx,\by)$ where the map $\by \mapsto \Aa(\bx,\by)$ is $Y$-periodic (where $Y=(0,1)^d$ is the unit cube), the coefficients of the homogenized tensor $\Aa^0$ are given by
\begin{equation*}
A^0_{ij}(\bx)=\int_Y\left(A_{ij}(\bx,\by)+\sum_{k=1}^d A_{ik}(\bx,\by)\frac{\partial \chi_j}{\partial y_k}(\bx,\by)\right)d\by,
\end{equation*}
where $\chi_j(\bx,\cdot)$ denotes the periodic solution (of mean zero) to the following cell (or corrector) problem:
\begin{equation}\label{pbcorr}
\nab_{\by}\cdot\big[\Aa(\bx,\cdot)(\be_j+\nab_{\by}\chi_j(\bx,\cdot))\big]=0 \quad \text{in $Y$}, \qquad \text{$\by \mapsto \chi_j(\bx,\by)$ is $Y$-periodic}.
\end{equation}
The vector set $\{\be_j\}_{1 \leq j \leq d}$ is the standard basis of $\RR^d$. This result can be formally obtained using the two-scale expansion $u^\eps(\bx)=u^0(\bx,\by)+\eps u^1(\bx,\by)+\eps^2 u^2(\bx,\by)+\dots$, where each map $\by \mapsto u^s(\bx,\by)$, for any $s \geq 0$, is assumed to be $Y$-periodic.

Furthermore, $u^\eps$ and $\nab u^\eps$ can be accurately approximated for small values of $\eps$ by
\begin{equation*}
u^\eps(\bx) \approx u^{\eps,1}(\bx) := u^0(\bx)+\eps\sum_{j=1}^d\chi_j(\bx,\by)\frac{\partial u^0}{\partial x_j}(\bx), \quad 
\nab u^\eps(\bx)\approx \nab u^0(\bx)+\sum_{j=1}^d \nab_{\by} \chi_j(\bx,\by)\frac{\partial u^0}{\partial x_j}(\bx),
\end{equation*}
in the sense that $\| u^\eps - u^{\eps,1} \|_1 = O(\sqrt{\eps})$. We note that the contribution of the first-order term $\dis \eps u^1(\bx,\by) = \eps\sum_{j=1}^d\chi_j(\bx,\by)\frac{\partial u^0}{\partial x_j}(\bx)$ is crucial for the approximation of $\nab u^\eps$: its gradient is of the same order as $\nab u^0$. The gradient of $u^\eps - u^0$ does not (strongly) converge to 0, but the gradient of $u^\eps - (u^0 + \eps u^1)$ does.

\begin{remark}
In general, the first-order approximation $u^0+\eps u^1$ of $u^\eps$ does not satisfy the homogeneous Dirichlet boundary conditions on $\Gamma_D$ that $u^\eps$ satisfies. These can be enforced through a first-order correction term $\eps \theta^\eps$, with $\theta^\eps$ verifying $\nab \cdot [\Aa^\eps \nab \theta^\eps]=0$ in $\Omega$ and $\theta^\eps=-u^1$ on $\Gamma_D$. The correction term is such that $\|\eps \theta^\eps \|_1 = O(\sqrt{\eps})$.
\end{remark}

For other structures of $\Aa^\eps$, dedicated numerical methods need to be devised to solve~\eqref{modelpb} with a computational cost comparable to that required to solve the homogenized problem on a coarse mesh of size $H$ independent of $\eps$. When only the homogenized solution is of interest (which implies that $\eps \ll 1$), numerical homogenization over Representative Volume Elements (RVE) may be performed. It consists in solving the equilibrium problem $\nab \cdot [\Aa^\eps \nab v^\eps]=0$ on the RVE, complemented with specific boundary conditions (Dirichlet, Neumann or periodic boundary conditions) which are in agreement with the Hill-Mandel lemma (equivalence of micro and macro energies). We first write that
\begin{multline*}
\frac{1}{|\text{RVE}(\bx,\ell)|} \int_{\text{RVE}(\bx,\ell)} \btau^\eps \cdot \nab v^\eps - \bQ^\eps(\bx,\ell) \cdot \bE^\eps(\bx,\ell)
\\=
\frac{1}{|\text{RVE}(\bx,\ell)|} \int_{\text{RVE}(\bx,\ell)} \left( \btau^\eps - \bQ^\eps(\bx,\ell) \right) \cdot \left( \nab v^\eps - \bE^\eps(\bx,\ell) \right)
\end{multline*}
with
\begin{equation*}
\bQ^\eps(\bx,\ell) = \frac{1}{|\text{RVE}(\bx,\ell)|} \int_{\text{RVE}(\bx,\ell)} \btau^\eps, \qquad \bE^\eps(\bx,\ell) = \frac{1}{|\text{RVE}(\bx,\ell)|} \int_{\text{RVE}(\bx,\ell)} \nab v^\eps,
\end{equation*}
where $\text{RVE}(\bx,\ell)$ is a RVE around the macroscopic point $\bx$ of diameter $\ell$, $v^\eps$ is the solution to the problem posed on the RVE and $\btau^\eps = \Aa^\eps \nab v^\eps$. We note that $\bQ^\eps(\bx,\ell)$ and $\bE^\eps(\bx,\ell)$ are macroscopic quantities obtained through spatial averaging. Using an integration by parts, we deduce that
\begin{multline*}
\frac{1}{|\text{RVE}(\bx,\ell)|} \int_{\text{RVE}(\bx,\ell)} \btau^\eps \cdot \nab v^\eps - \bQ^\eps(\bx,\ell) \cdot \bE^\eps(\bx,\ell)
\\
=
\frac{1}{|\text{RVE}(\bx,\ell)|}\int_{\partial_{\text{RVE}(\bx,\ell)}}\left(v^\eps - \bE^\eps(\bx,\ell)\cdot \bX\right) \left(\btau^\eps-\bQ^\eps(\bx,\ell)\right) \cdot \bn,
\end{multline*}
where $\partial_{\text{RVE}(\bx,\ell)}$ is the boundary of the RVE, $\bn$ is the outgoing normal vector and $\bX$ is the local position vector over the RVE. Using the boundary conditions imposed on the RVE problem, we eventually deduce the equivalence of the micro and macro energies:
$$
\frac{1}{|\text{RVE}(\bx,\ell)|} \int_{\text{RVE}(\bx,\ell)} \btau^\eps \cdot \nab v^\eps - \bQ^\eps(\bx,\ell) \cdot \bE^\eps(\bx,\ell) = 0.
$$
By linearity, we have that $\bQ^\eps(\bx,\ell) = \overline{\Aa}(\bx,\ell,\eps) \, \bE^\eps(\bx,\ell)$ for some matrix $\overline{\Aa}(\bx,\ell,\eps)$. The homogenization approach is valid provided the characteristic length $\ell$ of the RVEs is such that $\eps \ll \ell \ll L$. We thus expect that $\dis \lim_{\ell \to 0} \lim_{\eps \to 0} \overline{\Aa}(\bx,\ell,\eps) = \Aa^0(\bx)$. A single RVE may be considered, for instance when calibrating macroscopic phenomenological constitutive laws, but a set of RVEs can also be defined as in multilevel FEM~\cite{E03,FEY03}, where a specific RVE is attached to each quadrature point of the macroscopic mesh.

\medskip

We consider in the following an alternative multiscale numerical approach, the Multiscale Finite Element Method (MsFEM), which has the capability to approximate the full fluctuating solution $u^\eps$ (with fine-scale details) in a convenient manner. This is performed without any assumption on the structure of $\Aa^\eps$, and without any \textit{a priori} knowledge on the form of the homogenized model associated with the heterogeneous model of interest.

\section{Basics on MsFEM}\label{section:MSFEM}

We recall here some basic elements on several variants of MsFEM: the initial conforming variant, the oversampling variant (see Section~\ref{section:oversampling}) and the higher-order variant introduced by G. Allaire and R. Brizzi (see Section~\ref{section:higherorder}). The reader familiar with MsFEM may easily skip this section and directly proceed to the {\it a posteriori} error estimation described in Section~\ref{section:CRE}.

\subsection{General idea}

\begin{figure}[htbp]
\begin{center}
  \includegraphics[width=40mm]{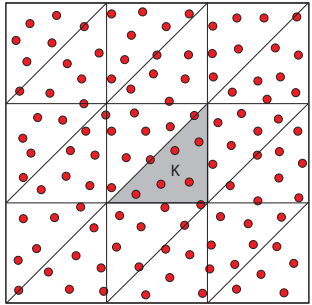}
\end{center}
\caption{A coarse mesh of the domain $\Omega$ is introduced, with elements of diameter $H$ much larger than the small characteristic size $\eps$ of the heterogeneities (here the diameter of the inclusions). On each coarse element, we solve a local problem, using in practice a discretization at the size $h$ adapted to the heterogeneities.}
\label{fig:maillages}
\end{figure}

MsFEM was introduced in~\cite{HOU97} as a multiscale approach for simulating problems in highly heterogeneous media which are intractable by conventional FEM.  The MsFEM approach is a particularly attractive tool to capture both microscale and macroscale features of the solution. It is self-consistent in the sense that it only uses the knowledge of $\Aa^\eps$. The main idea in the MsFEM approach is to construct a set $\{\phi^\eps_i\}_{1 \leq i \leq I}$ of local multiscale FE basis functions that encode small scale information within each element of a coarse mesh $\mT_H$ ($H\gg \eps$). The basis functions $\phi^\eps_i$, associated with each node $i$ of the coarse mesh $\mT_H$ (except those with prescribed Dirichlet boundary conditions), are adapted to the local properties of the operator. They are pre-computed in an \textit{offline} stage, over each element $K$ of the coarse mesh (see Figure~\ref{fig:maillages}), as solutions to local elliptic equations of the form
\begin{equation}\label{eq:localpbmsfem}
\nab\cdot \big[\Aa^\eps \nab \phi^\eps_i\big]=0 \quad \text{in $K$},
\end{equation}
complemented with various boundary conditions discussed below (see e.g.~\eqref{eq:linearBC}). We note that the basis functions are independent of the loads $f$ and $g$ present in the right-hand side of~\eqref{modelpbweak}. In addition, the problems~\eqref{eq:localpbmsfem} are decoupled one from each other. The MsFEM approach is thus well-suited to massively parallel computers. This is an essential feature of the approach.

\begin{remark}
In practice, the local problems~\eqref{eq:localpbmsfem} are solved numerically using a fine mesh of characteristic size $h < \eps$. 
\end{remark}

Once the multiscale basis functions are computed, the MsFEM strategy is the same as that of a classical FE approach. It consists in performing a Galerkin approximation of~\eqref{modelpbweak} on the finite dimensional space 
$$
V^\eps_H = \text{Span} \{\phi^\eps_i, \ 1 \leq i \leq I \}.
$$
In an \textit{online} stage, we thus look for the approximate solution $u^\eps_H \in V^\eps_H$ such that
\begin{equation}\label{eq:GalerkinMsFEM}
\forall v \in V^\eps_H, \quad B^\eps(u^\eps_H,v)=F(v).
\end{equation}
The small scale information incorporated in the basis functions is thus brought to the large scales through couplings in the global stiffness matrix. The assembly of this matrix is inexpensive since it reuses local matrices computed and stored in the \textit{offline} stage to solve~\eqref{eq:localpbmsfem}.

\medskip

As regards memory requirements and computational cost, it is straightforward to observe that they are considerably reduced when using the MsFEM approach compared to a direct FE approach. Indeed, considering regular grids with $N$ (resp. $M$) subdivisions in each dimension, let us denote by $N^d$ the number of elements $K$ in $\mT_H$ and by $M^d$ the number of subcell elements in each element $K$ that are used for the practical construction of the basis functions. There are thus $(MN)^d$ elements in the fine grid. Assume that we have $P$ available processors with similar computing load. We then see that: (i) the required computing memory on each processor for traditional FEM is of the order of $(MN)^d/P$; (ii) assuming that $P \geq N^d$ (i.e. that we have more processors than coarse elements in $\mT_H$), the required computing memory on each processor for MsFEM is of the order of $M^d+N^d/P$. The MsFEM approach therefore systematically brings a saving factor in terms of memory.

For a comparison on the computational cost, let us consider a single processor ($P=1$) and let us assume that the cost to solve a linear system with $\mathcal{N}$ unknowns on this processor scales as $\mathcal{N}^\alpha$ for some $\alpha > 1$ (the value of $\alpha$ depends on the solver used). The CPU cost is thus of the order of $(MN)^{d\alpha}$ for the FE approach and of the order of $N^{d\alpha}+ k N^d M^{d\alpha}$ for the MsFEM approach, where $k$ is the number of local problems of the form~\eqref{eq:localpbmsfem} to solve (i.e. the number of nodes) on each element $K$. The MsFEM approach thus brings a saving factor in terms of computational cost as well.

\begin{remark}
It is common in practice to perform multiple runs for the same medium, with different boundary conditions or loadings $f$ in~\eqref{modelpb}. Since the local problems~\eqref{eq:localpbmsfem} are independent of these, the \textit{offline} stage in MsFEM is done only once in such a multi-query context. Denoting by $n_R$ the number of runs, the CPU cost becomes of the order of $n_R(MN)^{d\alpha}$ for FEM and $n_RN^{d\alpha}+ k N^d M^{d\alpha}$ for MsFEM. The MsFEM approach is even more attractive in terms of computational costs in such a multi-query context.
\end{remark}

The philosophy of MsFEM, in which the fine-scale equations are solved inside coarser elements and are thus totally decoupled, is closely related to that of the VMS method~\cite{HUG98} and the GFEM method~\cite{STR01}. However, they differ on several points, in particular in the definition of the local fine-scale problems (computational domain, boundary conditions):
\begin{itemize}
\item The basic idea in VMS is to introduce an additive decomposition of the solution into coarse and fine scale contributions, under the form $u^\eps=\overline{u}+u^f$. It corresponds to the decomposition of the functional space $V$ in the direct sum $V=\overline{V}\oplus V^f$, where $\overline{V}$ is a coarse-scale space and $V^f$ is an infinite dimensional fine-scale space. This decomposition leads to two coupled (coarse/fine) problems. The part $u^f$, solution to the fine-scale problem involving the coarse-scale residual, is searched analytically using Green's functions (or element-wise approximations of these Green's functions; bubble functions are often used in practice). It is then eliminated from the coarse-scale problem. This leads to a modified coarse-scale equation, which takes the effects of the fine scales into account, and where the only unknown is $\overline{u}$. This coarse-scale contribution $\overline{u}$ is then numerically approximated using a finite-dimensional space $\overline{V}_H$. In the adaptive VMS proposed in~\cite{LAR07}, fine-scale equations are decoupled by a partition of unity (defined from the coarse mesh $\mT_H$) and solved on vertex patches with homogeneous Dirichlet boundary conditions.
\item In GFEM, some so-called handbook functions are constructed as fine-scale numerical solutions to local boundary value problems posed in vertex patches, and complemented with Neumann boundary conditions. These functions are then introduced, by means of a partition of unity defined from the coarse mesh $\mT_H$, in order to locally enrich the classical FEM basis relating to the topology of the material structure or the physical domain.
\end{itemize}

\subsection{Conforming MsFEM approximation}\label{section:conforming}

The boundary conditions associated to~\eqref{eq:localpbmsfem} in the construction of basis functions are of critical importance, since they determine the behavior of the numerical approximation $u_H^\eps$ on the edges $\dis \cup_{K \in \mT_H} \partial K$ of the mesh and how this behavior may be different from that of the exact solution $u^\eps$. It has been observed that the choice of boundary conditions in~\eqref{eq:localpbmsfem} can significantly affect the accuracy of the MsFEM approach (see e.g.~\cite{HOU97,LEB14}).

A first choice, leading to a conforming numerical discretization, is to impose a linear evolution of $\phi^\eps_i$ along $\partial K$ as for classical first-order FE basis functions. We thus introduce the piecewise affine FE basis functions $\phi^0_i$ over the coarse mesh $\mT_H$, $1 \leq i \leq I$, recall that $V_H^0 = \text{Span} \{\phi^0_i, \ 1 \leq i \leq I \}$ (where $V_H^0$ is defined by~\eqref{eq:def_VH}), and consider the local problems
\begin{equation}\label{eq:linearBC}
\nab\cdot \big[\Aa^\eps \nab \phi^\eps_i\big]=0 \quad \text{in $K$}, \qquad \phi^\eps_i = \phi^0_i \quad \text{on $\partial K$}.
\end{equation}
Note that the support of $\phi^\eps_i$ is identical to that of $\phi^0_i$.

\begin{remark}
Another choice of Dirichlet boundary conditions on $\partial K$, still leading to a conforming approach, consists in setting $\phi^\eps_i = \lambda^\eps_i$ on $\partial K$, where $\lambda^\eps_i$ is a rapidly oscillating function solution to a reduced elliptic problem on each edge of $\partial K$ (see~\cite{HOU99}). This reduced problem is obtained from $\nab\cdot \big[\Aa^\eps\nab \phi^\eps_i\big]=0$ by deleting terms with partial derivatives in the direction normal to $\partial K$. 
\end{remark}

It is straightforward to show that the multiscale basis functions $\phi^\eps_i$ defined by~\eqref{eq:linearBC} verify a partition of unity property. Indeed, from~\eqref{eq:linearBC}, we get for any element $K \in \mT_H$ that
\begin{equation*}
\nab \cdot \left[ \Aa^\eps \nab \left( \sum_i \phi^\eps_i \right) \right]=0 \quad \text{in $K$}, \qquad \sum_i \phi^\eps_i = 1 \quad \text{on $\partial K$},
\end{equation*}
where the sum on $i$ runs either on $1 \leq i \leq I$, or on the vertices of $K$. This implies that
\begin{equation}
\label{eq:pum}
\sum_i \phi^\eps_i=1 \quad \text{in $K$}.
\end{equation}
This property will be used in Section~\ref{equilitractionsconf} to derive the \textit{a posteriori} error estimate. 

\medskip

In the regime of interest $H \geq \eps$, and assuming that the oscillations in the material behavior are {\em periodic}, that is that $\Aa^\eps(\bx)=\Aa_{\rm per}(\bx/\eps)$ for a fixed periodic tensor $\Aa_{\rm per}$, the following {\it a priori} error estimate is shown in~\cite{HOU99}:
\begin{equation}\label{eq:aprioriboundMsFEM}
\| u^\eps - u^\eps_H \|_1 \leq C \left( H + \sqrt{\eps/H} + \sqrt{\eps} \right),
\end{equation}
where $C$ is a positive constant independent of $\eps$ and $H$ (note that some technical regularity assumptions are required on $\Aa_{\rm per}$ and on the homogenized solution for~\eqref{eq:aprioriboundMsFEM} to hold).

The estimate~\eqref{eq:aprioriboundMsFEM} shows that the MsFEM approach converges to the correct solution in the homogenization limit (namely when first $\eps \to 0$ and next $H \to 0$, which is the regime of interest). The ratio $\eps/H$ reflects the so-called resonance error, which occurs when $H \approx \eps$. In practice, this is often the leading error term.

The key to derive~\eqref{eq:aprioriboundMsFEM} is to use the periodic homogenization theory (recalled in Section~\ref{section:homogenization}), which provides (when $\eps \ll 1$) a detailed description of $u^\eps$ and of the basis functions $\phi^\eps_i$ in terms of a two-scale expansion. The basis functions being constructed using the {\em same} tensor $\Aa^\eps$ as that of the problem of interest, this implies that all the $\phi^\eps_i$ have the same asymptotic structure as $u^\eps$, which eventually allows for a small error between $u^\eps_H$ and $u^\eps$. The term proportional to $\dis \sqrt{\eps/H}$ (resp. $\dis \sqrt{\eps}$) in the estimate is due to the first-order correction term in the two-scale expansion of the basis functions $\phi^\eps_i$ (resp. of $u^\eps$), whereas the term proportional to $H$ comes from the classical {\it a priori} error estimate on the approximation of the homogenized solution $u^0$ by the piecewise affine basis functions $\phi^0_i$.

\subsection{Nonconforming MsFEM approximation: the oversampling technique}\label{section:oversampling}

The limited accuracy of the MsFEM approach presented in Section~\ref{section:conforming} is due to a mismatch  on the boundary conditions imposed in~\eqref{eq:linearBC}. Indeed, in the vicinity of $\partial K$, $u^\eps$ oscillates (as it does everywhere in $\Omega$), whereas the basis functions $\phi^\eps_i$ (and thus $u^\eps_H$) do not, as a consequence of the non-oscillatory boundary conditions $\phi^\eps_i = \phi^0_i$ on $\partial K$. More precisely, there exists a boundary layer around $\partial K$ on which the small scale information (encoded in $\Aa^\eps$) is not properly present in $\phi^\eps_i$, since this information is overshadowed by the incorrect boundary conditions.

\begin{figure}[htbp]
\begin{center}
  \includegraphics[width=40mm]{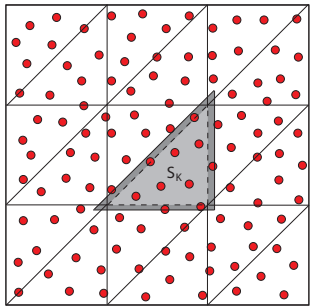}
\end{center}
\caption{Illustration of the oversampling technique: a domain $S_K$, which is slightly larger than $K$, is introduced to solve the local fine scale problems.}
\label{fig:maillages2}
\end{figure}

Since the boundary layer is thin, a possible strategy is to impose (incorrect) Dirichlet boundary conditions not on $\partial K$, but further away, on the boundary of a domain $S_K$ which is slightly larger than $K$, and to only use the interior information, on $K$, to construct the basis functions (see Figure~\ref{fig:maillages2}). This is the well-known oversampling technique. Any reasonable boundary conditions can be imposed on the boundary $\partial S_K$ of the sampling domains $S_K$. In practice (and this is our choice here), affine Dirichlet boundary conditions are often used. The approach thus consists in first solving the local problems
\begin{equation}\label{eq:localpbmsfem2}
\nab \cdot \big[\Aa^\eps\nab \psi^\eps_i\big]=0 \quad \text{in $S_K$}, \qquad \psi^\eps_i \ \ \text{is affine on $\partial S_K$}, \qquad \psi^\eps_i(\bs_j)=\delta_{ij},
\end{equation}
where $S_K$ is a domain which is homothetic to the mesh element $K$ and $\bs_j$ are the coordinates of the vertices of $S_K$. Second, the MsFEM basis functions are defined by
\begin{equation}\label{eq:localpbmsfem2_comp}
\phi^\eps_i = \psi^\eps_{i|K} \quad \text{in $K$}.
\end{equation}
This approach is more costly than the one presented in Section~\ref{section:conforming} since the local problems are set on larger domains.

Note that a unique basis function $\phi^\eps_i$ is associated to each mesh node $i$. However, the basis function $\phi^\eps_i$ may be discontinuous on the element boundaries, i.e. across $\partial K$. This makes the analysis of the oversampling method somewhat more technical. It is interesting and important to notice that, here again, the multiscale basis functions $\phi^\eps_i$ verify a partition of unity property. Indeed, from~\eqref{eq:localpbmsfem2}, we get that
\begin{equation*}
\nab\cdot \left[ \Aa^\eps \nab \left( \sum_i\psi^\eps_i \right) \right]=0 \quad \text{in $S_K$}, \qquad \sum_i \psi^\eps_i = 1 \quad \text{on $\partial S_K$}, 
\end{equation*}
which implies that $\dis \sum_i \psi^\eps_i= 1$ on $S_K$ and thus, in view of~\eqref{eq:localpbmsfem2_comp}, the property~\eqref{eq:pum}. Again, this property will be used in Section~\ref{equilitractionsnonconf} below to derive the \textit{a posteriori} error estimate. 

\begin{remark}
Rather than using~\eqref{eq:localpbmsfem2_comp}, it is also possible to reconstruct $\phi^\eps_i$ on $K$ as a linear combination of the functions $\psi^\eps_j$ by setting 
\begin{equation*}
\phi^\eps_i = \sum_j c_{ij} \ \psi^\eps_j \quad \text{in $K$}, 
\end{equation*}
where $c_{ij}$ are constants (independent of $\eps$) determined by the condition $\dis \sum_j c_{ij} \ \psi^0_j(\bx_k)=\delta_{ik}$, where $\bx_k$ are the coordinates of vertices of $K$ and $\psi_j^0$ are affine functions on $S_K$ such that $\psi^0_j(\bs_k)=\delta_{jk}$. We note that $c_{ij} = \left( \mathbb{G}^{-1} \right)_{ij}$ where the matrix $\mathbb{G}$ is defined by $\left( \mathbb{G} \right)_{ij} = \psi^0_i(\bx_j)$.

In this case, the multiscale basis functions $\phi^\eps_i$ still verify the partition of unity property~\eqref{eq:pum}. Indeed, we note that $\dis \phi_i^0 = \sum_j c_{ij} \ \psi_j^0$ in $K$ (both quantities are affine functions and they are equal on the vertices of $K$). We deduce that
\begin{equation*}
1 = \sum_i \phi_i^0 = \sum_j \left( \sum_i c_{ij} \right) \psi_j^0 \ \text{in $K$} \quad \Longrightarrow \quad \sum_j \left( \sum_i c_{ij} \right) \psi_j^0 = 1 \ \text{in $S_K$}.
\end{equation*}
By evaluation at the vertices of $S_K$, we deduce that, for any $j$, we have $\dis \sum_i c_{ij} = 1$. This eventually implies that $\dis \sum_i \phi^\eps_i = \sum_j \psi^\eps_j = 1$ in $K$.
\end{remark}

Once the multiscale basis functions $\phi^\eps_i$ have been computed, a Galerkin approximation of~\eqref{modelpbweak} on the finite dimensional space $V^\eps_H = \text{Span} \{\phi^\eps_i, \ 1 \leq i \leq I \}$ is performed. Note however that $V^\eps_H \not\subset \mH^1(\Omega)$. The global problem is then defined as
\begin{equation}
\label{eq:GalerkinMsFEM-Over}  
\text{Find $u_H^\eps\in V_H^\eps$ such that, for any $v\in V^\eps_H$,} \quad B^\eps_H(u_H^\eps,v)=F(v), 
\end{equation}
with
\begin{equation}
\label{eq:defBH}  
B^\eps_H(u,v) = \sum_{K\in \mT_H} \int_K [\Aa^\eps\nab u] \cdot \nab v =\intO [\Aa^\eps\nab_H u] \cdot \nab_H v,
\end{equation}
where we have set $\dis \nab_H u = \sum_{K\in \mT_H} 1_K \left. \nab u \right|_K$ and likewise for $v$. 

\medskip

In~\cite{EFE00}, an \textit{a priori} error estimate is derived for the oversampling variant of the MsFEM approach. It is based on general results (see~\cite{CIA78,STR77}) for non-conforming approximations. In the case $H \gg \eps$ and in the periodic setting, assuming that the distance between $K$ and $\partial S_K$ is of order $H$, the estimate reads
\begin{equation}\label{eq:boundover90}
\|u^\eps-u^\eps_H\|_{H,1} \le C_1 H + C_2 \frac{\eps}{H} + C_3 \sqrt{\eps},
\end{equation}
where $C_1$, $C_2$ and $C_3$ are independent of $\eps$ and $H$, and $C_2$ depends on the oversampling ratio. Since the approach is non-conforming, the error is measured in the broken $\mH^1$-norm $\dis \| \cdot \|_{H,1} = \left( \sum_{K \in \mT_H} \|  \cdot \|_{\mH^1(K)}^2\right)^{1/2}$.

The estimate~\eqref{eq:boundover90} is better than the one without oversampling (namely~\eqref{eq:aprioriboundMsFEM}), since the leading order term $\sqrt{\eps/H}$ (due to boundary effects on $\dis \cup_{K \in \mT_H} \partial K$) is now replaced by the term $\eps/H$. However, a resonance effect still exists (cell resonance error), caused by the mismatch between the mesh size and the perfect sample size. It appears as a higher order correction in the error expansion.

\medskip

The oversampling technique by itself does not remove the cell resonance error. It is claimed in~\cite{HOU04} that this error can be completely eliminated using a Petrov-Galerkin formulation. In this formulation, which may simplify implementation aspects, non-conforming multiscale basis functions are used for trial functions alone, while conforming linear functions are considered for test functions. In other words, denoting $V_H^\eps=\text{Span} \{\phi_i^\eps, \ 1 \leq i \leq I \}$ the non-conforming MsFEM space, the Petrov-Galerkin MsFEM formulation reads:
\begin{equation}\label{eq:Petrovformul}
\text{Find $u_H^\eps\in V_H^\eps$ such that, for any $v\in V^0_H$,}\quad B^\eps_H(u_H^\eps,v)=F(v).
\end{equation}

\begin{remark}
Another oversampling-like approach is introduced in~\cite{HEN13}, somewhat in the spirit of Section~\ref{section:recastMsFEM} below. It is claimed in that work that this approach is free from any resonance effect. 
\end{remark}

\subsection{Construction of higher-order MsFEM basis functions}\label{section:higherorder}

Originally, MsFEM was proposed for linear elements, i.e. linear coarse-scale FE shape functions $\phi_i^0$. For many applications, high-order elements are known to be advantageous in terms of accuracy and efficiency. In that spirit, the construction of higher-order oscillating MsFEM basis functions was performed in~\cite{ALL06} using a Taylor expansion and with the introduction of a composition rule (change of variables with the use of local Kozlov harmonic coordinates), allowing for a simple treatment of high-order MsFEM.

For each element $K$ of the coarse mesh $\mT_H$, the functions $w_j^{\eps,K}$ ($j=1,\dots,d$) are defined as the solutions to
\begin{equation}\label{eq:pbsolw}
\nab \cdot \big[\Aa^\eps\nab w_j^{\eps,K}\big]= 0 \ \ \text{in $K$}, \qquad 
w_j^{\eps,K} = \bx \cdot \be_j = x_j \ \ \text{on $\partial K$}.
\end{equation}
Conforming functions $w_j^\eps\in \mH^1(\Omega)$ are next defined by $w_{j|K}^\eps = w_j^{\eps,K}$ on each $K$. We denote $\bw^\eps=(w_1^\eps,\dots,w_d^\eps)^T$. 

\medskip
 
Introducing the set $\{\phi_i^0\}$ of the $P_k$ Lagrange FE basis functions, the local higher-order MsFEM basis $\{\phi^\eps_i\}$ is constructed as 
\begin{equation}\label{eq:compoAllaire}
\phi_{i|K}^\eps=\phi_{i|K}^0 \circ \bw^{\eps,K} \qquad \text{on each $K$}.
\end{equation}
Since $\bw^{\eps,K}(\bx)=\bx$ on $\partial K$, we have $\phi_{i|K}^\eps(\bx)=\phi_{i|K}^0(\bx)$ on $\partial K$. Consequently, a conforming higher-order MsFEM with space $V^\eps_H\subset \mH_{0,\Gamma_D}^1(\Omega)$ is obtained. The Galerkin approximation of~\eqref{modelpbweak} on $V^\eps_H$ yields an approximation $u^\eps_H = u^0_H \circ \bw^\eps \in V^\eps_H$ of $u^\eps$, and $u^0_H$ is expected to be an approximation of the homogenized solution $u^0$.

This high-order approach is easy to implement since the computation of the oscillating functions $w_j^\eps$ is independent of the order $k$ of the coarse mesh basis. 

In the case of a periodically oscillating coefficient $\Aa^\eps = \Aa_{\rm per}(\cdot/\eps)$, \textit{a priori} error estimates are available for this higher-order MsFEM variant. Under some regularity assumptions on $\Aa_{\rm per}$ and on the homogenized solution, it is shown in~\cite{ALL06} that there exists a constant $C$ independent of $\eps$ and $H$ such that
\begin{equation*}
\|u^\eps-u^\eps_H\|_1 \le C\left(H^k+\sqrt{\eps/H}+\sqrt{\eps}\right).
\end{equation*}

\begin{remark}
When using $P_1$ Lagrange FE basis functions, we recover the classical MsFEM given in~\cite{HOU97} and presented in Section~\ref{section:conforming}. Indeed, when $\phi_{i|K}^0 \in P_1(K)$, we infer from~\eqref{eq:compoAllaire} that oscillating basis functions can be written (using that $\bx \mapsto \phi_{i|K}^0(\bx)$ is affine) as
\begin{equation*}
\phi_{i|K}^\eps(\bx)=\phi_{i|K}^0(\bx)+\sum_{j=1}^d\left(w_j^{\eps,K}(\bx)-x_j\right)\frac{\partial \phi_{i|K}^0}{\partial x_j}(\bx).
\end{equation*}
A simple calculation shows that
\begin{equation*}
\nab \cdot[\Aa^\eps\nab \phi_{i|K}^\eps]=\sum_{j=1}^d \nab \cdot[\Aa^\eps\nab w_j^{\eps,K}]\frac{\partial \phi_{i|K}^0}{\partial x_j} =0 \ \ \text{in $K$}, \qquad \phi_{i|K}^\eps = \phi_{i|K}^0 \ \ \text{on $\partial K$},
\end{equation*}
and we thus recover~\eqref{eq:linearBC}.
\end{remark}

\begin{remark}
The domain of definition of $w_j^{\eps,K}$ is the element $K$. However, it is possible to define it on a larger domain $S_K \supset K$, and to only use its restriction on $K$ (oversampling technique). In that case, the support of $\phi_i^\eps$ may be different from that of coarse-scale Lagrange FE basis functions, and its nodal values may also be different. 
\end{remark}

\begin{remark}
Inspired by the framework proposed in~\cite{ALL06}, a more explicit multiscale FE space was constructed for high-order MsFEM in~\cite{HES14}. It does not use a composition rule any more. Local oscillating functions $\widetilde{w}^\eps$ are defined as the solutions to a variant of~\eqref{eq:pbsolw} posed in a domain $S_K \supset K$. These local problems are again solved using a high-order FEM. The multiscale FE basis functions are then defined by
\begin{equation*}
\phi_i^\eps(\bx) = \phi^0_i(\bx)+\left(\widetilde{\bw}^\eps(\bx)-\bx\right)\cdot \nab \phi^0_i(\bx),
\end{equation*}
where both $\phi^0_i\in \mH^1(\Omega)$ and $\widetilde{\bw}^\eps$ are higher-order functions. Note that $\widetilde{\bw}^\eps$ does not belong to $[\mH^1(\Omega)]^d$ in general. 

We note that, due to the use of oversampling, the MsFEM approach of~\cite{HES14} is non-conforming. The convergence result established in~\cite{HES14} is similar to that given in~\cite{EFE00} for the oversampling approach described in Section~\ref{section:oversampling}. On the other hand, if oversampling is not used, results similar to those given in~\cite{ALL06} are recovered. 
\end{remark}

\section{\textit{A posteriori} error estimation using CRE}\label{section:CRE}

In this section, we derive a guaranteed upper bound on the error in the energy norm $\vertiii{u^\eps-u^\eps_H}_H$ with $\vertiii{\cdot}_H = \sqrt{B^\eps_H(\cdot,\cdot)}$, where $B^\eps_H(\cdot,\cdot)$ is defined by~\eqref{eq:defBH}. The method is based on the Constitutive Relation Error (CRE) concept which has been used for decades as an effective and robust error estimation tool in FEM simulations~\cite{LAD83,DES99,LAD04}. Even though we consider this concept for the linear elliptic model of interest, it can be successfully extended and applied to a large class of mechanical models, both in linear (elasticity, visco-elasticity, transient thermal, \dots) and nonlinear (visco-plasticity, contact, \dots) contexts. A review can be found in~\cite{LAD04}. 

\subsection{Basics on CRE}

The CRE concept is based on a dual formulation of the problem~\eqref{modelpbweak}, that we briefly review here. We introduce the space of equilibrated fluxes
\begin{equation*}
W=\{\bp \in \mH(div,\Omega), \quad \nab\cdot \bp + f=0 \ \text{in $\Omega$}, \quad \bp \cdot \bn = g \;\text{on $\Gamma_N$} \},
\end{equation*}
with $\mH(div,\Omega):=\{\bp\in [L^2(\Omega)]^d, \ \ \nab\cdot \bp \in L^2(\Omega)\}$. Any flux field $\widehat{\bp} \in W$ is said statically admissible (SA) and verifies the following relation (weak form of the equilibrium):
\begin{equation}\label{eq:weakformequil}
\forall v \in V, \quad \intO \widehat{\bp} \cdot\nab v  = \intO fv + \int_{\Gamma_N} g \, v.
\end{equation}
Then, for any given approximation $\widehat{u}^\eps$ of $u^\eps$ which is kinematically admissible (KA), i.e. which verifies $\widehat{u}^\eps \in V$, and for any SA flux field $\widehat{\bp}$, we define the CRE functional $E_{CRE}$ as
\begin{multline}
  \label{eq:inde}
  \big( E_{CRE}\left( \widehat{u}^\eps,\widehat{\bp} \right) \big)^2
  =
  \intO (\Aa^\eps)^{-1}(\widehat{\bp}-\Aa^\eps \nab \widehat{u}^\eps)\cdot(\widehat{\bp}-\Aa^\eps \nab \widehat{u}^\eps)
  \\=
  \vertiii{\widehat{\bp}-\Aa^\eps \nab \widehat{u}^\eps}^2_{\mF}
  =
  2\big(J_1(\widehat{u}^\eps)+J_2(\widehat{\bp})\big),
\end{multline}
where $\dis \vertiii{\bp}_{\mF}=\left(\intO (\Aa^\eps)^{-1} \bp \cdot \bp \right)^{1/2}$ is the energy norm for flux fields, and $J_1$ (resp. $J_2$) is the potential (resp. complementary) energy functional defined by
\begin{equation*}
J_1(v) = \frac{1}{2} \intO \Aa^\eps \nab v \cdot \nab v - \intO fv - \int_{\Gamma_N} g \, v, \qquad J_2(\bp) = \frac{1}{2} \intO (\Aa^\eps)^{-1} \bp \cdot \bp.
\end{equation*}
Introduce the exact flux field $\bq^\eps=\Aa^\eps\nab u^\eps$, which is the unique minimizer of the complementary problem
$$
\inf \left\{ J_2(\widehat{\bp}), \ \ \widehat{\bp}\in W \right\}. 
$$
It is then easy to show that
\begin{equation}\label{eq:boundE}
\forall \widehat{\bp} \in W, \qquad \vertiii{u^\eps-\widehat{u}^\eps} = E_{CRE}(\widehat{u}^\eps,\bq^\eps) \le E_{CRE}(\widehat{u}^\eps,\widehat{\bp}),
\end{equation}
so that a guaranteed upper bound on the error $\vertiii{u^\eps-\widehat{u}^\eps}$ is obtained from any SA flux field $\widehat{\bp}$. Its quality is ensured by a suitable choice of $\widehat{\bp} \in W$.

Another version of~\eqref{eq:boundE} is the following Prager-Synge equality:
\begin{equation}\label{eq:Prager}
\forall \widehat{\bp} \in W, \qquad \vertiii{u^\eps-\widehat{u}^\eps}^2+\vertiii{\bq^\eps-\widehat{\bp}}^2_{\mF} = E^2_{CRE}(\widehat{u}^\eps,\widehat{\bp}),
\end{equation}
that holds due to the fact that $\dis \intO (\bq^\eps-\widehat{\bp})\cdot \nab (u^\eps-\widehat{u}^\eps) =0$. 

\begin{remark}
An additional property of the CRE functional, directly derived from~\eqref{eq:Prager}, reads as follows (see~\cite{LAD04}):
\begin{equation}\label{eq:Prager2}
\forall \widehat{\bp} \in W, \qquad 4\vertiii{\bq^\eps-\widehat{\bp}^{\ast}}^2_{\mF} = E^2_{CRE}(\widehat{u}^\eps,\widehat{\bp}) \ \ \text{where} \ \ \widehat{\bp}^{\ast}:=\frac{1}{2}(\widehat{\bp}+\Aa^\eps\nab \widehat{u}^\eps).
\end{equation}
\end{remark}

\begin{remark}
In the case when $\widehat{\bp}$ is not in $W$ but in the larger space $W'=\{\bp \in \mH(div,\Omega), \quad \bp\cdot \bn = g \ \text{on $\Gamma_N$} \} \supset W$, the result~\eqref{eq:boundE} should be changed in
\begin{equation}\label{eq:bound2}
\forall \widehat{\bp} \in W', \quad \vertiii{u^\eps-\widehat{u}^\eps} \le E_{CRE}(\widehat{u}^\eps,\widehat{\bp}) + \frac{C_{\Omega}}{\sqrt{\alpha}} \|\nab\cdot \widehat{\bp} + f \|_0,
\end{equation}
where $C_{\Omega}$ is the Poincar\'e constant of $\Omega$ and $\alpha$ is the ellipticity constant of $\Aa^\eps$ (see~\eqref{eq:elliptic}). Indeed, for any flux $\widehat{\bp} \in \mH(div,\Omega)$, we have
\begin{equation*}
\begin{aligned}
\vertiii{u^\eps-\widehat{u}^\eps}^2 &= \intO(\bq^\eps - \Aa^\eps \nab \widehat{u}^\eps)\cdot \nab(u^\eps-\widehat{u}^\eps) \\
&= \intO(\widehat{\bp}-\Aa^\eps \nab \widehat{u}^\eps)\cdot \nab(u^\eps-\widehat{u}^\eps) + \intO(f+\nab\cdot \widehat{\bp})(u^\eps-\widehat{u}^\eps)+\int_{\Gamma_N}(g-\widehat{\bp} \cdot \bn)(u^\eps-\widehat{u}^\eps). 
\end{aligned}
\end{equation*}
\end{remark}

\begin{remark}
As pointed out previously, the CRE concept can be extended to more complex problems (visco-plasticity, damage, contact, \dots) using duality arguments. The CRE functional reads in the general case
\begin{equation*}
E_{CRE}^2(X,Y)=\intO\big(\Phi(X)+\Phi^{\ast}(Y)-\langle X, Y\rangle\big),
\end{equation*}
where $\Phi$ and $\Phi^{\ast}$ are dual convex potentials (in the Legendre-Fenchel sense) describing the material behavior, and $\langle X, Y\rangle$ is the duality bracket~\cite{LAD04}. This CRE functional should be integrated over the time domain when considering time-dependent problems.
\end{remark}

When using a conforming MsFEM, we can choose $\widehat{u}^\eps=u^\eps_H$ so that the CRE concept directly provides an upper bound on the error $\vertiii{u^\eps-u^\eps_H}$. On the other hand, in the case a non-conforming MsFEM version is used (such as the one described in Section~\ref{section:oversampling}), one cannot choose $\widehat{u}^\eps=u^\eps_H$ since $u^\eps_H$ is not kinematically admissible ($u^\eps_H \notin \mH^1(\Omega)$). The procedure should be amended as follows. The error is decomposed into its conforming and non-conforming parts: we write
$$
u^\eps-u^\eps_H = e_C + e_{NC},
$$
where $e_C \in V$ is defined by
\begin{equation*}
\forall v \in V, \quad B^\eps(e_C,v)=B_H^\eps(u^\eps-u^\eps_H,v)=F(v)-\intO[\Aa^\eps \nab_Hu^\eps_H] \cdot \nab v.
\end{equation*}
We then check that
$$
\forall v \in V, \quad B_H^\eps(e_{NC},v) = B_H^\eps(u^\eps-u^\eps_H - e_C,v) = B_H^\eps(u^\eps-u^\eps_H,v) - B_H^\eps(e_C,v) = 0.
$$
We thus have that
\begin{equation}\label{eq:errorsplitting}
\vertiii{u^\eps-u^\eps_H}^2_H = \vertiii{e_C + e_{NC}}^2_H = \vertiii{e_C}^2 + \vertiii{e_{NC}}^2_H.
\end{equation}
We have $u^\eps- e_C - e_{NC} = u^\eps_H$, and hence $u^\eps- e_C$ is the orthogonal projection of $u^\eps_H$ on the conforming space $V$, and $\dis \vertiii{e_{NC}}_H = \min_{v \in V}\vertiii{v-u^\eps_H}_H$. The idea is then to approximate the error splitting~\eqref{eq:errorsplitting} by building a conforming (kinematically admissible) function $\widehat{u}^\eps_H$, close to $u^\eps_H$. This is performed by a local nodal averaging on the edges $\Gamma$ of $\mT_H$, at each node of the fine mesh used to compute the multiscale basis functions. Consequently, we get
\begin{equation*}
\vertiii{u^\eps-u^\eps_H}_H \le \vertiii{u^\eps-\widehat{u}^\eps_H} + \vertiii{\widehat{u}^\eps_H-u^\eps_H}_H.
\end{equation*}
The approximation $\vertiii{\widehat{u}^\eps_H-u^\eps_H}_H$ of the non-conforming part of the error is computable, whereas the approximation $\vertiii{u^\eps-\widehat{u}^\eps_H}$ of the conforming part of the error can be estimated using the CRE functional:
\begin{equation*}
\forall \widehat{\bp} \in W, \quad \vertiii{u^\eps-\widehat{u}^\eps_H} \le E_{CRE}(\widehat{u}^\eps_H,\widehat{\bp}).
\end{equation*} 

\subsection{Construction of an equilibrated flux field}\label{section:constructSA}

The technical point in the CRE concept is the construction of a relevant admissible flux $\widehat{\bp} \in W$, that we now discuss. 

\subsubsection{General procedure}

The most accurate strategy would consist in using a Galerkin approximation of the dual formulation of the problem, arising from the minimization (see~\eqref{eq:inde} and~\eqref{eq:boundE}) of the complementary energy $J_2$ (FEM with equilibrated flux basis functions, see~\cite{FRA65}). Nevertheless, this strategy is not suited to current simulation softwares (building finite dimensional subspaces of $W$ is challenging in practice) and requires the costly solution to an additional global problem. This may be not reasonable for error estimation purposes.

We use an alternative method, referred as the Hybrid-Flux Technique, or Element Equilibration Technique (EET) in the recent literature (see~\cite{PLE11}). This technique enables the construction of an equilibrated flux field verifying~\eqref{eq:weakformequil}, denoted $\widehat{\bq}^\eps_H$ in the following, from a post-processing of the MsFEM flux $\bq^\eps_H = \Aa^\eps \nab u^\eps_H$ at hand and from the solutions to independent elementary problems. This method, which is an extension to MsFEM of the one exposed in~\cite{LAD83,LAD96,LAD04} for FEM simulations, is made of two steps:
\begin{itemize}
\item[$\bullet$] \textbf{Step 1:} construction of tractions $\widehat{g}_K$ along the edges of each element $K \in \mT_H$, with $\widehat{g}_K=g$ on $\Gamma_N$. These tractions should verify the equilibrium at the element level:
\begin{equation}\label{eq:equilelement}
\forall K, \quad \int_K f  + \int_{\partial K}\widehat{g}_K =0.
\end{equation}
We will see in what follows that~\eqref{eq:equilelement} indeed holds.

In order to satisfy the continuity of normal fluxes along element edges, tractions are defined as $\widehat{g}_{K|\Gamma} = \eta_K^{\Gamma} \, \widehat{g}_{|\Gamma}$ over each edge $\Gamma$ of $\partial K$, with $\eta_K^{\Gamma}=\pm 1$ and where $\widehat{g}_{|\Gamma}$ has a well-defined value on $\Gamma$. After ordering the elements of the mesh $\mT_H$, a possible setting for $\eta_K^{\Gamma}$ is the following:
\begin{equation}\label{eq:signeta}
\eta_K^{\Gamma}=\left\{
\begin{array}{l}
+1\quad \text{if}\; \Gamma = \partial K \cap\partial K', \quad K>K', \\ 
-1 \quad \text{if}\; \Gamma = \partial K \cap\partial K', \quad K<K', \\
+1 \quad \text{if}\; \Gamma \in \partial \Omega.
\end{array}
\right. 
\end{equation}
\item[$\bullet$] \textbf{Step 2:} local construction of $\widehat{\bq}^\eps_{H|K}$, over each element $K \in \mT_H$, such that
\begin{equation}\label{eq:localpbelement}
-\nab \cdot \widehat{\bq}^\eps_{H|K} = f \ \ \text{in $K$}, \qquad \widehat{\bq}^\eps_{H|K}\cdot \bn_K = \widehat{g}_K \ \ \text{on $\partial K$},
\end{equation}
where $\widehat{g}_K$ are the tractions built in Step 1, and $\bn_K$ is the outgoing normal on $\partial K$. 
\end{itemize}



\medskip

Details on these two steps are given in the following sections. In Section~\ref{equilitractionsconf}, we discuss the construction of equilibrated tractions (Step 1 of the above procedure), in the case of the conforming MsFEM approach of Section~\ref{section:conforming} and in the case of the Petrov-Galerkin approach~\eqref{eq:Petrovformul}. In Section~\ref{equilitractionsnonconf}, we again discuss the construction of equilibrated tractions, but now in the case of the non-conforming MsFEM approach~\eqref{eq:GalerkinMsFEM-Over} of Section~\ref{section:oversampling}. In Section~\ref{elementaryproblems}, we discuss the construction of equilibrated fluxes (Step 2 of the above procedure). 

We note that the hybrid construction method we use is based on fine-scale computations (Step 2) only at the element level (and not on patches of elements, as for alternative reconstruction methods). It is therefore particularly suited to the MsFEM framework. Furthermore, we show below that it is actually possible to perform most of the computations of Step 2 in the \textit{offline} phase of MsFEM, independently of the loadings $f$ and $g$.

\subsubsection{Computation of equilibrated tractions for conforming or Petrov-Galerkin non-conforming MsFEM formulations}\label{equilitractionsconf}

We assume here that the numerical solution $u_H^\eps$ is either provided by the conforming MsFEM approach of Section~\ref{section:conforming} or by the Petrov-Galerkin approach~\eqref{eq:Petrovformul}.

An arbitrary but convenient manner to derive equilibrated tractions $\widehat{g}_K$ by post-processing $\bq^\eps_H$ is the enforcement of the so-called \textit{prolongation condition} (energy condition). We require that, for each element $K$ and each node $i$ connected to $K$,
\begin{equation}\label{eq:prolongation-pre}
\int_K (\widehat{\bq}_H^\eps-\bq_H^\eps) \cdot \nab \phi_i^\star =0,
\end{equation}
where $\phi_i^\star$ are the test functions of the MsFEM formulation, i.e. $\phi_i^\star=\phi_i^\eps$ for conforming MsFEM and $\phi_i^\star=\phi_i^0$ for Petrov-Galerkin (possibly non-conforming) MsFEM. Since $\widehat{\bq}_H^\eps \in W$, this is equivalent to requiring that 
\begin{equation}\label{eq:prolongation}
\int_{\partial K} \widehat{g}_K \, \phi_i^\star = \int_K (\bq^\eps_H \cdot \nab \phi^\star_i - f \phi^\star_i) =: Q_i^K,
\end{equation}
where $\widehat{g}_K = \widehat{\bq}_H^\eps \cdot \bn_K$ on $\partial K$. Note that the quantity $Q_i^K$ is fully computable. Enforcing~\eqref{eq:prolongation} naturally provides for the local equilibrium condition~\eqref{eq:equilelement} due to the partition of unity property $\sum_i\phi_{i|K}^\star=1$ (see~\eqref{eq:pum}).

\begin{figure}[htbp]
\begin{center}
\includegraphics[width=110mm]{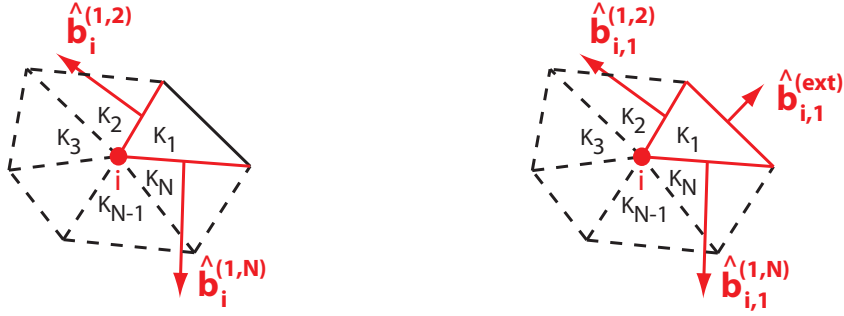}
\end{center}
\caption{Illustration of the prolongation condition applied to elements connected to a node $i$ in a 2D mesh: conforming case (left), non-conforming case (right).}
\label{fig:system12}
\end{figure}

Collecting the relations~\eqref{eq:prolongation} for all the elements connected to a given node $i$ yields a local system where the unknowns are the projections $\dis \int_{\partial K} \widehat{g}_K \, \phi_i^\star$ of the tractions $\widehat{g}_K$ on the basis function $\phi_i^\star$. The local system spreads over the support of $\phi_i^\star$, i.e. the patch $\Omega_i$ of elements connected to node $i$. As an example, consider an internal vertex node $i$ connected to $N$ elements in a 2D mesh $\mT_H$ made of 3-nodes triangle elements (see left side of Figure~\ref{fig:system12}). Since $\phi_{i|\partial \Omega_i}^\star=0$ and $\phi_i^\star$ is continuous across $\Gamma_{jk} = \partial K_j \cap\partial K_k$ for the two MsFEM formulations considered in this section, the local system reads
\begin{equation}\label{eq:systemi}
\begin{aligned}
\widehat{b}_i^{(1,2)}-\widehat{b}_i^{(N,1)} &= Q_i^{K_1}, \\
\widehat{b}_i^{(2,3)}-\widehat{b}_i^{(1,2)} &= Q_i^{K_2}, \\
\dots &= \dots \\
\widehat{b}_i^{(N-1,N)}-\widehat{b}_i^{(N-2,N-1)} &= Q_i^{K_{N-1}}, \\
\widehat{b}_i^{(N,1)}-\widehat{b}_i^{(N-1,N)} &= Q_i^{K_N}, 
\end{aligned}
\end{equation}
with $\dis \widehat{b}_i^{(j,k)} = \int_{\Gamma_{jk}} \eta_{K_j}^{\Gamma_{jk}} \ \widehat{g}_{|\Gamma_{jk}} \, \phi_i^\star$. There are $N$ equations and $N$ unknowns. The $N \times N$ matrix corresponding to the linear system~\eqref{eq:systemi} is however not invertible. Solutions to~\eqref{eq:systemi} yet exist since the right-hand side of~\eqref{eq:systemi} satisfies the appropriate compatibility condition:
\begin{equation*}
\sum_{j=1}^NQ_i^{K_j}=\int_{\Omega_i} (\bq^\eps_H \cdot \nab \phi^\star_i - f\phi^\star_i)=\int_{\Omega} (\bq^\eps_H \cdot \nab \phi^\star_i - f\phi^\star_i)=0,
\end{equation*}
where the last equality above follows from the fact that the MsFEM flux $\bq^\eps_H$ verifies the weak equilibrium~\eqref{eq:GalerkinMsFEM} (in the case of conforming MsFEM) or~\eqref{eq:Petrovformul} (in the case of Petrov-Galerkin MsFEM). Recall that we consider here for simplicity an internal vertex node, so that the term $\dis \int_{\Gamma_N} g \, \phi^\star_i$ in $F(\phi^\star_i)$ vanishes.

A unique solution to the linear system~\eqref{eq:systemi} may be obtained by minimizing a given cost function. In practice, the following least squares cost function is chosen:
\begin{equation}\label{eq:costfunction}
\sum_{j,k, \ \text{$K_j$ and $K_k$ neighbors}}\frac{\Big(\widehat{b}_i^{(j,k)}-\overline{b}_i^{(j,k)}\Big)^2}{|\Gamma_{jk}|^2} \quad \text{with} \quad \overline{b}_i^{(j,k)} = \frac{1}{2} \int_{\Gamma_{jk}} \phi^\star_i \, (\bq^\eps_{H|K_j}+\bq^\eps_{H|K_k})\cdot \bn_{K_j}.
\end{equation}

\begin{remark}
The previous constrained minimization problem can essentially be analytically solved. We set
\begin{equation*}
\begin{aligned}
\widehat{\bX} &= \left[\widehat{b}_i^{(1,2)},\widehat{b}_i^{(2,3)},\dots,\widehat{b}_i^{(N,1)}\right]^T, \qquad
\overline{\bX} = \left[\overline{b}_i^{(1,2)},\overline{b}_i^{(2,3)},\dots,\overline{b}_i^{(N,1)}\right]^T, \\
\bX_0 &= \left[0,Q_i^{K_2},Q_i^{K_2}+Q_i^{K_3},\dots,\sum_{j=2}^{N}Q_i^{K_j} \right]^T, \quad
\bA=[1,1,\dots,1]^T.
\end{aligned}
\end{equation*}
Then the solutions to~\eqref{eq:systemi} are the vectors $\widehat{\bX}(s)=s \bA +\bX_0$ for any $s \in \RR$. The cost function~\eqref{eq:costfunction} reads
\begin{equation*}
  f(s)
  =
  (\widehat{\bX}(s)-\overline{\bX})^T\D(\widehat{\bX}(s)-\overline{\bX}) 
  =
  s^2\bA^T\D \bA-2s\bA^T\D(\overline{\bX}-\bX_0)+(\overline{\bX}-\bX_0)^T\D(\overline{\bX}-\bX_0)
\end{equation*}
with $\D=diag(1/|\Gamma_{12}|^2,1/|\Gamma_{23}|^2,\dots,1/|\Gamma_{N1}|^2)$. The minimization of $f$ leads to
$$
s^\star=(\bA^T\D\bA)^{-1}\bA^T\D(\overline{\bX}-\bX_0).
$$
The selected solution $\widehat{\bX}(s^\star)$ to~\eqref{eq:systemi} is hence obtained explicitly with little computational effort. 
\end{remark}

\begin{remark}
For a node $i$ lying on the boundary $\partial \Omega$, the procedure to recover $\widehat{b}_i^{(j,k)}$ uses similar ideas. We refer to~\cite{LAD04} for details.
\end{remark}

\begin{remark}
Some variants of the previously introduced procedure exist in the literature for FEM. We refer to~\cite{FLO02,PLE12} and also to~\cite{LAD10c,PLE11}.
\end{remark}

\begin{figure}[htbp]
\begin{center}
\includegraphics[width=120mm]{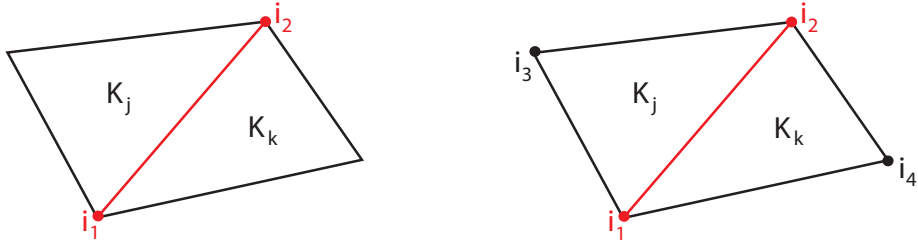}
\end{center}
\caption{Environment of the interface $\Gamma_{jk}$ in the conforming case (left) and the non-conforming case (right).}
\label{fig:interface1}
\end{figure}

Once projections $\widehat{b}_i^{(j,k)}$ have been obtained for all nodes $i$ of the mesh $\mT_H$, tractions $\widehat{g}_{\Gamma_{jk}}$ can be recovered. We choose to search $\widehat{g}_{\Gamma_{jk}}$ as an affine combination of the functions $\phi^\star_{i|\Gamma_{jk}}$. For the 2D example of the left side of Figure~\ref{fig:interface1}, this reads
\begin{equation}\label{eq:construction1}
\widehat{g}_{\Gamma_{jk}} = \eta_{K_j}^{\Gamma_{jk}} \, \langle \bq^\eps_H \rangle_{\Gamma_{jk}}\cdot \bn_{K_j} + \alpha_1 \, \phi^\star_{i_1|\Gamma_{jk}} + \alpha_2 \, \phi^\star_{i_2|\Gamma_{jk}},
\end{equation}
where $\dis \langle \bq^\eps_H \rangle_{\Gamma_{jk}} := \frac{1}{2} (\bq^\eps_{H|K_j} + \bq^\eps_{H|K_k})$,
where the nodes $i_1$ and $i_2$ are the two vertex nodes connected to the edge $\Gamma_{jk}$ (see Figure~\ref{fig:interface1}). Enforcing that $\dis \widehat{b}_i^{(j,k)} = \int_{\Gamma_{jk}} \eta_{K_j}^{\Gamma_{jk}} \ \widehat{g}_{|\Gamma_{jk}} \, \phi_i^\star$ leads to the simple $2\times 2$ linear system
\begin{equation*}
\eta_{K_j}^{\Gamma_{jk}} \, \left[
\begin{array}{c c}
\int_{\Gamma_{jk}}\phi^\star_{i_1}\phi^\star_{i_1} & \int_{\Gamma_{jk}}\phi^\star_{i_1}\phi^\star_{i_2} \\
\int_{\Gamma_{jk}}\phi^\star_{i_1}\phi^\star_{i_2} & \int_{\Gamma_{jk}}\phi^\star_{i_2}\phi^\star_{i_2}
\end{array}
\right]
\left[
\begin{array}{c}
\alpha_1 \\ \alpha_2
\end{array}
\right] = 
\left[
\begin{array}{c}
\widehat{b}_{i_1}^{(j,k)}-\overline{b}_{i_1}^{(j,k)} \\ \widehat{b}_{i_2}^{(j,k)}-\overline{b}_{i_2}^{(j,k)}
\end{array}
\right], 
\end{equation*}
for which the $2\times 2$ matrix $\M^{\Gamma_{jk}}$ can be computed and inverted in the \textit{offline} stage of MsFEM.

\medskip

\begin{remark}
The previous system is invertible, since the Cauchy-Schwarz inequality between two linearly independent functions $\phi^\star_{i_1}$ and $\phi^\star_{i_2}$ leads to a strictly positive determinant for $\M^{\Gamma_{jk}}$. 
\end{remark}

\begin{remark}
A usual procedure used in FEM (see~\cite{LAD04}) is to search $\widehat{g}_{\Gamma_{jk}}$ as a linear combination of the functions $\phi^\star_{i|\Gamma_{jk}}$. For the 2D example of the left side of Figure~\ref{fig:system12}, this corresponds to writing $\widehat{g}_{\Gamma_{jk}} = \beta_1 \, \phi^\star_{i_1|\Gamma_{jk}} + \beta_2 \, \phi^\star_{i_2|\Gamma_{jk}}$ for some coefficients $\beta_1$ and $\beta_2$. However, for the two MsFEM formulations considered here, the test functions $\phi^\star_i$ are affine on the element edges. This construction thus leads to an affine evolution of $\widehat{g}_{\Gamma_{jk}}$ on the edges, and may result in a poor admissible flux field $\widehat{\bq}^\eps_H$ that does not capture small scale oscillations near element edges. For this reason, we prefer to recover $\widehat{g}_{\Gamma_{jk}}$ under the form~\eqref{eq:construction1}.
\end{remark}

\subsubsection{Computation of equilibrated tractions for the oversampling MsFEM formulation}\label{equilitractionsnonconf}

We assume here that $u_H^\eps$ is the solution to the oversampling MsFEM approach~\eqref{eq:GalerkinMsFEM-Over}. This non-conforming approach (described in Section~\ref{section:oversampling}) is not a classical discontinuous Galerkin method since there is only one degree of freedom per node and no additional jump term in the weak formulation~\eqref{eq:GalerkinMsFEM-Over}. Consequently, recovering equilibrated tractions is not as easy as for classical discontinuous Galerkin methods (see e.g.~\cite{DES98b,ERN07}), where test functions with support over a single element can be considered. Nevertheless, when using Galerkin MsFEM formulation with oversampling, the main aspects of the previous strategy to construct equilibrated tractions can be reused. In particular, enforcing the prolongation condition~\eqref{eq:prolongation-pre} (or equivalently~\eqref{eq:prolongation}) with non-conforming test functions $\phi_i^\eps$ still provides tractions $\widehat{g}_K$ satisfying the local equilibrium condition~\eqref{eq:equilelement} due to the partition of unity property~\eqref{eq:pum}. However, some minor changes are required.

Considering again a 2D example, and noticing that $\phi^\eps_{i|\partial \Omega_i} \neq 0$ and $\phi^\eps_{i|K_j}\neq \phi^\eps_{i|K_k}$ on $\Gamma_{jk}$, the prolongation condition provides for the following system (see the right side of Figure~\ref{fig:system12}): 
\begin{equation}\label{eq:systemi2}
\begin{aligned}
\widehat{b}_{i,1}^{(1,2)}-\widehat{b}_{i,1}^{(N,1)} + \widehat{b}_{i,1}^{(ext)} &= Q_i^{K_1}, \\
\widehat{b}_{i,2}^{(2,3)}-\widehat{b}_{i,2}^{(1,2)} + \widehat{b}_{i,2}^{(ext)} &= Q_i^{K_2}, \\
\dots &= \dots \\
\widehat{b}_{i,N-1}^{(N-1,N)}-\widehat{b}_{i,N-1}^{(N-2,N-1)} + \widehat{b}_{i,N-1}^{(ext)} &= Q_i^{K_{N-1}}, \\
\widehat{b}_{i,N}^{(N,1)}-\widehat{b}_{i,N}^{(N-1,N)} + \widehat{b}_{i,N}^{(ext)} &= Q_i^{K_N}, 
\end{aligned}
\end{equation}
with $\dis \widehat{b}_{i,l}^{(j,k)} = \int_{\Gamma_{jk}} \eta_{K_j}^{\Gamma_{jk}} \ \widehat{g}_{|\Gamma_{jk}} \, \phi^\eps_{i|K_l}$, $\dis \widehat{b}_{i,l}^{(ext)} = \int_{\Gamma_{ext}^{i,l}} \eta_{K_l}^{\Gamma_{ext}^{i,l}} \ \widehat{g}_{|\Gamma_{ext}^{i,l}} \, \phi^\eps_{i|K_l}$ and $\Gamma_{ext}^{i,l} = \partial K_l \cap \partial \Omega_i$. There are $N$ equations and $3N$ unknowns in~\eqref{eq:systemi2}.

In order to come back to the former local system~\eqref{eq:systemi}, we additionally enforce the following conditions:
\begin{equation*}
\begin{aligned}
\widehat{b}_{i,j}^{(j,k)} &= \widehat{b}_{i,k}^{(j,k)} =: \widehat{b}_i^{(j,k)} \quad \text{(the projections of $\widehat{g}_{|\Gamma_{jk}}$ over $\phi^\eps_{i|K_j}$ and $\phi^\eps_{i|K_k}$ are equal)}, \\
\widehat{b}_{i,j}^{(ext)} &= 0.
\end{aligned}
\end{equation*}
The last condition enables to keep a local system, i.e. to not couple~\eqref{eq:systemi2} (written in the patch $\Omega_i$) with similar systems written for neighboring patches.

The system~\eqref{eq:systemi2} then becomes identical to~\eqref{eq:systemi}, and we can compute $\widehat{b}_i^{(j,k)}$ following the procedure detailed in Section~\ref{equilitractionsconf}, using again the quantities
\begin{equation*}
\overline{b}_i^{(j,k)} = \frac{1}{2}\int_{\Gamma_{jk}}(\phi^\eps_{i|K_j}\bq^\eps_{H|K_j}+\phi^\eps_{i|K_k}\bq^\eps_{H|K_k})\cdot \bn_{K_j}.
\end{equation*}

\medskip

Next, on each edge $\Gamma_{jk}$, $\widehat{g}_{\Gamma_{jk}}$ is searched as an affine combination of the multiscale functions $\phi^\eps_{i|\Gamma_{jk}}$. For the 2D example of the right side of Figure~\ref{fig:interface1}, this reads
\begin{equation}\label{eq:construction2}
\begin{aligned}
\widehat{g}_{\Gamma_{jk}} = & \eta_{K_j}^{\Gamma_{jk}} \ \langle \bq^\eps_H \rangle_{\Gamma_{jk}}\cdot \bn_{K_j} + \alpha_{1,j} \, \phi^\eps_{i_1|K_j,\Gamma_{jk}} 
+ \alpha_{1,k} \, \phi^\eps_{i_1|K_k,\Gamma_{jk}} \\
&+ \alpha_{2,j} \, \phi^\eps_{i_2|K_j,\Gamma_{jk}} 
+ \alpha_{2,k} \, \phi^\eps_{i_2|K_k,\Gamma_{jk}} 
+ \alpha_{3} \, \phi^\eps_{i_3|K_j,\Gamma_{jk}}
+ \alpha_{4} \, \phi^\eps_{i_4|K_k,\Gamma_{jk}},
\end{aligned}
\end{equation}
where $\dis \langle \bq^\eps_H \rangle_{\Gamma_{jk}}$ is defined below~\eqref{eq:construction1}, where the nodes $i_1$ and $i_2$ are the two vertex nodes connected to the edge $\Gamma_{jk}$, and the nodes $i_3$ and $i_4$ are the other vertex nodes of the elements $K_j$ and $K_k$, respectively (see Figure~\ref{fig:interface1}). The coefficients $(\alpha_{1,j},\alpha_{1,k},\alpha_{2,j},\alpha_{2,k},\alpha_{3},\alpha_{4})$ are obtained from the following relations:
\begin{equation*}
\begin{aligned}
\widehat{b}_{i_1,j}^{(j,k)} &= \widehat{b}_{i_1,k}^{(j,k)} = \widehat{b}_{i_1}^{(j,k)}, \\
\widehat{b}_{i_2,j}^{(j,k)} &= \widehat{b}_{i_2,k}^{(j,k)} = \widehat{b}_{i_2}^{(j,k)}, \\
\widehat{b}_{i_3,j}^{(ext)} &= 0, \\
\widehat{b}_{i_4,k}^{(ext)} &= 0.
\end{aligned}
\end{equation*}
This leads to the inversion of a simple $6\times 6$ linear system. The matrix $\M^{\Gamma_{jk}}$ of that system can again be computed and inverted in the \textit{offline} stage of MsFEM.

\subsubsection{Solution to elementary problems}
\label{elementaryproblems}

For given tractions $\widehat{g}_K$, the optimal admissible flux field $\widehat{\bq}^\eps_{H|K}$ in the element $K$ is the one that minimizes
$$
E^2_{CRE|K}(\widehat{u}^\eps_H,\widehat{\bq})
=
\int_K (\Aa^\eps)^{-1}(\widehat{\bq}-\Aa^\eps \nab \widehat{u}^\eps)\cdot(\widehat{\bq}-\Aa^\eps \nab \widehat{u}^\eps)
=
\vertiii{\widehat{\bq}-\Aa^\eps \nab \widehat{u}^\eps}^2_{\mF|K}
$$
among all flux fields $\widehat{\bq}$ verifying~\eqref{eq:localpbelement}. For such fields, we compute
\begin{equation*}
E^2_{CRE|K}(\widehat{u}^\eps_H,\widehat{\bq}) = \vertiii{\widehat{\bq}-\Aa^\eps\nab\widehat{u}^\eps_H}_{\mF|K}^2 = \vertiii{\widehat{\bq}}^2_{\mF|K} + \vertiii{\Aa^\eps\nab\widehat{u}^\eps_H}_{\mF|K}^2 - 
2\left[\int_K f \, \widehat{u}^\eps_H  + \int_{\partial K}\widehat{g}_K \, \widehat{u}^\eps_H \right],
\end{equation*}
and it is thus equivalent to minimizing $\vertiii{\widehat{\bq}}^2_{\mF|K}=2J_{2|K}(\widehat{\bq})$ among fluxes $\widehat{\bq}$ verifying~\eqref{eq:localpbelement}. Introducing the functional spaces
\begin{equation*}
\begin{aligned}
W^K &:= \{\bp \in \mH(div,K), \quad \nab \cdot \bp + f =0 \ \text{in $K$}, \quad \bp \cdot \bn_K = \widehat{g}_K \ \text{on $\partial K$} \}, \\
W^K_0 &:= \{\bp \in \mH(div,K), \quad \nab \cdot \bp = 0 \ \text{in $K$}, \quad \bp \cdot \bn_K = 0 \ \text{on $\partial K$} \},
\end{aligned}
\end{equation*}
the optimal admissible flux field satisfies
\begin{equation}\label{eq:minJ2}
\widehat{\bq}^\eps_{H|K} \in W^K \quad \text{and} \quad \forall \bp \in W^K_0, \quad \int_K \widehat{\bq}^\eps_{H|K} \cdot (\Aa^\eps)^{-1} \bp = 0.
\end{equation}
We describe now two ways to compute $\widehat{\bq}^\eps_{H|K}$. 

\medskip

As a first case, assume that the load $f$ and the tractions $\widehat{g}_K$ are polynomial. Then $\widehat{\bq}^\eps_{H|K} \in W^K$ may be constructed analytically with polynomial form (in the spirit of~\cite{LAD97}) or numerically solving~\eqref{eq:minJ2} by means of a non-conventional FE method in fluxes (in the spirit of~\cite{FRA65}).

\medskip

Another technique, more appealing for implementation purposes, consists in a dual approach. Indeed, duality shows that the minimization of $J_{2|K}(\widehat{\bq})$ among fluxes $\widehat{\bq}$ satisfying~\eqref{eq:localpbelement} is equivalent to finding a primal field $w^\eps \in \mH^1(K)$ such that
\begin{equation}\label{eq:dualelempb}
\widehat{\bq}^\eps_{H|K} = \Aa^\eps\nab w^\eps \quad \text{and} \quad \forall v \in \mH^1(K), \quad B^\eps_{|K}(w^\eps,v)=\int_K f v  + \int_{\partial K}\widehat{g}_K \, v,
\end{equation}
where $\dis B^\eps_{|K}(u,v) = \int_K \left[ \Aa^\eps \, \nab u \right] \cdot \nab v$. In view of~\eqref{eq:equilelement}, the above local Neumann problem~\eqref{eq:dualelempb} is well-posed (up to the addition of a constant). An accurate approximation $\widetilde{w}^\eps$ of its solution $w^\eps$ is in practice computed using a higher-order basis within the element $K$. In order to reuse the \textit{offline} MsFEM computations at hand and to limit the computational costs, we construct the local higher-order basis using the technique proposed in~\cite{ALL06} and detailed in Section~\ref{section:higherorder}. We use polynomial shape functions $\phi_{i|K}^0$ of degree $1+k$ ($k=3$ in practice) in the composition rule~\eqref{eq:compoAllaire}. Furthermore, when the degree $1+k$ is large, it may be effective to use hierarchical subspaces of the full space $V_{H|K}^{1+k}$, orthogonal with respect to the inner product defined by $B^\eps_{|K}$ (see~\cite{AIN00}). We do not pursue in that direction.

\begin{remark}
When using the numerical approximation $\widetilde{\bq}^\eps_{H|K} = \Aa^\eps \nab \widetilde{w}^\eps$ of $\widehat{\bq}^\eps_{H|K} = \Aa^\eps \nab w^\eps$, the strict upper error bound given by the CRE concept is not valid any more. The extra polynomial degree $k$ is in practice chosen high enough so that an asymptotic regime is reached. We refer to~\cite{BAB94b} for studies of that question.
\end{remark}

In the case when the body loading $f$ is regular in each element $K$, computations associated to the resolution of~\eqref{eq:dualelempb} can mostly be performed \textit{offline}, thus avoiding fine-scale computations in the \textit{online} phase, and thus respecting the MsFEM paradigm in multi-query contexts.

Assume first that $f$ is constant in each coarse element $K$. In view of~\eqref{eq:construction1} and~\eqref{eq:construction2}, we notice that the tractions $\widehat{g}_K$ are defined as linear combinations of the functions $\phi_i^\star$ and $\nab \phi_i^\star \cdot \bn_K$ (where $\phi_i^\star=\phi_i^\eps$ or $\phi_i^\star=\phi_i^0$ depending on the MsFEM formulation). This leads to solving the following set of well-posed Neumann problems, in the \textit{offline} phase of the approach. For each element $K$, each edge $\Gamma_{jk}$ of $K$ and each node $i$ of $K$, we look for $\theta^{jk}_i \in \mH^1(K)$ such that
\begin{equation}\label{eq:dualelempb12_a}
  \forall v \in \mH^1(K), \qquad B^\eps_{|K}(\theta^{jk}_i,v) = \int_{\Gamma_{jk}} \phi_i^\star \, v - \int_K\left(\frac{1}{|K|}\int_{\Gamma_{jk}}\phi_i^\star\right) v
\end{equation}
and for $\mu^{jk}_i \in \mH^1(K)$ such that
\begin{equation}\label{eq:dualelempb12_b}
\forall v \in \mH^1(K), \quad B^\eps_{|K}(\mu^{jk}_i,v) = \int_{\Gamma_{jk}} \nab \phi_i^\star \cdot \bn_K \, v - \int_K\left(\frac{1}{|K|}\int_{\Gamma_{jk}}\nab\phi_i^\star \cdot \bn_K\right)v.
\end{equation}
Of course, $\theta^{jk}_i$ and $\mu^{jk}_i$ are only defined up to the addition of a constant. Since $f$ is constant in $K$, the solution $w^\eps$ to~\eqref{eq:dualelempb} can then be written as a linear combination of $\theta^{jk}_i$ and $\mu^{jk}_i$.

For more complex variations of $f$ in $K$, we write $f = \overline{f}_{|K} + \delta f_{|K}$ in $K$, where $\delta f_{|K}$ is self-equilibrated in $K$ (i.e. $\dis \int_K \delta f_{|K} =0$). A flux $\widehat{\bq}^\eps_{\overline{f}|K}$ which is equilibrated with the tractions $\widehat{g}_K$ (the computation of which is inexpensive, whatever the form of $f$) and the constant, first term $\overline{f}_{|K}$ can be computed using the above approach (replacing $f$ with $\overline{f}_{|K}$ in~\eqref{eq:dualelempb} and using the \textit{offline} procedure~\eqref{eq:dualelempb12_a}--\eqref{eq:dualelempb12_b}). The mean-free term $\delta f_{|K}$ can be handled in two manners. It can either be taken into account in the CRE bound (see~\eqref{eq:bound2}). It can also give rise to a flux $\widehat{\bq}^\eps_{\delta f|K}$, equilibrated with $\delta f_{|K}$, and computed as
\begin{equation}\label{eq:dualelempb2}
\widehat{\bq}^\eps_{\delta f|K}=\Aa^\eps\nab \rho \quad \text{and} \quad \forall v \in \mH^1(K), \quad B^\eps_{|K}(\rho,v) = \int_K \delta f_{|K} \, v.
\end{equation}
It is of course possible to expand $\delta f_{|K}$ over a given polynomial basis made of self-equilibrated basis functions $r_{j|K}$. Then most of the computations associated to~\eqref{eq:dualelempb2} can be performed \textit{offline}, by computing beforehand the elementary solutions $\rho_j$ to
\begin{equation}\label{eq:dualelempb22}
\forall v \in \mH^1(K), \quad B^\eps_{|K}(\rho_j,v) = \int_K r_{j|K} \, v.
\end{equation}
Recall now that $f$ is a macroscopic function (it does not depend on $\eps$). On each coarse element $K$, $\delta f_{|K}$ may thus be accurately approximated with a limited number of basis functions $r_{j|K}$.
 
\section{Adaptive strategy}\label{section:adaptive}

In the previous section, we have developed a CRE estimate that yields a fully computable upper bound quantifying the overall MsFEM error. In this section, we develop error indicators that assess the various error sources and enable an automatic adaptation of the MsFEM parameters (the coarse mesh size $H$, the fine mesh size $h$ and the oversampling ratio) in order to reach a given error tolerance while keeping the computational cost minimal. 
The error indicators are defined from (hierarchical) auxiliary reference problems, which are detailed below, and the associated CRE properties.

\subsection{Recast of MsFEM}\label{section:recastMsFEM}

From the MsFEM formulations presented in Section~\ref{section:MSFEM}, it is not straightforward to set up a clear adaptive strategy of MsFEM parameters. This is due to the fact that there is no hierarchy appearing in the MsFEM approximation, going from the reference model~\eqref{modelpbweak} to the numerical approximation. Such a hierarchy was recently enlightened in~\cite{HEN14}, using similarities between the MsFEM approach and the Variational Multiscale (VMS) Method~\cite{HUG98}, and enabling \textit{a posteriori} adaptive strategies similar to those presented in~\cite{LAR07} for VMS. Following~\cite{HEN14}, we first recast MsFEM so that the hierarchical structure becomes apparent.

\medskip

We recall that $V = \mH^1_{0,\Gamma_D}(\Omega)$ (see~\eqref{eq:def_V}). Let us introduce a projection operator $I_H : V \to V^0_H \subset V$ (satisfying $I_H \circ I_H = I_H$) such that its image is $V^0_H$: $\text{Im}\, I_H = V^0_H$ (we can e.g. take the Cl\'ement interpolant), where $V_H^0$ is the P1-Lagrange finite dimensional space~\eqref{eq:def_VH}. We also introduce the fine-scale space $V^f \subset V$ defined by
\begin{equation*}
V^f=\big\{\phi \in V, \ \ I_H(\phi)=0\big\}.
\end{equation*}
We then have
\begin{equation*}
V = V^0_H \oplus V^f.
\end{equation*}
For each function $\varphi_H \in V^0_H$, we define the corrector $Q^\eps(\varphi_H )\in V^f$ as the solution to the problem:
\begin{equation}\label{eq:correctorpb2}
\forall v \in V^f, \quad \intO \Aa^\eps\nab\big(\varphi_H+Q^\eps(\varphi_H)\big)\cdot \nab v  = \intO fv  + \int_{\Gamma_N} g \, v.
\end{equation}
The associated operator $Q^\eps:V^0_H \to V^f$ is an affine operator:
\begin{equation}
  \label{eq:affine_1}
  Q^\eps(\varphi_H) = Q^\eps_0 + Q_\ell^\eps(\varphi_H)
\end{equation}
with
\begin{equation}
  \label{eq:affine_2}
Q^\eps_0 \in V^f \quad \text{s.t.} \quad 
\forall v \in V^f, \quad \intO \Aa^\eps\nab Q^\eps_0\cdot \nab v  = \intO fv  + \int_{\Gamma_N} g \, v
\end{equation}
and where $Q_\ell^\eps$ is a linear operator defined by
\begin{equation}
  \label{eq:affine_3}
Q_\ell^\eps(\varphi_H) \in V^f \quad \text{s.t.} \quad
\forall v \in V^f, \quad \intO \Aa^\eps\nab\big(\varphi_H+Q_\ell^\eps(\varphi_H)\big)\cdot \nab v = 0.
\end{equation}

We then define a reconstruction operator $R^\eps:V^0_H \to V$ as $R^\eps=I+Q^\eps$, as well as the multiscale functional space $\overline{V}^\eps_H := \{R^\eps(\varphi_H), \ \ \varphi_H \in V^0_H\} = Q^\eps_0 + \text{Span} \big\{\phi^0_i+Q_\ell^\eps(\phi^0_i)\big\}$, where $\{\phi^0_i\}$ is the set of the Lagrange basis functions of $V^0_H$. We note that $\overline{V}^\eps_H$ is an affine space.

\medskip

Let us now consider the following MsFEM-like problem: 
\begin{equation}\label{eq:formulMsFEM0}
\text{Find $\overline{u}^\eps_H \in \overline{V}^\eps_H$ such that, for any $v \in V^0_H$,}\quad \intO \Aa^\eps \nab \overline{u}^\eps_H\cdot \nab v =\intO f v  + \int_{\Gamma_N} g \, v.
\end{equation}
Problem~\eqref{eq:formulMsFEM0} is well-posed. Indeed, introducing the vector space $\widetilde{V}^\eps_H = \text{Span} \big\{\phi^0_i+Q_\ell^\eps(\phi^0_i)\big\}$ and writing $\overline{u}^\eps_H = Q^\eps_0 + \widetilde{u}^\eps_H$, problem~\eqref{eq:formulMsFEM0} can be recast as
\begin{equation*}
\text{Find $\widetilde{u}^\eps_H \in \widetilde{V}^\eps_H$ such that, for any $v \in V^0_H$,} \quad \intO \Aa^\eps \nab \widetilde{u}^\eps_H\cdot \nab v = \intO f v  + \int_{\Gamma_N} g \, v - \intO \Aa^\eps \nab Q^\eps_0 \cdot \nab v.
\end{equation*}
For any $v \in V^0_H$, we have $Q_\ell^\eps(v) \in V^f$. Thus, in view of~\eqref{eq:affine_3}, the above problem is equivalent to:
\begin{equation*}
\text{Find $\widetilde{u}^\eps_H \in \widetilde{V}^\eps_H$ such that, for any $v \in \widetilde{V}^\eps_H$,} \ \ \intO \Aa^\eps \nab \widetilde{u}^\eps_H\cdot \nab v = \intO f v_H  + \int_{\Gamma_N} g \, v_H - \intO \Aa^\eps \nab Q^\eps_0 \cdot \nab v_H,
\end{equation*}
where $v_H \in V_H^0$ is such that $v = v_H + Q_\ell^\eps(v_H)$. 
The above problem is obviously well-posed thanks to the Lax-Milgram lemma.

It is next straightforward to show, from~\eqref{eq:affine_2}, \eqref{eq:affine_3} and~\eqref{eq:formulMsFEM0}, that $\overline{u}^\eps_H$ satisfies the variational formulation~\eqref{modelpbweak} of the reference problem. We hence have $\overline{u}^\eps_H = u^\eps$.

The formulation~\eqref{eq:formulMsFEM0} is as challenging to use as the formulation~\eqref{modelpbweak} (note in particular that $V^f$ is not of finite dimension, and that computing $Q^\eps_\ell(\varphi_H)$ amounts to solving a highly oscillatory problem over the whole domain $\Omega$). It is however useful as explained in what follows. 

\begin{remark}
The solution $\overline{u}^\eps_H$ to~\eqref{eq:formulMsFEM0} is also solution to the following problem:
\begin{equation}\label{eq:formulMsFEM02}
\text{Find $\overline{u}^\eps_H \in \overline{V}^\eps_H$ such that, for any $v \in \overline{V}^\eps_H$,}\quad \intO \Aa^\eps \nab \overline{u}^\eps_H\cdot \nab v =\intO f v  + \int_{\Gamma_N} g \, v.
\end{equation}
This problem is a Galerkin type problem, whereas~\eqref{eq:formulMsFEM0} is a Petrov-Galerkin type problem. The problems~\eqref{eq:formulMsFEM0} and~\eqref{eq:formulMsFEM02} are equivalent, but they will lead to different approximations, as shown below.
\end{remark}

Starting from the formulations~\eqref{eq:formulMsFEM0} or~\eqref{eq:formulMsFEM02}, we can construct a fully discrete MsFEM formulation over the MsFEM space $V^\eps_H$ by making the following simplifying changes (denoted \textbf{S\#}):
\begin{itemize}
\item \textbf{S1}: Replace the right-hand side term with $0$ in the fine-scale problem~\eqref{eq:correctorpb2}, which makes the operator $Q^\eps$ linear and independent of the load ($Q^\eps_0=0$ and $Q^\eps \equiv Q^\eps_\ell$), and leads to a low-dimensional vector space $\text{Span} \{R^\eps(\phi^0_i)\} = \text{Span} \big\{\phi^0_i+Q_\ell^\eps(\phi^0_i)\big\}$.
\item \textbf{S2}: Restrict to $\Omega_i$ the correction $Q^\eps(\phi^0_i)$, where $\Omega_i$ is the support of $\phi_i^0$, i.e. the set of elements connected to node $i$. This leads to considering
\begin{equation}
  \label{eq:def_VH1}
  \overline{V}^\eps_{H,1} = \text{Span} \big\{\phi^0_i+Q_\ell^\eps(\phi^0_i)_{|\Omega_i} \big\}.
\end{equation}
\item \textbf{S3}: Localize the fine-scale problem~\eqref{eq:affine_3} by replacing $\Omega$ by smaller domains. In practice, these are subdomains $S_K$ including the element $K$, with prescribed homogeneous Dirichlet boundary conditions on $\partial S_K$. More precisely, for each $K \in \mT_H$ and each basis function $\phi^0_i \in V^0_H$, a local correction $\cQ^\eps_K(\phi^0_{i|K}) \in V^f(S_K)$ is defined as the solution to: 
\begin{equation*}
\forall v \in V^f(S_K), \quad \int_{S_K} \Aa^\eps \nab \big(\pi(\phi^0_{i|K}) + \cQ^\eps_K(\phi^0_{i|K})\big)\cdot \nab v =0,
\end{equation*}
where
$$
V^f(S_K) = H^1_0(S_K)
$$
and where $\pi(\phi^0_{i|K})$ is the affine function on $S_K$ which is equal to $\phi^0_{i|K}$ on $K$. Local correctors $\cQ^\eps_K$ are then connected together in a global corrector $\cQ^\eps$: for any $\phi^0_i \in V^0_H$, we set $\dis \cQ^\eps(\phi^0_i) = \sum_{K \in \mT_H} \cQ^\eps_K(\phi^0_{i|K})_{|K}$. The MsFEM functional space eventually reads
\begin{equation}
\label{eq:def_VH2}
\overline{V}^\eps_{H,2}=\text{Span} \{\cR^\eps(\phi^0_i)\} = \text{Span} \{\phi^0_i+\cQ^\eps(\phi^0_i)\}
\end{equation}
where the global reconstruction operator $\cR^\eps$ is defined by $\cR^\eps(\phi^0_i)=\phi^0_i+\cQ^\eps(\phi^0_i)$.

\item \textbf{S4}: Replace $V^f(S_K)$ by a discrete $P_1$-Lagrange FE space $V^f_h(S_K)$ associated to a fine grid $\mT_h$ obtained from a regular refinement of $\mT_H$. The mesh size $h$ should be chosen so that oscillations of the data can be accurately captured.
\end{itemize}

\begin{remark}
  \label{rem:justif_S1}
Simplification \textbf{S1} can be motivated by the following heuristics. Assume for instance that $I_H$ is the Cl\'ement interpolant. Consider $v \in V^f$. By definition of the Cl\'ement interpolant, we infer from $I_H(v) = 0$ that, for any $i$, $\dis \int_{\Omega_i} v \, p = 0$ for any continuous and piecewise affine function $p$ on $\Omega_i$. The function $v$ then satisfies a Poincar\'e inequality of the type $\dis \| v \|_{L^2(\Omega_i)} \leq C H \| \nabla v \|_{L^2(\Omega_i)}$. Since the mesh is regular, we hence get $\dis \| v \|_{L^2(\Omega)} \leq C H \| \nabla v \|_{L^2(\Omega)}$. It is thus expected that the left- and the right-hand sides of the corrector problem~\eqref{eq:correctorpb2} have different orders of magnitude.
\end{remark}

From simplifications \textbf{S1}-\textbf{S4}, a general MsFEM formulation with possible oversampling is obtained. For each $K \in \mT_H$ and each basis function $\phi^0_i \in V^0_H$, the local correction $\cQ^\eps_{h,K}(\phi^0_{i|K}) \in V^f_h(S_K)$ is defined as the solution to:
\begin{equation*}
\forall v \in V^f_h(S_K), \quad \int_{S_K}\Aa^\eps \nab\big(\pi(\phi^0_{i|K})+\cQ^\eps_{h,K}(\phi^0_{i|K})\big)\cdot \nab v =0.
\end{equation*}
Local correctors $\cQ^\eps_{h,K}$ are then connected together in a global corrector $\cQ^\eps_h$ defined by $\dis \cQ^\eps_h(\phi^0_i) = \sum_{K \in \mT_H} \cQ^\eps_{h,K}(\phi^0_i)_{|K}$ for any $\phi^0_i \in V^0_H$. The MsFEM functional space eventually reads $V^\eps_H=\text{Span} \{\cR^\eps_h(\phi^0_i)\} = \text{Span} \{\phi^0_i+\cQ^\eps_h(\phi^0_i)\}$ where the global reconstruction operator $\cR^\eps_h$ is defined by $\cR^\eps_h(\phi^0_i)=\phi^0_i+\cQ^\eps_h(\phi^0_i)$.

\medskip

The MsFEM approximation of $u^\eps$ is $u^\eps_H=\cR^\eps_h(u_H) \in V^\eps_H$ solution to the problem
\begin{equation*}
\forall v \in V^0_H, \quad \sum_{K\in \mT_H}\int_K \Aa^\eps \nab u^\eps_H\cdot \nab v = \intO f v +  \int_{\Gamma_N} g \, v,
\end{equation*}
which is a Petrov-Galerkin formulation (in the spirit of~\eqref{eq:formulMsFEM0}), or to the problem
\begin{equation}
  \label{eq:lundi}
\forall v \in V^\eps_H, \quad \sum_{K\in \mT_H}\int_K \Aa^\eps \nab u^\eps_H\cdot \nab v = \intO f v + \int_{\Gamma_N} g \, v,
\end{equation}
which is a Galerkin formulation (in the spirit of~\eqref{eq:formulMsFEM02}). Even though the continuous formulations~\eqref{eq:formulMsFEM0} and~\eqref{eq:formulMsFEM02} are equivalent, the two above discrete formulations are a priori not equivalent.

\begin{remark}
When choosing $S_K=K$, the Galerkin formulation~\eqref{eq:lundi} is exactly the MsFEM formulation described in Section~\ref{section:conforming}.
\end{remark}

\begin{remark}
In~\cite{HEN14}, a conforming projection $P_{H,h}$ that maps piecewise continuous functions on $\mT_H$ to elements of $V^f_h$ is introduced, in order to define a conforming global corrector $\widetilde{\cQ}^\eps_h=P_{H,h}(\cQ^\eps_h)$. This projector, which can be constructed from simple local averaging on element edges, is not considered here. 
\end{remark}

\subsection{Definition of error indicators using CRE}\label{section:defindicators}

The previous formulation of MsFEM indicates that 3 parameters can be tuned in an adaptive procedure in order to improve the accuracy:
\begin{itemize}
\item[$\bullet$] the local size $H_K$ of the coarse mesh $\mT_H$. Simplification \textbf{S1} amounts to ignoring the influence of $f$ and $g$ in the micro-scale equations. This assumption is valid provided $H_K$ is comparable to the characteristic length of spatial variations of the external loading and macroscale solution (capability of $\mT_H$ to represent long wavelength phenomena). Remark~\ref{rem:justif_S1} also indicates that, the smaller $H_K$ is, the better Simplification \textbf{S1} is justified. The validity of simplification \textbf{S2} is also related to $H_K$. If the set $\{H_K\}$ is not chosen correctly, the coarse-scale discretization error can dominate the overall approximation error.
\item[$\bullet$] the size $d_K$ (minimum patch radius) of the computational domains $S_K$ used to solve the micro-scale equations obtained after the simplification \textbf{S3}. If the set $\{d_K\}$ is not chosen correctly, the artificial boundary conditions set on $\partial S_K$ may strongly dominate the overall approximation error.
\item[$\bullet$] the local size $h_K$ of the fine mesh used to solve the micro-scale equations (in each $S_K$) after applying simplification \textbf{S4}. If the set $\{h_K\}$ is not chosen correctly, a fine-scale discretization error (say of order $O(h/\eps)$ if a P1 approximation is used to solve the problems in $S_K$) can dominate the overall approximation error. 
\end{itemize}

\medskip

We wish to define a robust adaptive algorithm that is able: (i) to detect regions where the MsFEM approximation is not sufficiently accurate; (ii) to adjust locally and optimally the relevant parameters among the three ones mentioned above. The first item can be directly addressed with a spatial decomposition, over the macro mesh $\mT_H$, of the \textit{a posteriori} error estimate $\Delta_{MsFEM}=E_{CRE}(\widehat{u}^\eps_H,\widehat{\bq}^\eps_H)$ built in Section~\ref{section:CRE}, in order to derive local error estimates in each element $K$. As regards the second item, it requires the definition of specific error indicators associated to each error source.

\medskip

In order to set up such error indicators, we introduce two intermediate reference problems, denoted \textbf{PR1} and \textbf{PR2} below. 

The first intermediate reference problem, denoted \textbf{PR1}, is obtained from the initial reference problem~\eqref{modelpbweak} (recast as~\eqref{eq:formulMsFEM0} or~\eqref{eq:formulMsFEM02} and denoted \textbf{PR0}) by applying simplifications \textbf{S1} and \textbf{S2}. The functional space associated with \textbf{PR1} is thus the space $\overline{V}^\eps_{H,1}$ defined by~\eqref{eq:def_VH1}, and the solution to \textbf{PR1} is denoted $\overline{u}_{H,1}^\eps$. We note that, in \textbf{PR1}, fine-scale computations are still defined over the whole domain $\Omega$ (infinite oversampling). Assuming for instance that we adopt as \textbf{PR0} the Galerkin formulation~\eqref{eq:formulMsFEM02}, the problem \textbf{PR1} reads:
\begin{equation*}
\text{Find $\overline{u}^\eps_{H,1} \in \overline{V}^\eps_{H,1}$ such that, for any $v \in \overline{V}^\eps_{H,1}$,}\quad \sum_{K\in \mT_H}\int_K \Aa^\eps \nab \overline{u}^\eps_{H,1} \cdot \nab v =\intO f v  + \int_{\Gamma_N} g \, v.
\end{equation*}

\medskip

The second intermediate reference problem, denoted \textbf{PR2}, is obtained from the initial reference problem \textbf{PR0} by applying simplifications \textbf{S1}, \textbf{S2} and \textbf{S3}. The functional space associated with \textbf{PR2} is thus the space $\overline{V}^\eps_{H,2}$ defined by~\eqref{eq:def_VH2}, and the solution to \textbf{PR2} is denoted $\overline{u}_{H,2}^\eps$. In \textbf{PR2}, the fine-scale computations are defined over the subdomains $S_K$ (and not over $\Omega$) and are performed exactly, without any discretization. If \textbf{PR0} is again taken as the Galerkin formulation~\eqref{eq:formulMsFEM02}, the problem \textbf{PR2} reads:
\begin{equation*}
\text{Find $\overline{u}^\eps_{H,2} \in \overline{V}^\eps_{H,2}$ such that, for any $v \in \overline{V}^\eps_{H,2}$,}\quad \sum_{K\in \mT_H}\int_K \Aa^\eps \nab \overline{u}^\eps_{H,2} \cdot \nab v =\intO f v  + \int_{\Gamma_N} g \, v.
\end{equation*}
Note that the solutions $\overline{u}_{H,1}^\eps$ and $\overline{u}_{H,2}^\eps$ may be nonconforming.

\medskip

We next write
\begin{equation*}
u^\eps-u^\eps_H = (u^\eps-\overline{u}_{H,1}^\eps)+(\overline{u}_{H,1}^\eps-\overline{u}_{H,2}^\eps)+(\overline{u}_{H,2}^\eps-u^\eps_H), \ \ \text{hence} \ \ \vertiii{u^\eps-u^\eps_H}_H \le e_{macro}+e_{over}+e_{micro},
\end{equation*}
where:
\begin{itemize}
\item[$\bullet$] $e_{macro}=\vertiii{u^\eps-\overline{u}_{H,1}^\eps}_H$ is the part of the error due to the coarse-scale discretization;
\item[$\bullet$] $e_{over}=\vertiii{\overline{u}_{H,1}^\eps-\overline{u}_{H,2}^\eps}_H$ is the part of the error due to the localization on $S_K$ of the problems~\eqref{eq:affine_3};
\item[$\bullet$] $e_{micro}=\vertiii{\overline{u}_{H,2}^\eps-u^\eps_H}_H$ is the part of the error due to the fine-scale discretization.
\end{itemize}
These three error components can also be understood as
\begin{multline*}
e_{macro} = \lim_{S_K = \Omega, \ h_K \to 0} \vertiii{u^\eps-u^\eps_H}_H, \quad
e_{over} = \lim_{H_K \to 0, \ h_K \to 0}\vertiii{u^\eps-u^\eps_H}_H,
\\
e_{micro} = \lim_{S_K = \Omega, \ H_K \to 0}\vertiii{u^\eps-u^\eps_H}_H.
\end{multline*}

\medskip

\begin{figure}[htbp]
\begin{center}
\includegraphics[width=95mm]{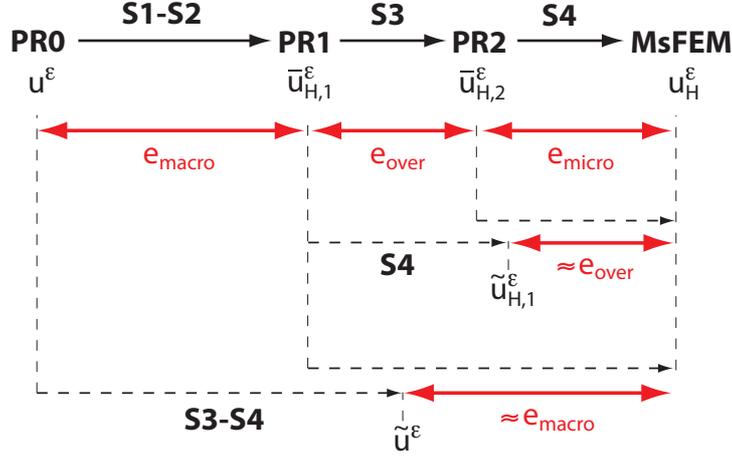}
\end{center}
\caption{Illustration of the strategy we use to define indicators for each error source.}
\label{fig:errorsources}
\end{figure}

The strategy to assess $e_{macro}$, $e_{over}$ and $e_{micro}$ from error indicators is illustrated in Figure~\ref{fig:errorsources}. First, the part $e_{micro}$ of the error can be assessed from the CRE concept considering admissibility in the sense of reference problem \textbf{PR2}. The associated error indicator thus reads $\Delta_{micro}=E_{CRE}(\widehat{u}^\eps_{H,micro},\widehat{\bq}_{micro})$, where $\widehat{u}^\eps_{H,micro} \in \overline{V}^\eps_{H,2}$ and $\widehat{\bq}_{micro}$ is a flux field that satisfies
\begin{equation}
  \label{eq:def_q2}
\forall v \in \overline{V}^\eps_{H,2}, \quad \sum_{K\in \mT_H}\int_K \widehat{\bq}_{micro} \cdot \nab v  = \intO f v  + \int_{\Gamma_N} g \, v.
\end{equation}
A first possibility would be to look for a relevant $\widehat{u}^\eps_{H,micro} \in \overline{V}^\eps_{H,2}$ and to find a field $\widehat{\bq}_{micro}$ solving~\eqref{eq:def_q2}.

In practice, we proceed in a slightly different manner. Recall indeed that our aim is to get an approximation of $e_{micro}$ and not to get a certified estimator for the error stemming from using a finite value $h_K$ for the fine mesh size. We use that flexibility to simplify the procedure and circumvent the construction of a field $\widehat{\bq}_{micro}$ that would exactly solve~\eqref{eq:def_q2}. In particular, our approach only needs the resolution of some local problems (namely, the problems~\eqref{eq:def_wK} below). There is no need to find equilibrated tractions on the coarse mesh edges, as in the estimation of the global error.

Our approach is based on the fact that, in the conforming case (i.e. when $S_K=K$), the difference $\cR^\eps(\phi^0_i) - \cR^\eps_h(\phi^0_i)$ between the basis functions of $\overline{V}^\eps_{H,2}$ and those of $V^\eps_H$ is a bubble function: it belongs to $\dis \oplus_K \widetilde{X}_0(K)$ with $\dis \widetilde{X}_0(K) = \text{Span} \left\{ (\cR^\eps(\phi^0_i) - \cR^\eps_h(\phi^0_i))_{|K} \right\}$. Rather than working in $\overline{V}^\eps_{H,2}$, our idea is to work in $\dis V^\eps_H \oplus \left( \oplus_K \widetilde{X}_0(K) \right)$. To improve robustness, we actually enrich $V^\eps_H$ by an orthogonal space. We thus define, for any coarse element $K$, the space
$$
X_{micro}(K) = \left\{ w \in \widetilde{X}_0(K), \quad \forall v \in V^\eps_H, \ \ \int_K (\nabla v)^T \Aa^\eps \nabla w = 0 \right\},
$$
consider the space $\dis \overline{W}^\eps_{H,2} = V^\eps_H \oplus \left( \oplus_K X_{micro}(K) \right)$ and make the approximation $\Delta_{micro} \approx E_{CRE}(\widehat{u}^\eps_{H,micro},\widehat{\bq}_{micro})$, where $\widehat{u}^\eps_{H,micro} \in \overline{W}^\eps_{H,2}$ and $\widehat{\bq}_{micro}$ is a flux field that satisfies
\begin{equation}
  \label{eq:def_q2bis}
\forall v \in \overline{W}^\eps_{H,2}, \quad \sum_{K\in \mT_H}\int_K \widehat{\bq}_{micro} \cdot \nab v = \intO f v  + \int_{\Gamma_N} g \, v.
\end{equation}
We choose $\widehat{u}^\eps_{H,micro} = u^\eps_H \in V^\eps_H \subset \overline{W}^\eps_{H,2}$. To build $\widehat{\bq}_{micro}$, we proceed as follows. For any coarse element $K$, we define $w_K \in X_{micro}(K)$ solution to:
\begin{equation}
  \label{eq:def_wK}
\forall v \in X_{micro}(K), \qquad \int_K (\nabla v)^T \Aa^\eps \nabla w_K = \int_K f v - \int_K (\nabla v)^T \Aa^\eps \nabla u^\eps_H.
\end{equation}
This problem is in practice solved using a fine discretization (i.e. finer than the current value $h_K$ of the fine mesh on $K$). We next set $\dis \widehat{\bq}_{micro} = \Aa^\eps \nabla u^\eps_H + \sum_K 1_K \, \Aa^\eps \nabla w_K$. It is easy to check that~\eqref{eq:def_q2bis} is satisfied. Consider indeed, with obvious notations, a generic function $\dis v = v_H + \sum_K v_K$ in $\overline{W}^\eps_{H,2}$. We compute
\begin{align*}
  & \sum_{K\in \mT_H}\int_K \widehat{\bq}_{micro} \cdot \nab v
  \\
  &=
  \sum_{K\in \mT_H}\int_K (\nabla v_H + \nabla v_K)^T \Aa^\eps (\nabla u^\eps_H + \nabla w_K)
  \\
  &=
  \sum_{K\in \mT_H} \int_K (\nabla v_H)^T \Aa^\eps \nabla u^\eps_H
  +
  \sum_{K\in \mT_H} \int_K (\nabla v_H)^T \Aa^\eps \nabla w_K
  +
  \sum_{K\in \mT_H} \int_K (\nabla v_K)^T \Aa^\eps (\nabla u^\eps_H + \nabla w_K).
\end{align*}
For the first term, we use the variational formulation satisfied by $u^\eps_H$. For the second term, we use the orthogonality condition of the definition of $X_{micro}(K)$. For the third term, we use the definition~\eqref{eq:def_wK} of $w_K$. We thus obtain that
$$
\sum_{K\in \mT_H}\int_K \widehat{\bq}_{micro} \cdot \nab v
=
\intO f \, v_H + \int_{\Gamma_N} g \, v_H + \sum_{K\in \mT_H} \int_K f \, v_K
=
\intO f \, v + \int_{\Gamma_N} g \, v,
$$
which is exactly~\eqref{eq:def_q2bis}.

\begin{remark}
If we formally assume that $h_K = 0$, then $\Delta_{micro} = 0$. Indeed, in that case, we have $\dis \widetilde{X}_0(K) = \{ 0 \}$ and thus $\dis \overline{W}^\eps_{H,2} = V^\eps_H$. We then have $\dis \widehat{\bq}_{micro} = \Aa^\eps \nabla u^\eps_H$, which leads to $\Delta_{micro} = 0$.
\end{remark}

Second, the part $e_{over}$ of the error can be assessed assuming that $\vertiii{\overline{u}_{H,1}^\eps-\overline{u}_{H,2}^\eps}_H \approx \vertiii{\widetilde{u}_{H,1}^\eps-u^\eps_H}_H$, where $\widetilde{u}_{H,1}^\eps$ is solution to the problem obtained from \textbf{PR1} and applying simplification \textbf{S4} (that is, discretization on a fine mesh of size $h$ of the corrector problems). The associated functional space is
\begin{equation}
  \label{eq:def_VH1_tilde}
  \widetilde{V}^\eps_{H,1} = \text{Span} \big\{\phi^0_i+Q_{\ell,h}^\eps(\phi^0_i)_{|\Omega_i} \big\}
\end{equation}
where $Q_{\ell,h}^\eps(\phi^0_i)$ is the solution to~\eqref{eq:affine_3} after discretization on a fine mesh of size $h$.

An error indicator for $e_{over}$ could then be constructed from the CRE concept considering admissibility in the sense of the previous problem. It reads $\Delta_{over}=E_{CRE}(\widehat{u}^\eps_{H,over},\widehat{\bq}_{over})$, where $\widehat{u}^\eps_{H,over} \in \widetilde{V}^\eps_{H,1}$ and $\widehat{\bq}_{over}$ is a flux field that satisfies
\begin{equation*}
\forall v \in \widetilde{V}^\eps_{H,1}, \quad \intO \widehat{\bq}_{over} \cdot \nab v  = \intO f v  + \int_{\Gamma_N} g \, v.
\end{equation*}
In practice, as for $e_{micro}$, we proceed differently, and work in $\widetilde{W}^\eps_{H,1}$ rather than in $\widetilde{V}^\eps_{H,1}$, where $\dis \widetilde{W}^\eps_{H,1} = V^\eps_H \oplus \left( \oplus_K X_{over}(K) \right)$. The enrichment space $X_{over}(K)$ is defined by
$$
X_{over}(K) = \left\{ w_{|K}, \quad w \in \widetilde{V}^\eps_{H,1}, \quad \forall v \in V^\eps_H, \ \ \int_K (\nabla v)^T \Aa^\eps \nabla w = 0 \right\},
$$
where $\widetilde{V}^\eps_{H,1}$ is defined by~\eqref{eq:def_VH1_tilde} (in practice, the corrector problems~\eqref{eq:affine_3} are solved on a large domain $\widetilde{S}_K$, larger than the current oversampling domain $S_K$, but small enough for~\eqref{eq:affine_3} to be affordable). We choose $\widehat{u}^\eps_{H,over} = u^\eps_H \in V^\eps_H \subset \widetilde{W}^\eps_{H,1}$ and $\dis \widehat{\bq}_{over} = \Aa^\eps \nabla u^\eps_H + \sum_K 1_K \, \Aa^\eps \nabla w_K$, where $w_K \in X_{over}(K)$ is the solution to
$$
\forall v \in X_{over}(K), \qquad \int_K (\nabla v)^T \Aa^\eps \nabla w_K = \int_K f v + \int_{\Gamma_N} g \, v - \int_K (\nabla v)^T \Aa^\eps \nabla u^\eps_H.
$$
We check that $\dis \sum_{K\in \mT_H} \int_K \widehat{\bq}_{over} \cdot \nab v = \intO f v  + \int_{\Gamma_N} g \, v$ for any $v \in \widetilde{W}^\eps_{H,1}$ and make the approximation $\Delta_{over} \approx E_{CRE}(\widehat{u}^\eps_{H,over},\widehat{\bq}_{over})$.

\begin{remark}
If we formally assume that $S_K = \Omega$, then $\Delta_{over} = 0$. Indeed, in that case, we have $\widetilde{V}^\eps_{H,1} = V^\eps_H$. The functions in $X_{over}(K)$ are thus piecewise constant, and hence $\dis \widehat{\bq}_{over} = \Aa^\eps \nabla u^\eps_H$, which leads to $\Delta_{over} = 0$.
\end{remark}

Third, the part $e_{macro}$ of the error can be assessed assuming that $\vertiii{u^\eps-\overline{u}_{H,1}^\eps}_H \approx \vertiii{\widetilde{u}^\eps-u^\eps_H}_H$, where $\widetilde{u}^\eps$ is solution to the problem obtained from \textbf{PR0} and applying simplifications \textbf{S3} and \textbf{S4} (localization and discretization of the corrector problems). The associated functional space is denoted $\widetilde{V}^\eps$. An error indicator for $e_{macro}$ could then be constructed from the CRE concept considering admissibility in the sense of the previous problem. It reads $\Delta_{macro}=E_{CRE}(\widehat{u}^\eps_{H,macro},\widehat{\bq}_{macro})$, where $\widehat{u}^\eps_{H,macro} \in \widetilde{V}^\eps$ and $\widehat{\bq}_{macro}$ is a flux field that satisfies
\begin{equation*}
\forall v \in \widetilde{V}^\eps, \quad \intO \widehat{\bq}_{macro} \cdot \nab v  = \intO f v  + \int_{\Gamma_N} g \, v.
\end{equation*}
In practice, and since $\widetilde{V}^\eps$ differs from $V_H^\eps$ by macroscale functions, we consider the space $\widetilde{W}^\eps = V_H^\eps \oplus X_{macro}$ where the enrichment space $X_{macro}$ is made of monoscale Lagrange basis functions of degree 4 over the macro mesh $\mT_H$. In other words, $X_{macro}$ is obtained from a higher-degree space $V^0_H$ and orthogonalization with respect to $V_H^\eps$. We take $\widehat{u}^\eps_{H,macro} = u^\eps_H \in V^\eps_H \subset \widetilde{W}^\eps$ and $\dis \widehat{\bq}_{macro} = \Aa^\eps \nabla u^\eps_H + \delta \widehat{\bq}_{macro}$, where $\delta \widehat{\bq}_{macro}$ is a solution to
$$
\forall v \in X_{macro}, \qquad \intO \delta \widehat{\bq}_{macro} \cdot  \nabla v = \intO f v + \int_{\Gamma_N} g \, v - \intO (\nabla v)^T \Aa^\eps \nabla u^\eps_H.
$$
Such a field $\delta \widehat{\bq}_{macro}$ is obtained as follows, using inexpensive classical equilibration procedures in the monoscale case. Introduce the mean of $\Aa^\eps$ over each element $K$, that we denote $\overline{\Aa^\eps}$. We first compute a primal field $\delta u_{macro}$ in the degree 1 space $V^0_H$ such that
$$
\forall v \in V^0_H, \qquad \intO (\nabla v)^T \, \overline{\Aa^\eps} \, \delta u_{macro} = \intO f v + \int_{\Gamma_N} g \, v - \intO (\nabla v)^T \Aa^\eps \nabla u^\eps_H.
$$
This is a monoscole problem, which is hence inexpensive. We next recover equilibrated tractions, which vary in an affine manner along the element edges. Local equilibration using the high-degree space $X_{macro}$ is eventually performed inside each element $K$ to define $\delta \widehat{\bq}_{macro|K}$. We then make the approximation $\Delta_{macro} \approx E_{CRE}(\widehat{u}^\eps_{H,macro},\widehat{\bq}_{macro})$.

\bigskip

The error indicators $\Delta_{macro}$, $\Delta_{over}$ and $\Delta_{micro}$ associated to each error source can be split into local contributions $\{\Delta^K_{macro},\Delta^K_{over},\Delta^K_{micro}\}$ in each element $K$ of the mesh $\mT_H$, and used in a greedy adaptive algorithm which is detailed in the next section. Local indicators enable to mark elements of the coarse mesh for further enrichment. Furthermore, most of the computations which are required to get the indicators can again be performed in the \textit{offline} stage of MsFEM.

\subsection{Greedy adaptive algorithm with \textit{offline}/\textit{online} strategy}\label{section:greedyalgorithm}

Using the global error estimate $\Delta_{MsFEM}$, the error indicators $\Delta_{macro}$, $\Delta_{over}$ and $\Delta_{micro}$, and their decompositions into local contributions, we design a greedy algorithm associated to adaptive MsFEM computations. The overall algorithm is presented in Figure~\ref{fig:adaptivescheme}. We emphasize that the algorithm is initialized with no oversampling ($S_K = K$) and a rough fine mesh size (namely $h_K = H_K$).

Each iteration in the algorithm is made of three groups of computations: (i) microscale (possibly costly) computations that can be performed in the \textit{offline} phase of the iteration; (ii) macroscale computations that can be performed in this \textit{offline} phase; (iii) macroscale computations that need to be performed in the \textit{online} phase of the iteration, but which are relatively cheap.

\begin{figure}[htbp]
\begin{center}
\boxed{
\begin{minipage}{0.95\textwidth}
{\small    
\begin{enumerate}
\item[1.] Initialization 
\begin{enumerate}
\item[1.1] Set the tolerance $TOL$ on the error and the refinement parameter $\gamma_r$
\item[1.2] Define a regular mesh $\mT^{(0)}_H$ with uniform mesh size $H_K=H$
\item[1.3] Set $\mR^{(0)}_{loc}=\mR^{(0)}_{glob}=\{\text{all elements $K$ of $\mT^{(0)}_H$} \}$
\item[1.4] Set $S_K=K$ and $h_K=H_K$ for all elements $K$ of $\mT^{(0)}_H$
\item[1.5] Set $k=0$
\end{enumerate}
\item[2.] 
Increment $k \rightarrow k+1$
\item[3.] \textbf{MICROSCALE OFFLINE COMPUTATIONS}
\begin{enumerate}
\item[3.1] For all $K \in \mR^{(k-1)}_{loc}$, solve~\eqref{eq:pbsolw} over $S_K$ using a P1 discretization with the mesh size $h_K$, in order to get $\bw^{\eps,K(k)}$ (see Section~\ref{section:higherorder})
\item[3.2] Store the associated local stiffness matrices
\end{enumerate}
\item[4.] \textbf{MACROSCALE OFFLINE COMPUTATIONS}
\begin{enumerate}
\item[4.1] Define the new mesh $\mT_H^{(k)}$ by refining each element $K \in \mR^{(k-1)}_{glob}$ using the local mesh size $H_K$
\item[4.2] For all $K \in \mR^{(k-1)}_{glob}$, define the new MsFEM basis functions $\phi_i^{\eps(k)}$ from $\bw^{\eps,K(k)}$
\item[4.3] Define $\eta_K$ following~\eqref{eq:signeta}
\item[4.4] Compute the $2\times 2$ (or $6\times 6$) local matrices $\M^{\Gamma_{jk}}$ (see Sections~\ref{equilitractionsconf} and~\ref{equilitractionsnonconf})
\item[4.5] Solve the local problems~\eqref{eq:dualelempb12_a}--\eqref{eq:dualelempb12_b} (and possibly~\eqref{eq:dualelempb22}) to get $\theta_i^{jk}$, $\mu_i^{jk}$ and possibly $\rho_j$
\item[4.6] Solve the local problems associated to the computation of the error indicators (see Section~\ref{section:defindicators})
\end{enumerate}
\item[5.] \textbf{MACROSCALE ONLINE COMPUTATIONS}
\begin{enumerate}
\item[5.1] Solve the global MsFEM problem~\eqref{eq:GalerkinMsFEM} (or~\eqref{eq:Petrovformul}) to get $u^{\eps(k)}_H$, and recover $\widehat{u}^{\eps(k)}_H$ 
\item[5.2] Compute $Q_i^K$ and the traction projections $\widehat{b}_i^{(j,k)}$ solving~\eqref{eq:systemi} or~\eqref{eq:systemi2}
\item[5.3] Deduce $\widehat{g}_K$ from $\eta_K^{\Gamma}$, $\M_{|\Gamma_{jk}}$ and $\widehat{b}_i^{(j,k)}$
\item[5.4] Compute $\widehat{\bq}^{\eps(k)}_{H|K}$ from $\widehat{g}_K$, $\theta_i^{jk}$, $\mu_i^{jk}$ and $\rho_i$
\item[5.5] Compute the global estimate $\Delta_{MsFEM}^{(k)}=E_{CRE}(\widehat{u}^{\eps(k)}_H,\widehat{\bq}^{\eps(k)}_H)$. If $\Delta_{MsFEM}^{(k)} \le TOL$, then STOP; otherwise go to Step 5.6
\item[5.6] Compute the indicators $\Delta^{(k)}_{macro}$, $\Delta^{(k)}_{over}$, $\Delta^{(k)}_{micro}$ and $M=\max \big(\Delta^{(k)}_{macro},\Delta^{(k)}_{over},\Delta^{(k)}_{micro} \big)$
\item[5.7] If $M=\Delta^{(k)}_{micro}$, set $\dis \mR^{(k)}_{loc} = \mR^{(k)}_{glob} = \left\{ K \in \mT_H^{(k)}, \ \frac{\Delta^{K(k)}_{micro}}{|K|} \geq \gamma_r \, \max_{\overline{K}} \frac{\Delta^{\overline{K}(k)}_{micro}}{|\overline{K}|} \right\}$, modify $h_K$ for all $K \in \mR^{(k)}_{loc}$ and go to Step 2
\item[5.8] If $M=\Delta^{(k)}_{over}$, set $\dis \mR^{(k)}_{loc} = \mR^{(k)}_{glob} = \left\{ K \in \mT_H^{(k)}, \ \frac{\Delta^{K(k)}_{over}}{|K|} \geq \gamma_r \, \max_{\overline{K}} \frac{\Delta^{\overline{K}(k)}_{over}}{|\overline{K}|} \right\}$, modify $d_K$ (i.e. enlarge $S_K$) for all $K\in\mR^{(k)}_{loc}$ and go to Step 2
\item[5.9] If $M=\Delta^{(k)}_{macro}$, set $\mR^{(k)}_{loc} = \emptyset$ and $\dis \mR^{(k)}_{glob} = \left\{ K \in \mT_H^{(k)}, \ \frac{\Delta^{K(k)}_{macro}}{|K|} \geq \gamma_r \, \max_{\overline{K}} \frac{\Delta^{\overline{K}(k)}_{macro}}{|\overline{K}|} \right\}$, modify $H_K$ for all $K\in\mR^{(k)}_{glob}$ and go to Step 2
\end{enumerate}
\end{enumerate}
}
\end{minipage}
}
\caption{Adaptive algorithm to drive the MsFEM computations.}
\label{fig:adaptivescheme}
\end{center}
\end{figure}

The adaptive procedure is associated with the following considerations:
\begin{itemize} 
\item[$\bullet$] the algorithm is initialized with the coarsest MsFEM configuration that can be employed (a regular coarse mesh, no oversampling, fine-scale problems solved with a single element);
\item[$\bullet$] when modifying the parameters $h_K$ from their initial value $h_K=H_K$ (in Step 5.7), two values are considered in the adaptive process: $h_K=\eps/5$ and $h_K=\eps/20$ (finest mesh size at the microscale);
\item[$\bullet$] when modifying the parameters $d_K$ from their initial value $d_K=H_K$ i.e. $S_K=K$ (in Step 5.8), the oversampling size is determined by adding layers of progressive thickness $\eps$, $2 \eps$, $3 \eps$, \dots around the element $K$ in the adaptive process. Such an adaptation is not performed when $h_K=H_K$;
\item[$\bullet$] refining the mesh $\mT_H$, i.e. modifying the parameters $H_K$ (in Step 5.9), is performed using a quadtree (or octree in 3D) method with nested elements. This requires to handle hanging nodes;
\item[$\bullet$] the adapted parameters $h_K$, $d_K$ and $H_K$ are chosen using classical adaptive strategies (see~\cite{LAD04}) based on convergence rates given by the \textit{a priori} estimates~\eqref{eq:aprioriboundMsFEM} or~\eqref{eq:boundover90}, and with an error target $\Delta_{macro}=\Delta_{over}=\Delta_{micro}=TOL/3$;
\item[$\bullet$] the technique proposed in~\cite{ALL06} and based on composition rules (see Section~\ref{section:higherorder}) is beneficially used throughout the adaptive process. It enables independent computations without coming back to the fine scale \textit{offline} computations (no additional costly computations). The technique is used in its $p$-refinement version in Steps 4.5 and 4.6, while it is used in its $H$-refinement version in Step 4.2. In short, this technique allows us to refine the coarse mesh without the need to solve new fine-scale problems of the type~\eqref{eq:pbsolw}. The functions $\bw^{\eps,K}$ solutions to~\eqref{eq:pbsolw} are left unchanged when the coarse mesh is refined. Note that it would not be the case if the coarse mesh were to be de-refined. 
\end{itemize}

\medskip

\begin{remark}
At initialization, we set $h_K=H_K$. The multiscale basis functions $\phi^\eps_i$ are thus equal to the P1 basis functions $\phi^0_i$. Assuming for simplicity that the tensor $\Aa^\eps = \Aa_{\rm per}(\cdot/\eps)$ is periodic and that $\dis \eps \ll \min_K H_K$, we observe that the diffusion tensor that appears when solving the global problem is the average $\overline{\Aa}$ of $\Aa_{\rm per}$ over its periodic cell. No multiscale feature is encoded in the numerical solution $u^\eps_H$. We show here that our {\it a posteriori} error estimation procedure is nevertheless still able to observe that the approximation $u^\eps_H$ is inaccurate. For simplicity, we argue on the one-dimensional case, assuming homogeneous Neumann boundary conditions and that the mean of $f$ over $\Omega$ vanishes.

Assuming that $\dis \max_K H_K \ll 1$, the MsFEM solution $u^\eps_H$ is essentially equal to the solution to $\dis -\frac{d}{dx} \left[ \overline{\Aa} \, \frac{du^\eps_H}{dx} \right] = f$. We thus have $\dis -\overline{\Aa} \, \frac{du^\eps_H}{dx}(x) = F(x) := \int_0^x f$. The CRE procedure amounts to finding an equilibrated flux $\widehat{\bq}$. In the one-dimensional context, $\widehat{\bq} = -F$. The error is estimated by $E_{CRE}(u^\eps_H,\widehat{\bq})$, with
$$
\big( E_{CRE}\left( u_H^\eps,\widehat{\bq} \right) \big)^2
=
\intO (\Aa^\eps)^{-1} \, \left( \widehat{\bq}-\Aa^\eps \frac{du_H^\eps}{dx} \right)^2
=
\intO F^2 \, (\Aa^\eps)^{-1} \, \left( 1 - \frac{\Aa^\eps}{\overline{\Aa}} \right)^2.
$$
Again using the fact that $\eps \ll 1$, we see that $\dis \big( E_{CRE}\left( u_H^\eps,\widehat{\bq} \right) \big)^2 \approx \intO F^2 \, \left( \frac{1}{\Aa^0} - \frac{1}{\overline{\Aa}} \right)$, where $\Aa^0$ is the homogenized coefficient (i.e. the harmonic average of $\Aa_{\rm per}$ here). Unless $\Aa_{\rm per}$ is constant, we have $\Aa^0 \neq \overline{\Aa}$ and thus a non-vanishing error estimate.
\end{remark}

The adaptive algorithm enables to adjust automatically (and hopefully optimally) the calculation parameters during the simulation. It is expected to provide MsFEM parameters which ensure a correct compromise between computational cost and solution accuracy. In particular, it should lead to $H_K \approx \eps$ to resolve the fine scale in the zones where scale separation and structures are missing, whereas inexpensive computations (with $h_K>\eps/20$ to solve the local problems in the \textit{offline} phase) should be sufficient in the zones where the fine-scale features of the solution are not activated. This general idea, which consists in computing fine-scale functions with minimal but sufficient accuracy (with respect to a prescribed global error tolerance), was first investigated in~\cite{STR07} for GFEM.

\section{Numerical results}\label{section:results}

In this section, we investigate the performances of the proposed methodology for \textit{a posteriori} error estimation and adaptivity in MsFEM computations. In Section~\ref{sec:num_1D}, we consider a one-dimensional non-periodic test case. We next turn in the subsequent sections to two-dimensional numerical experiments with an isotropic diffusion tensor $\Aa^\eps$, which is chosen periodic in Section~\ref{sec:num_2D_per} and non-periodic in Section~\ref{sec:num_2D_non-per}. Recall indeed that the MsFEM approach is meant to address problems where $\Aa^\eps$ has no specific (periodic, quasi-periodic, \dots) structure. It is thus important to assess its accuracy (and the accuracy of our \textit{a posteriori} estimator) in such non-periodic contexts. We eventually consider an example with a crack in Section~\ref{sec:num_2D_crack}.

For all numerical experiments, the micro mesh is obtained by a local subdivision inside each macro element (nested meshes). This way, technical difficulties due to overlapping micro and macro meshes are avoided. Regarding the two-dimensional numerical experiments, the coarse mesh and the fine meshes of each coarse element are made of quadrangles.

\subsection{Illustrative 1D example in a non-periodic setting}
\label{sec:num_1D}

\begin{figure}[htbp]
\begin{center}
\includegraphics[width=80mm]{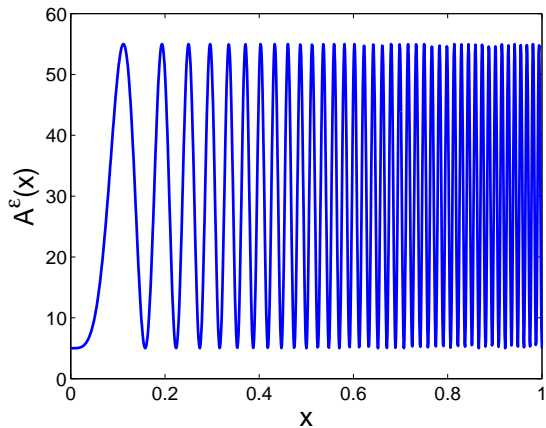}
\end{center}
\caption{Evolution of $A^\eps$ for the considered 1D non-periodic case.}
\label{fig:evolAE1Dnonper}
\end{figure}

We start with a one-dimensional problem defined in $\Omega=(0,1)$, with homogeneous Dirichlet boundary conditions at $x=0$ and $x=1$, and the non-uniform load $f(x)=x^2$. We consider a non-periodic diffusion coefficient of the form
\begin{equation*}
A^\eps(x)=5+50\sin^2(\pi x^2/\eps),
\end{equation*}
and we set $\eps=0.025$. The evolution of $A^\eps(x)$ is represented in Figure~\ref{fig:evolAE1Dnonper}.

\begin{figure}[htbp]
\begin{center}
\includegraphics[width=50mm]{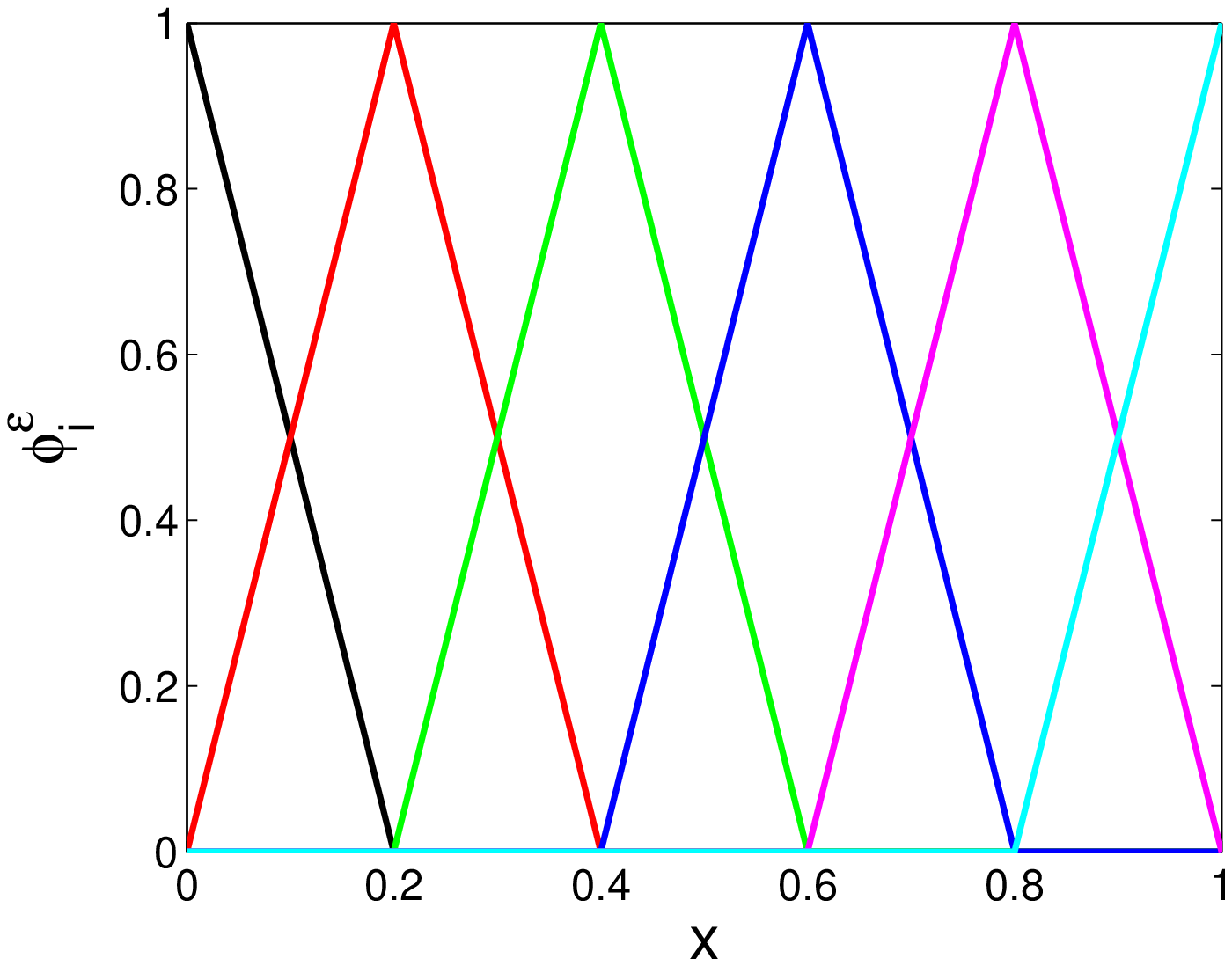}
\includegraphics[width=50mm]{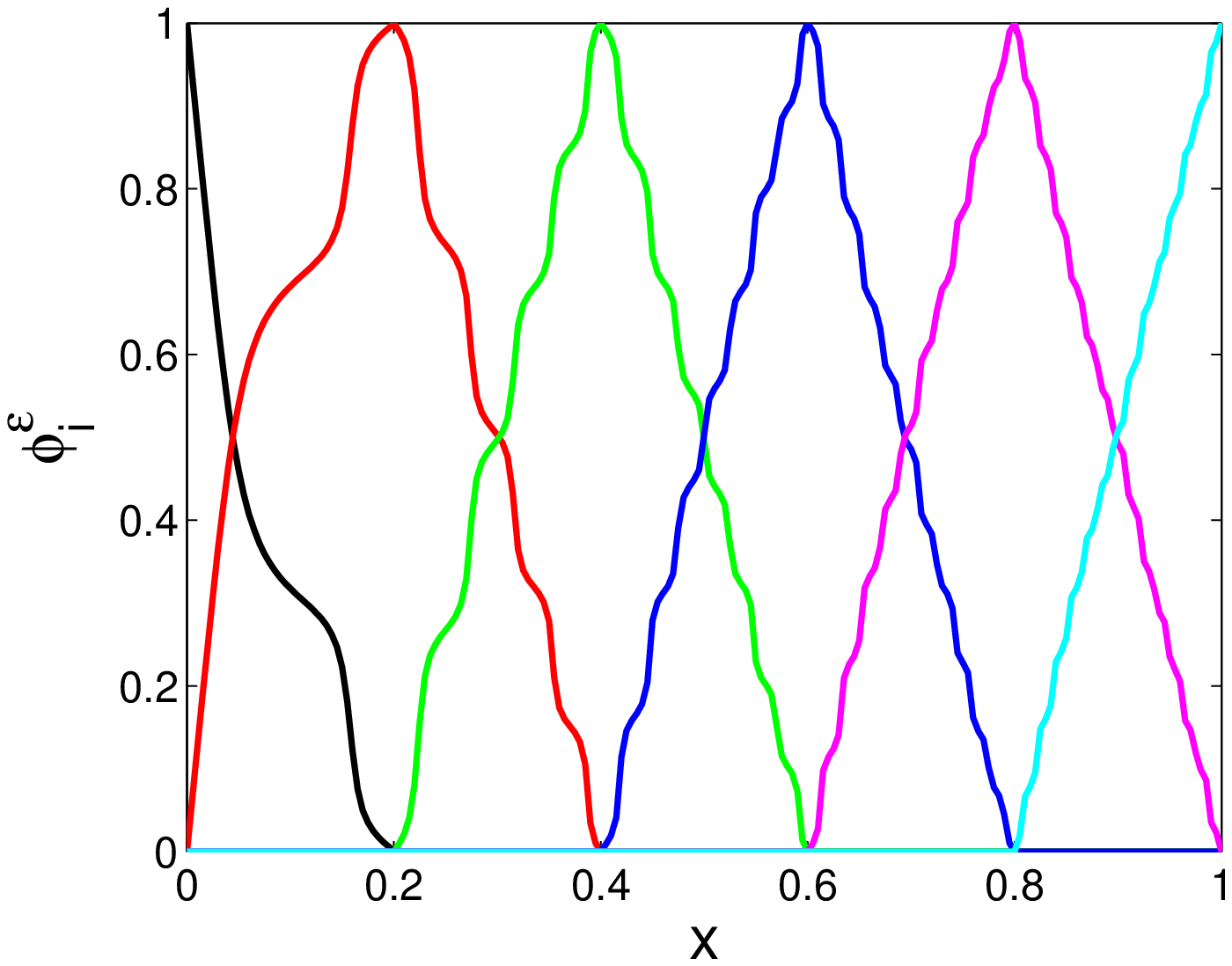}
\includegraphics[width=50mm]{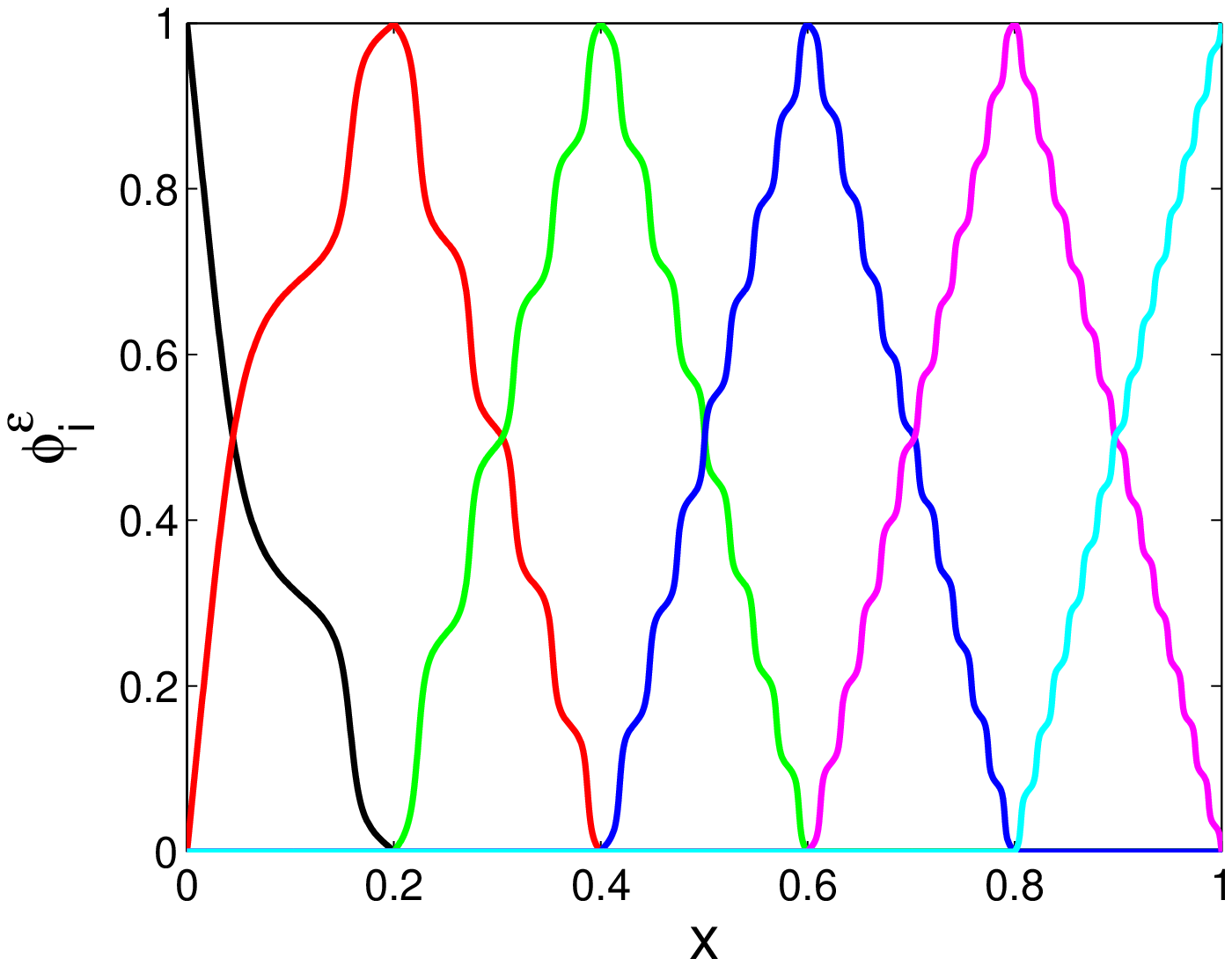}
\end{center}
\caption{Shape functions used in the MsFEM approximation, for $h_K=H_K$ (left), $h_K=\eps/5$ (center) and $h_K=\eps/20$ (right).}
\label{fig:shapefcts1Dnonper}
\end{figure}

We compute an approximate MsFEM solution using an initial coarse mesh $\mT_H$ composed of 5 elements with equal size ($H_K=0.2$). We recall that methods with or without oversampling are identical in the one-dimensional context. We show in Figure~\ref{fig:shapefcts1Dnonper} the computed MsFEM shape functions considering different fine mesh sizes (that are identical over all coarse elements $K$): $h_K=H_K$, $h_K=\eps/5$ and $h_K=\eps/20$.

\begin{figure}[htbp]
\begin{center}
\includegraphics[width=50mm]{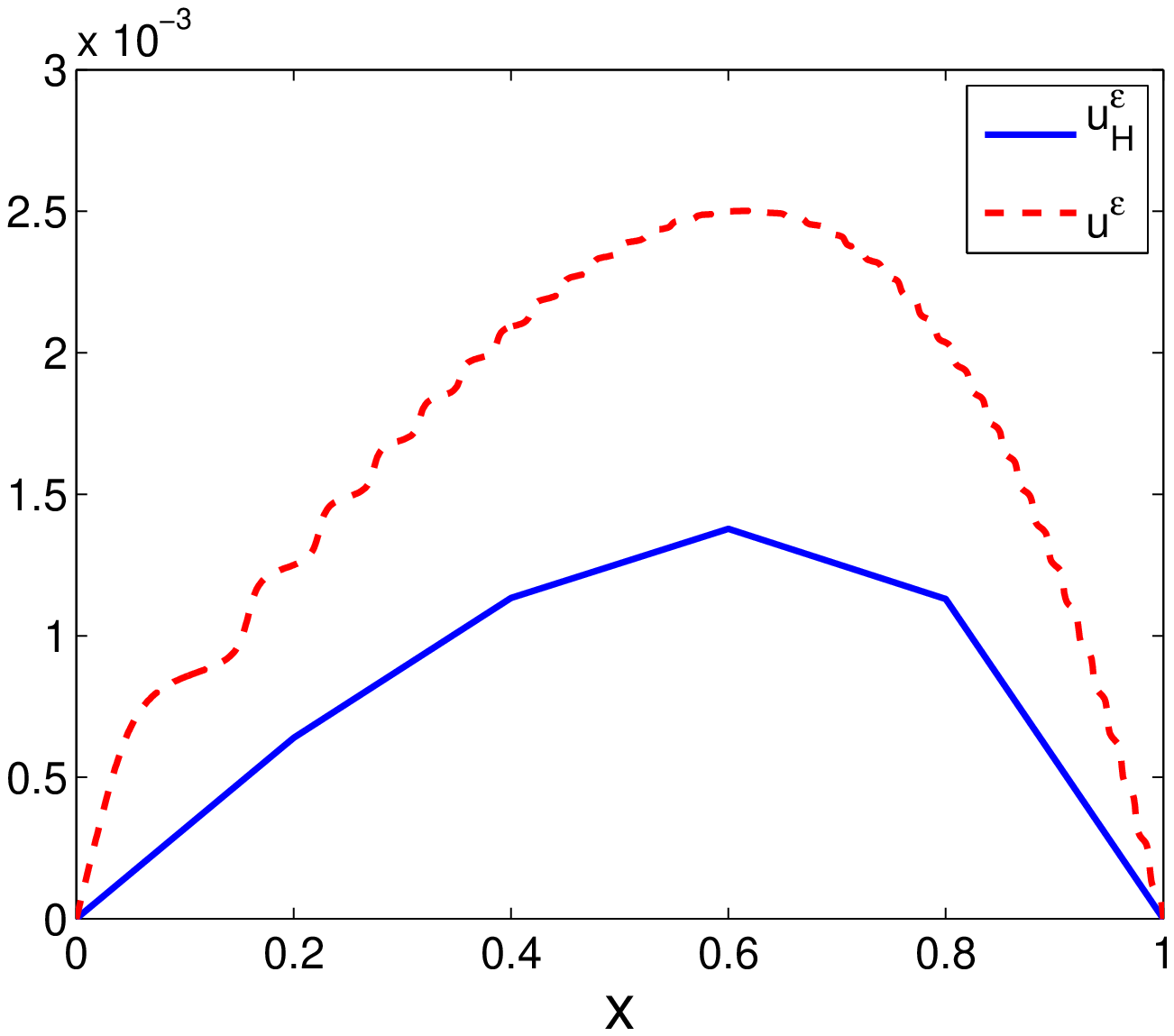}
\includegraphics[width=50mm]{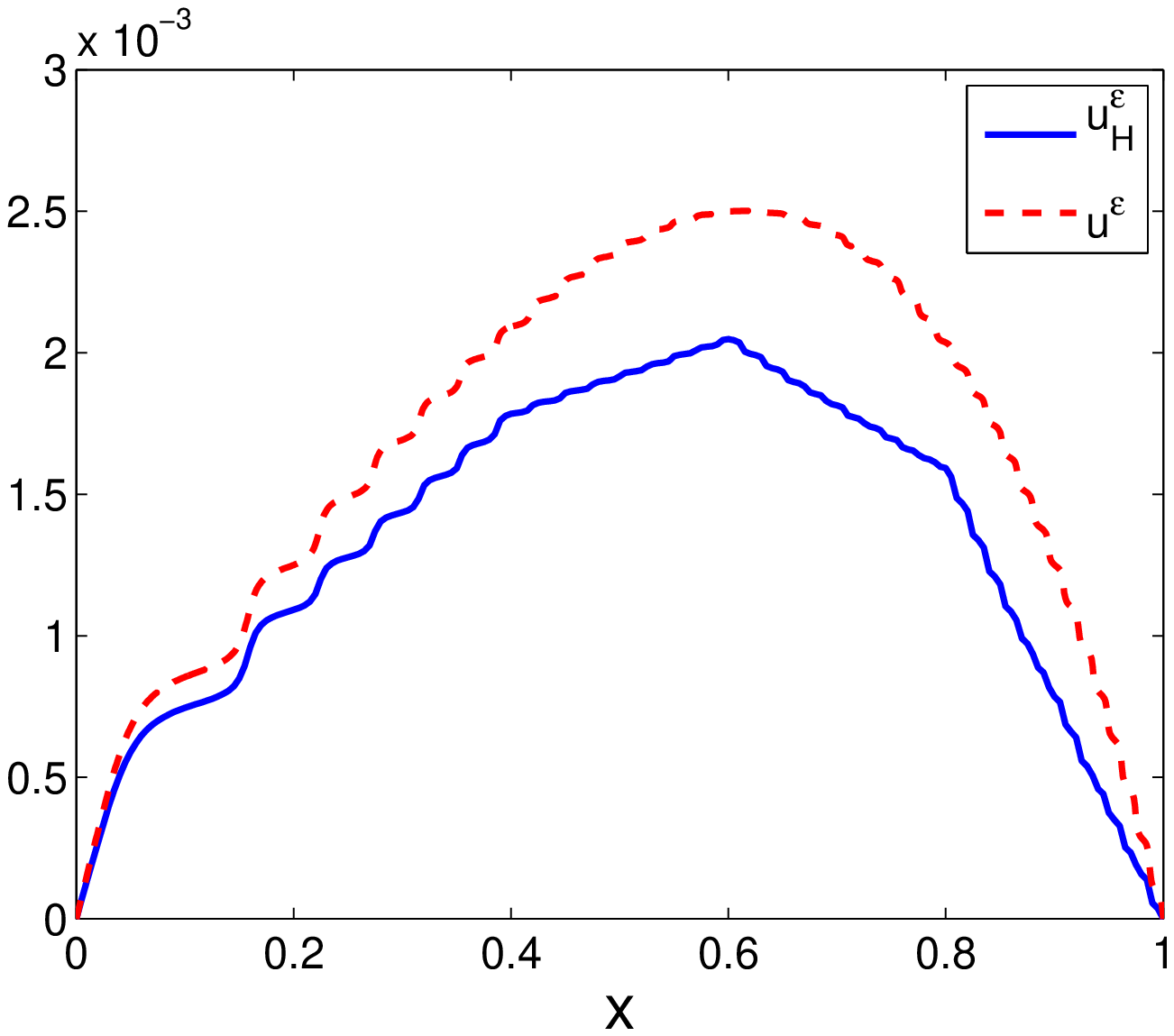}
\includegraphics[width=50mm]{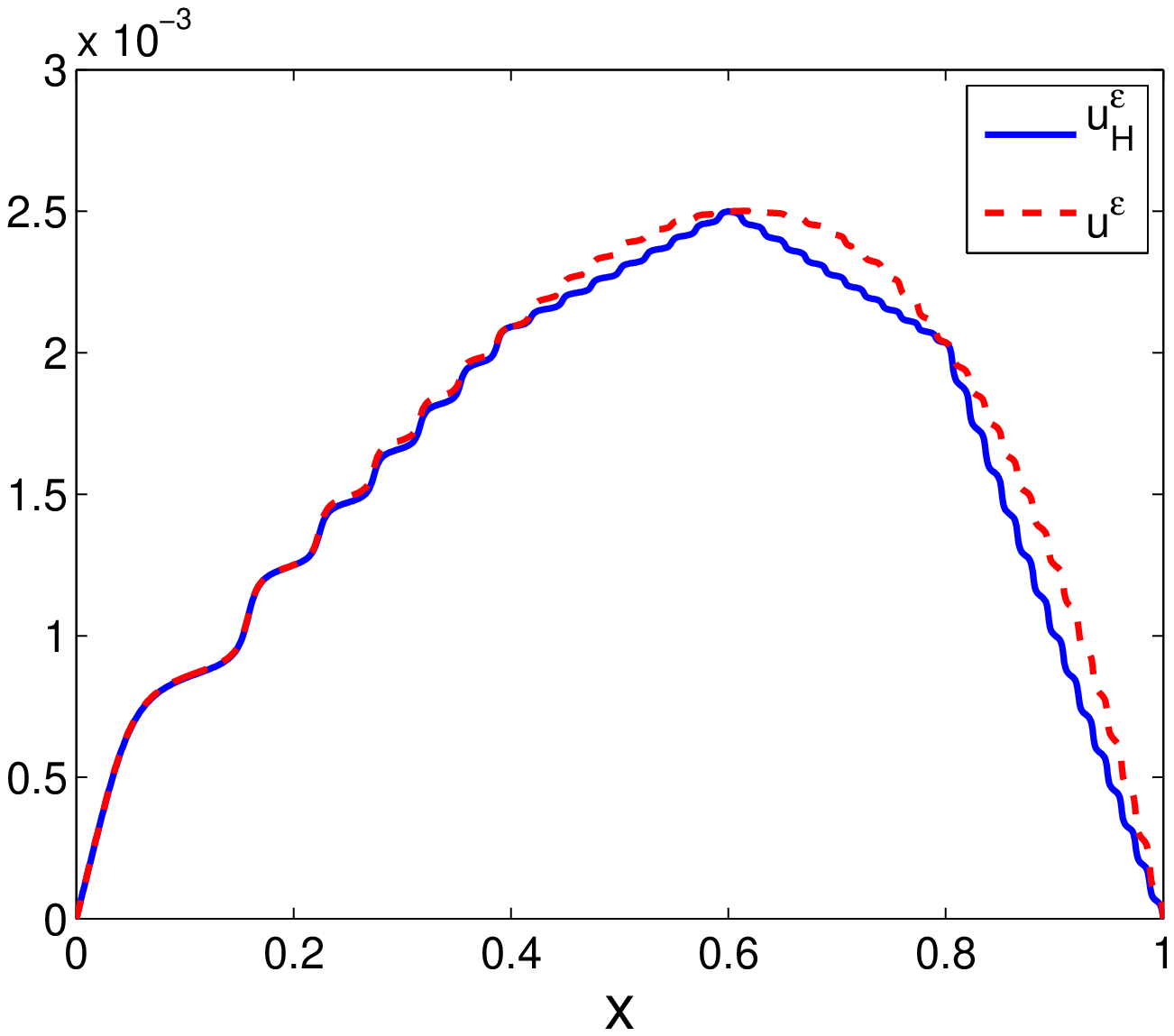}
\end{center}
\caption{Comparison between the exact solution $u^\eps$ and the approximate MsFEM solution $u^\eps_H$ when the MsFEM basis functions are computed with $h_K=H_K$ (left), $h_K=\eps/5$ (center) and $h_K=\eps/20$ (right).}
\label{fig:solMsFEM1Dnonper}
\end{figure}

The MsFEM solutions $u^\eps_H$ obtained using the different values of $h_K$ are represented in Figure~\ref{fig:solMsFEM1Dnonper}, where we also show the exact solution $u^\eps$. Obviously, the MsFEM approximation $u^\eps_H$ is not accurate when $h_K=H_K$ and when $h_K=\eps/5$. We obtain a relative error $\vertiii{u^\eps-u^\eps_H}/\vertiii{u^\eps}$ of 90\% for $h_K=\eps/5$, while it falls to 26\% for $h_K=\eps/20$. 

\begin{figure}[htbp]
\begin{center}
\includegraphics[width=75mm]{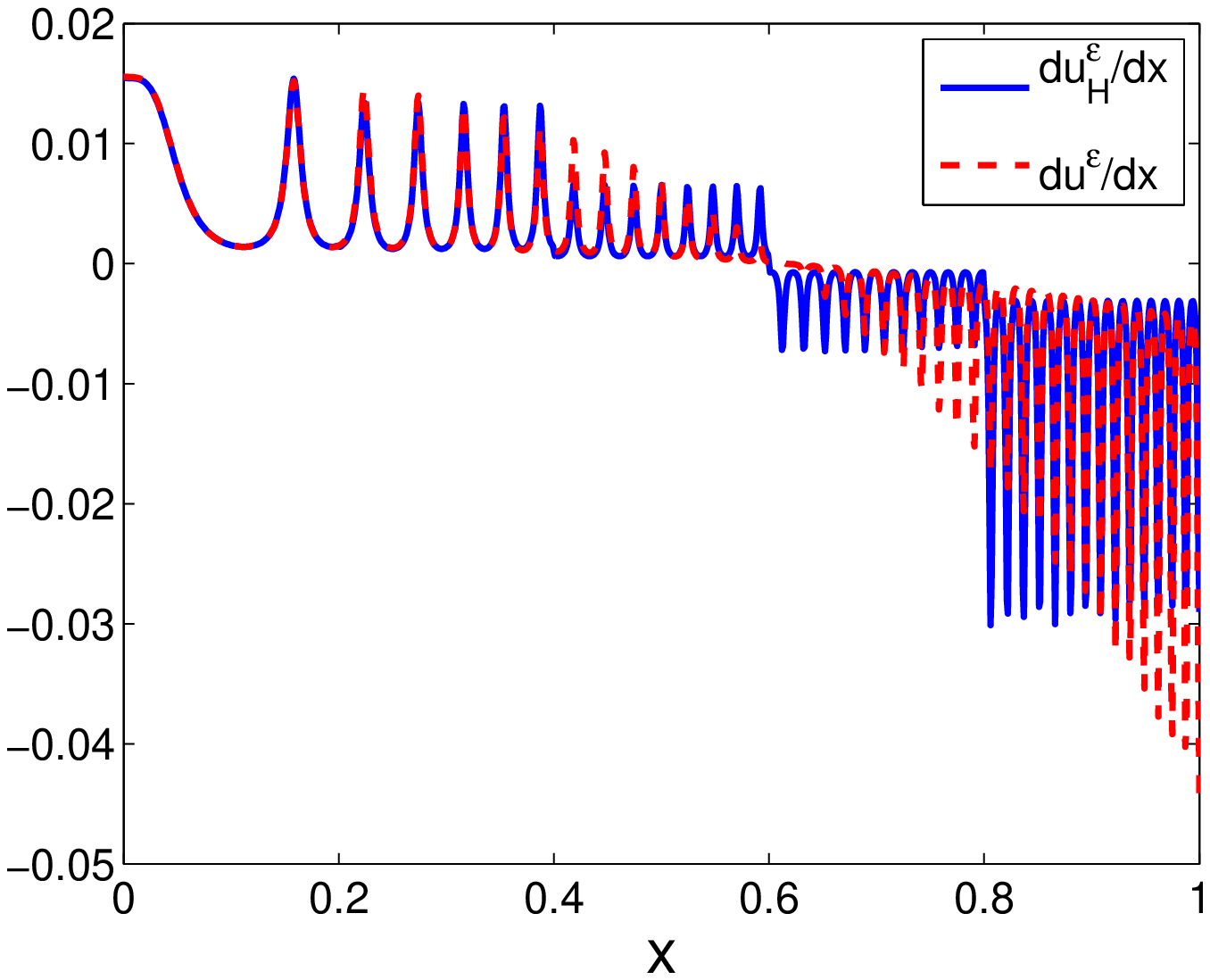}
\includegraphics[width=75mm]{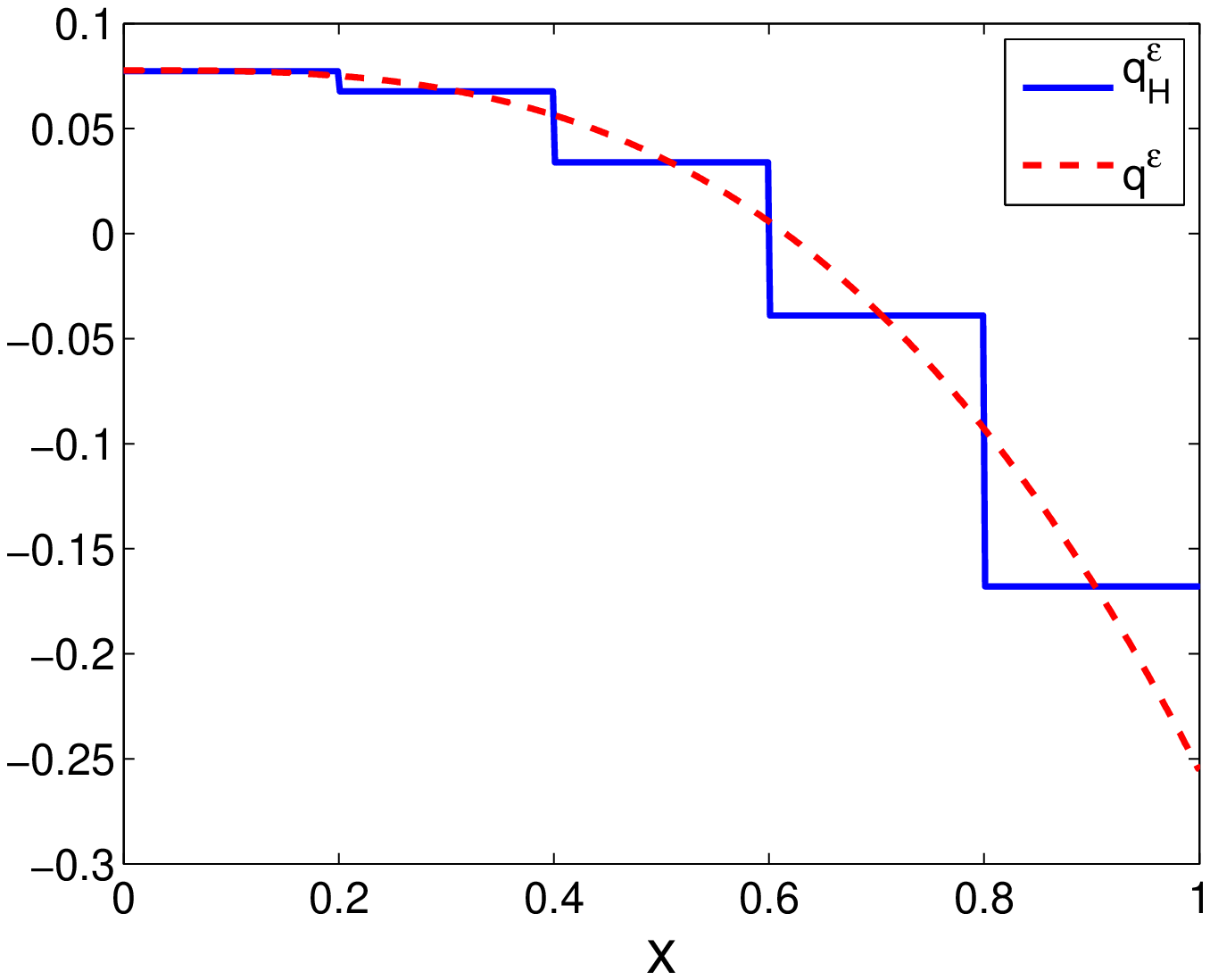}
\end{center}
\caption{Comparison between the exact solution and MsFEM solution for $h_K=\eps/20$: gradient (left) and flux (right).}
\label{fig:solMsFEM21Dnonper}
\end{figure}

For $h_K=\eps/20$, we represent in Figure~\ref{fig:solMsFEM21Dnonper} the gradient $du^\eps_H/dx$ and the flux $q^\eps_H = A^\eps du^\eps_H/dx$ of the MsFEM solution, as well as those of the exact solution. The flux of the MsFEM solution is constant in each coarse element, by definition of the MsFEM basis functions.

We now consider our adaptive algorithm described in Section~\ref{section:greedyalgorithm}. We hence start with $h_K = H_K$ and $S_K = K$. The two first iterations of the adaptive algorithm indicate that the error indicator $\Delta_{micro}$ is the highest. This leads to choosing $h_K = \eps/20$ for the five coarse elements in $\mT_H$.

\begin{figure}[htbp]
\begin{center}
\includegraphics[width=50mm]{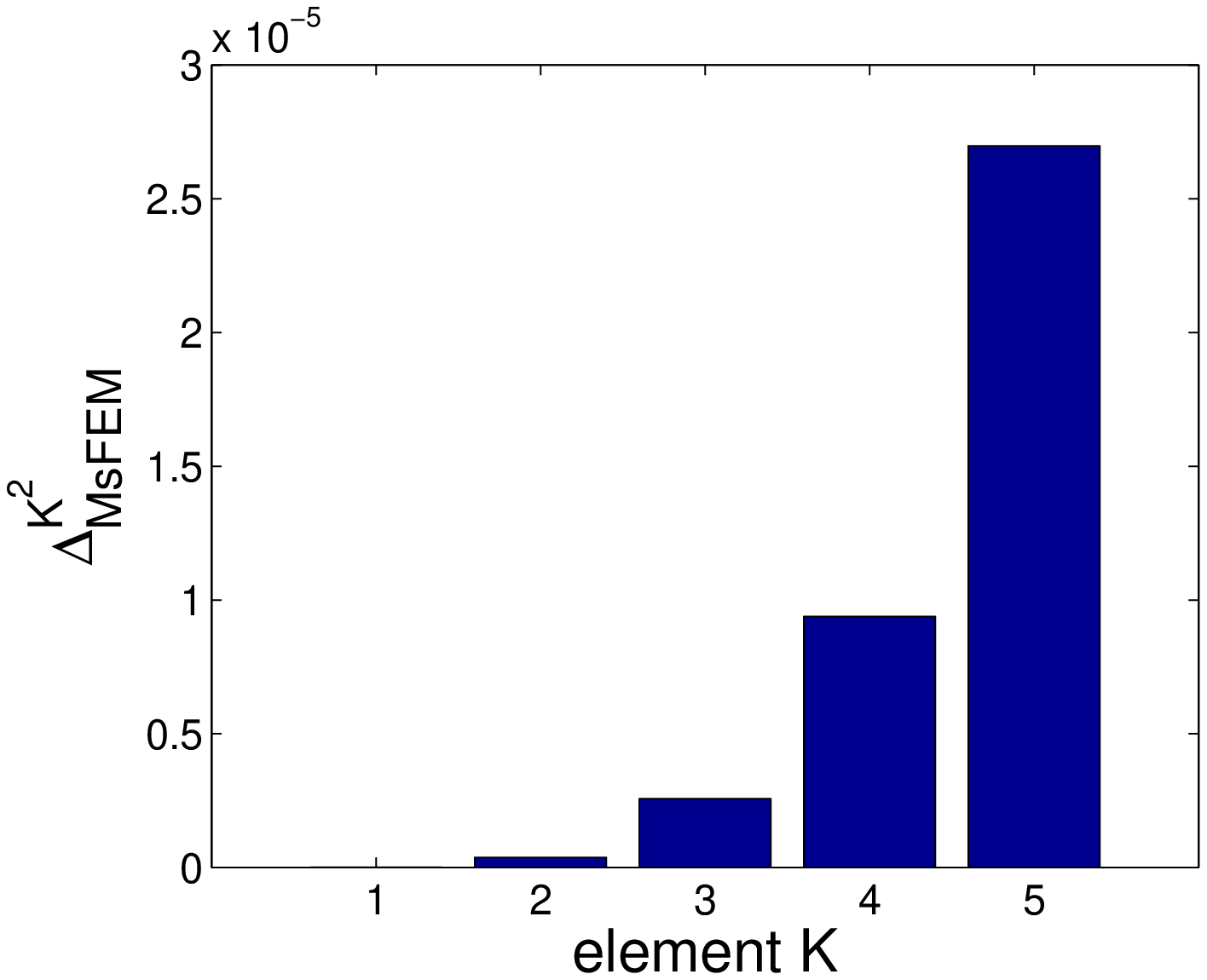}
\includegraphics[width=50mm]{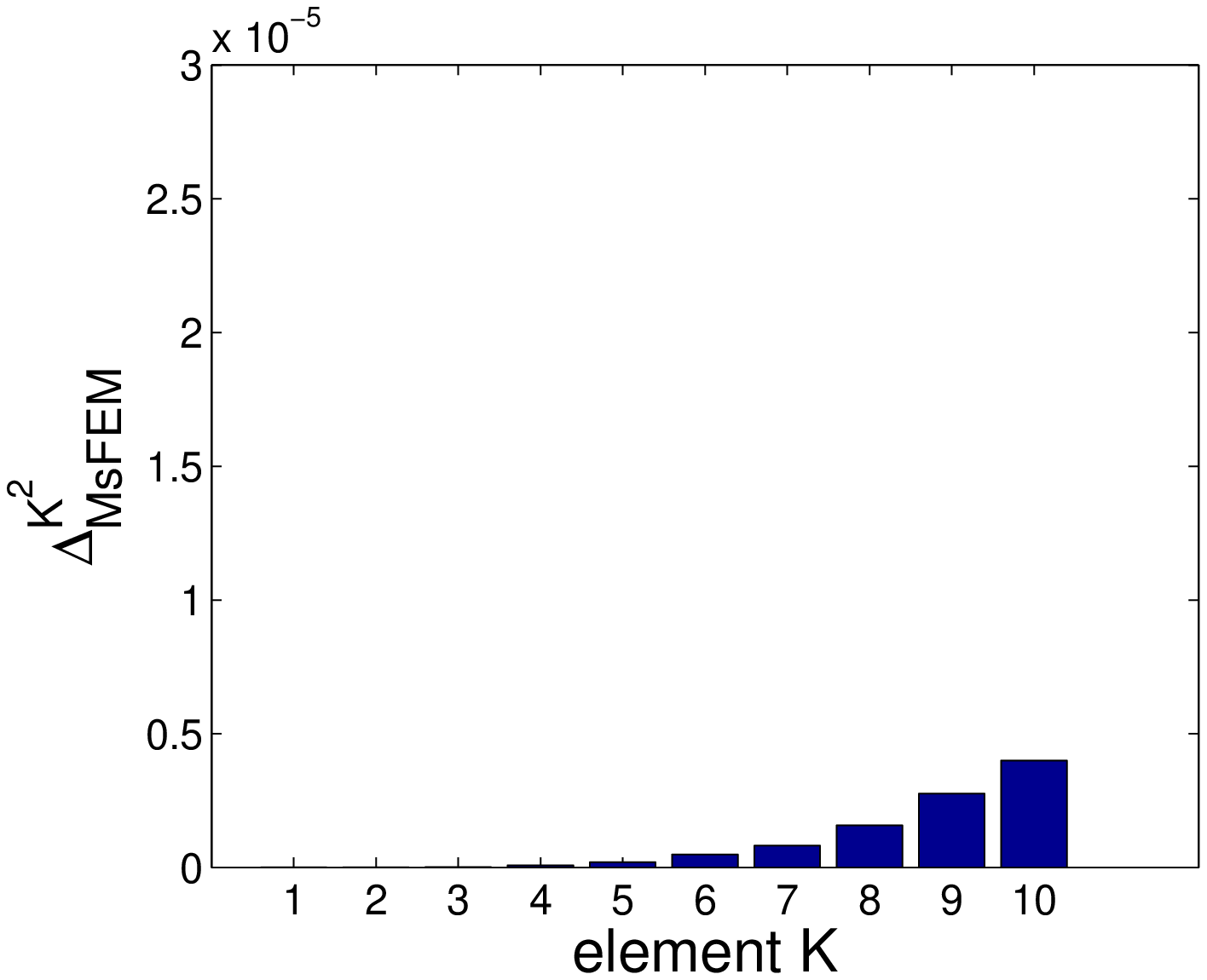}
\includegraphics[width=50mm]{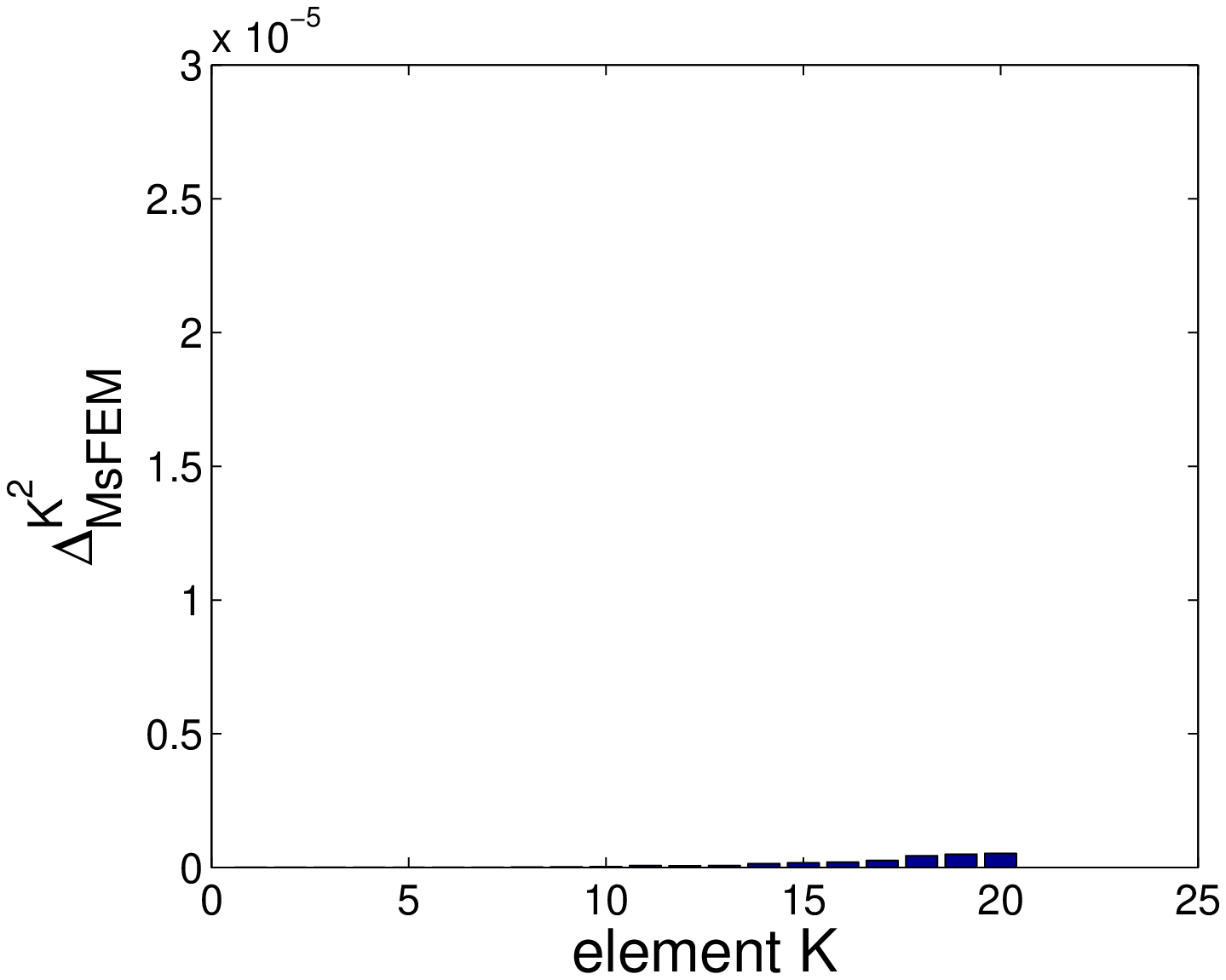}
\end{center}
\caption{Values of $(\Delta^{K}_{MsFEM})^2$ for $H_K=0.2$ (left), $H_K=0.1$ (center) and $H_K=0.05$ (right).}
\label{fig:errorMsFEM1Dnonper}
\end{figure}

To decrease the overall error, there are (at least) two possibilities. First, it is possible to uniformly refine the coarse mesh $\mT_H$. We show on Figure~\ref{fig:errorMsFEM1Dnonper} the spatial distribution of the error estimator $(\Delta_{MsFEM})^2$ for a uniform coarse mesh of size $H_K=0.2$, $H_K=0.1$ and $H_K=0.05$ (i.e. with 5, 10 and 20 coarse elements), in the case when $h_K = \eps/20$. Let us notice that the computed value of the estimate $\Delta_{MsFEM}$, considering $h_K=\eps/20$, is actually very close to the value of the exact error $\vertiii{u^\eps-u^\eps_H}$ here. This is due to the fact that in this simple 1D case, the recovered admissible flux $\widehat{q}^\eps_H$ (almost) corresponds to the exact flux $q^\eps$.

\begin{figure}[htbp]
\begin{minipage}[c]{.46\linewidth}
\includegraphics[width=75mm]{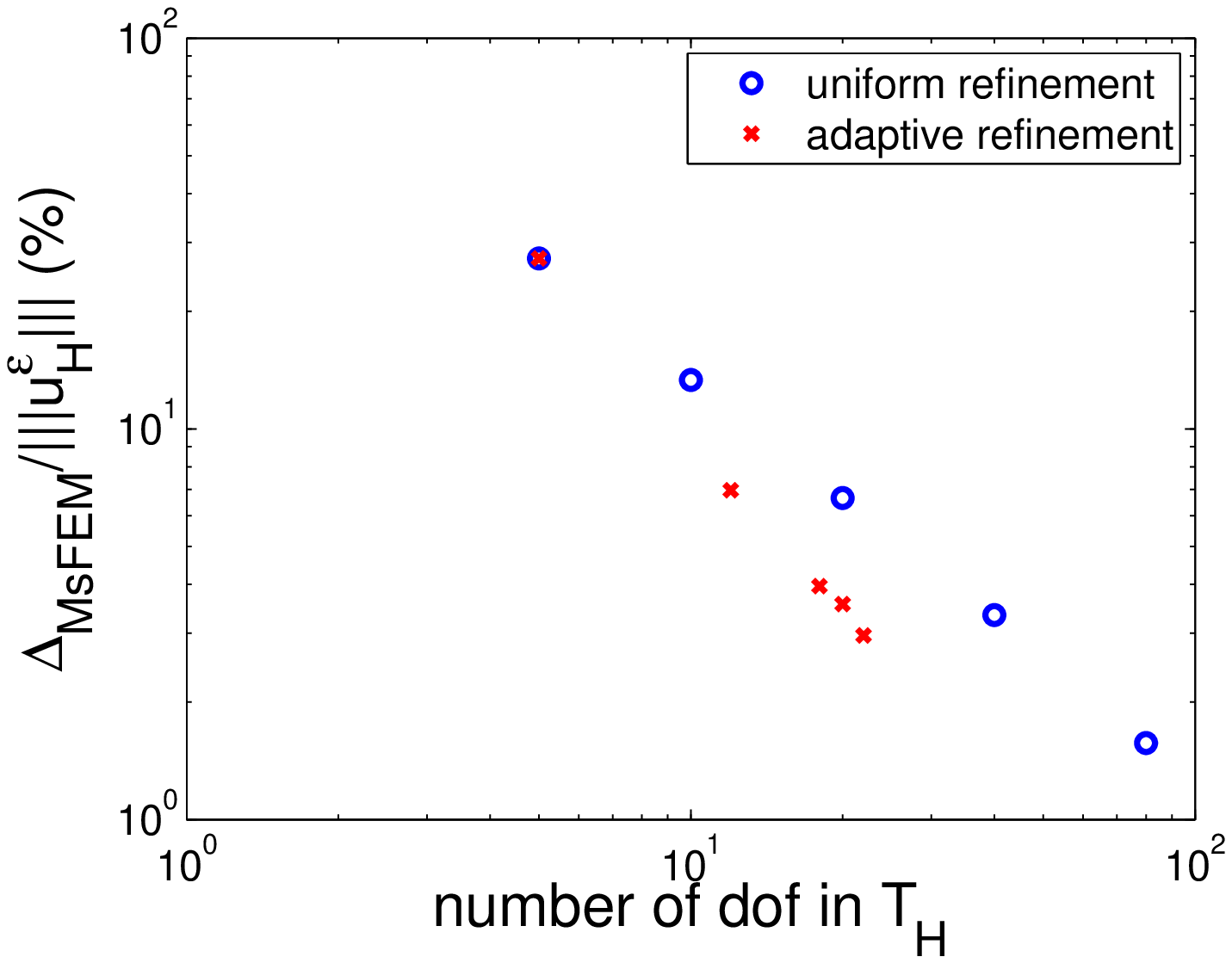} 
   \end{minipage} \hfill
\begin{minipage}[c]{.46\linewidth}
\includegraphics[width=75mm]{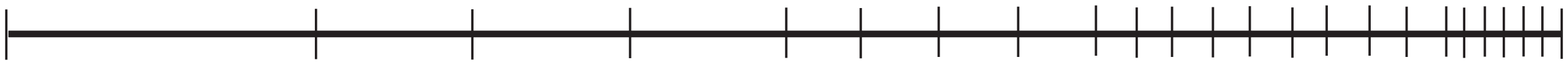}
\end{minipage}
\caption{Left: error convergence for uniform or adaptive coarse mesh refinement. Right: the coarse mesh which is obtained at the end of the adaptive algorithm.}
\label{fig:adaptive1D}
\end{figure}

As an alternative, it is possible to use the adaptive algorithm that we introduced in Section~\ref{section:greedyalgorithm}. For the third iteration and the subsequent ones, the algorithm leads to adapting the sizes $H_K$ of the coarse mesh elements, leaving the fine mesh size $h_K$ unchanged. We show in Figure~\ref{fig:adaptive1D} the convergence of the relative error estimate $\Delta_{MsFEM}/\vertiii{u^\eps_H}$ with respect to the number of degrees of freedom in $\mT_H$, for both the uniform refinement and the adaptive refinement. We clearly observe the beneficial use, in terms of convergence rate, of the adaptive algorithm. We also show in Figure~\ref{fig:adaptive1D} the adapted mesh which enables to reach a prescribed error tolerance of 3\%.

\subsection{2D example with a periodic coefficient}
\label{sec:num_2D_per}

We now consider a two-dimensional problem defined in the unit square domain $\Omega=(0,1)^2$. It is similar to the one discussed in~\cite{HOU97,ALL06}. We consider homogeneous Dirichlet boundary conditions on $\partial \Omega$, and we assume that a uniform load $f=-1$ is applied in $\Omega$. The diffusion tensor is taken as
\begin{equation*}
\Aa^\eps(\bx) = \frac{1}{\left(2+P\cos(2\pi (x_1-0.5)/\eps)\right)\left(2+P\cos(2\pi (x_2-0.5)/\eps)\right)} \ \II_2 = A^\eps(\bx) \, \II_2,
\end{equation*}
with $\eps=0.04$. The contrast parameter $P$ is chosen as $P=1.8$ here. 
For such a diffusion tensor, the corresponding homogenized tensor is explicit and reads $\Aa^0=1/2(4-P^2)^{1/2} \, \II_2$.

\begin{figure}[htbp]
\begin{center}
\includegraphics[width=75mm]{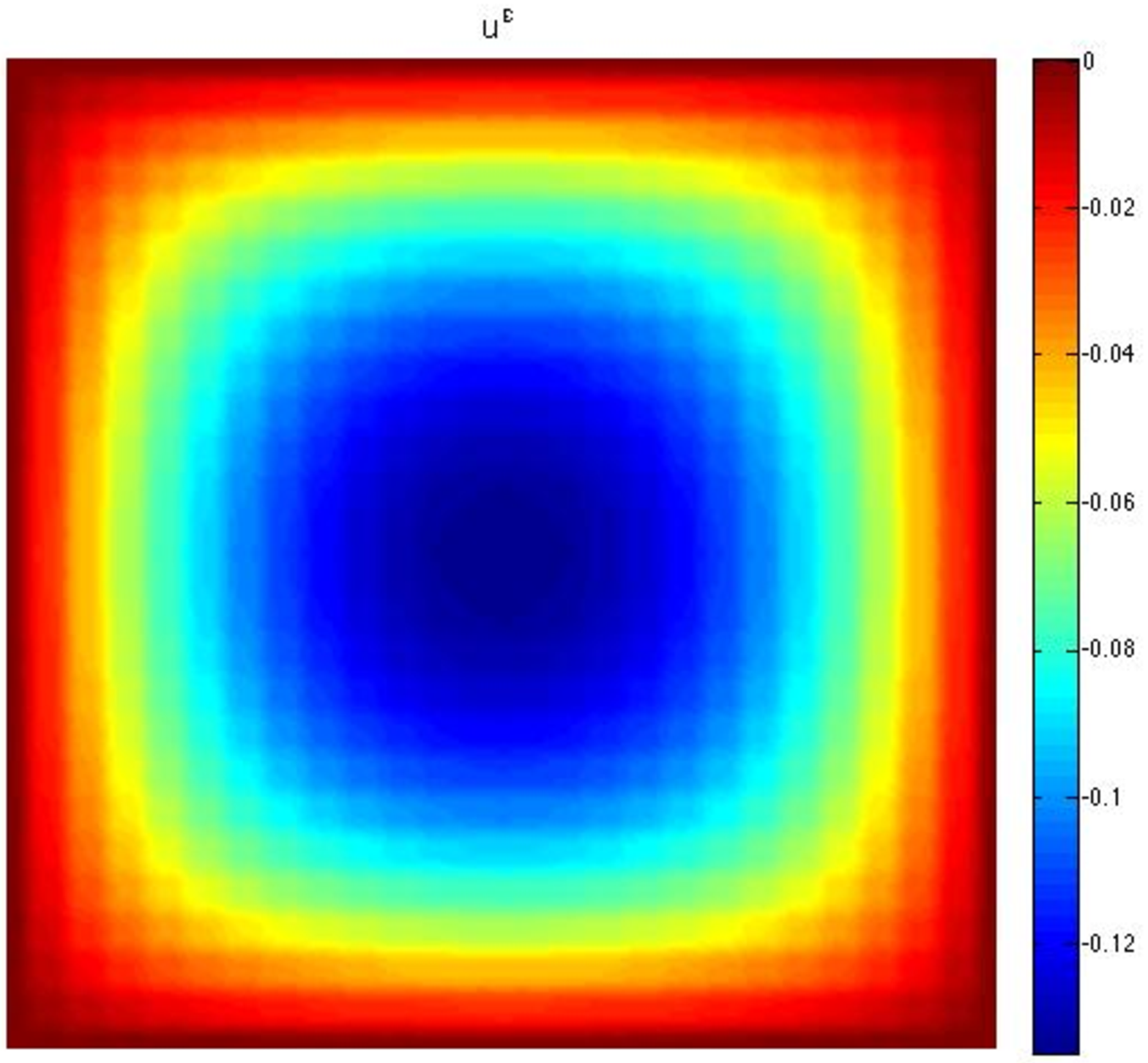}
\includegraphics[width=75mm]{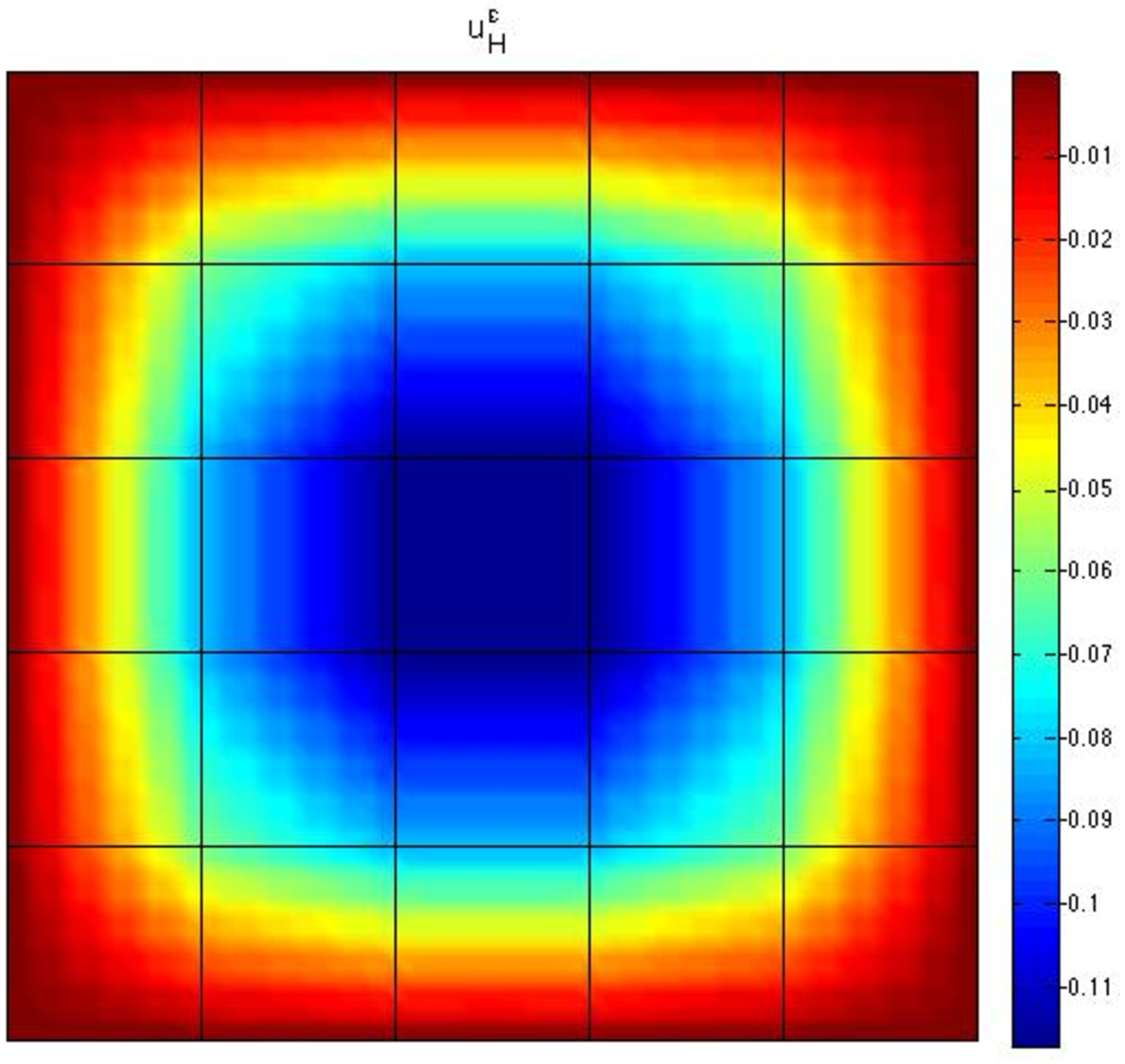}
\end{center}
\caption{2D periodic example: exact solution $u^\eps$ (left) and MsFEM solution $u^\eps_H$ (right).}
\label{fig:solMsFEM2Dperu}
\end{figure}
\begin{figure}[htbp]
\begin{center}
\includegraphics[width=75mm]{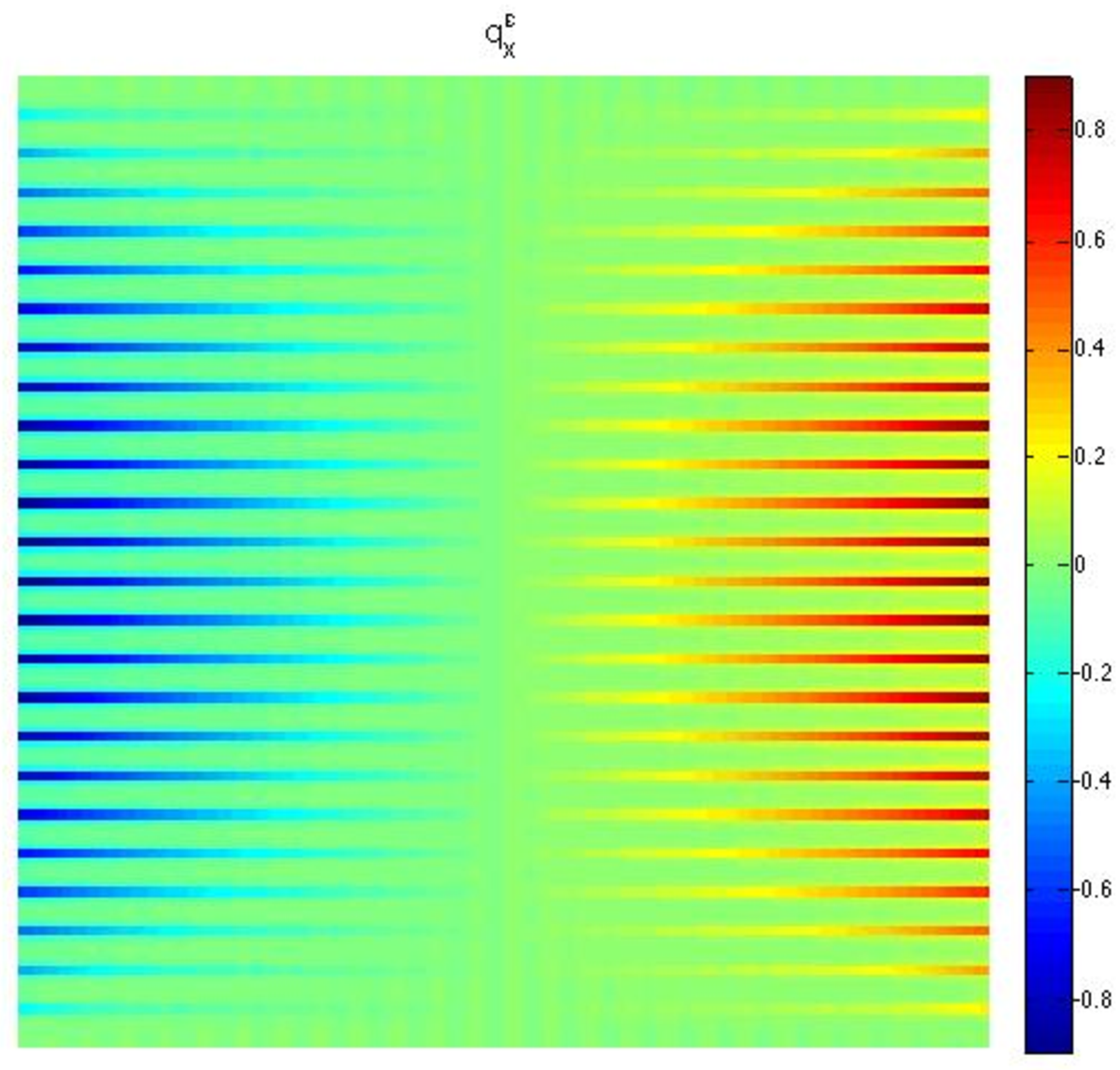}
\includegraphics[width=75mm]{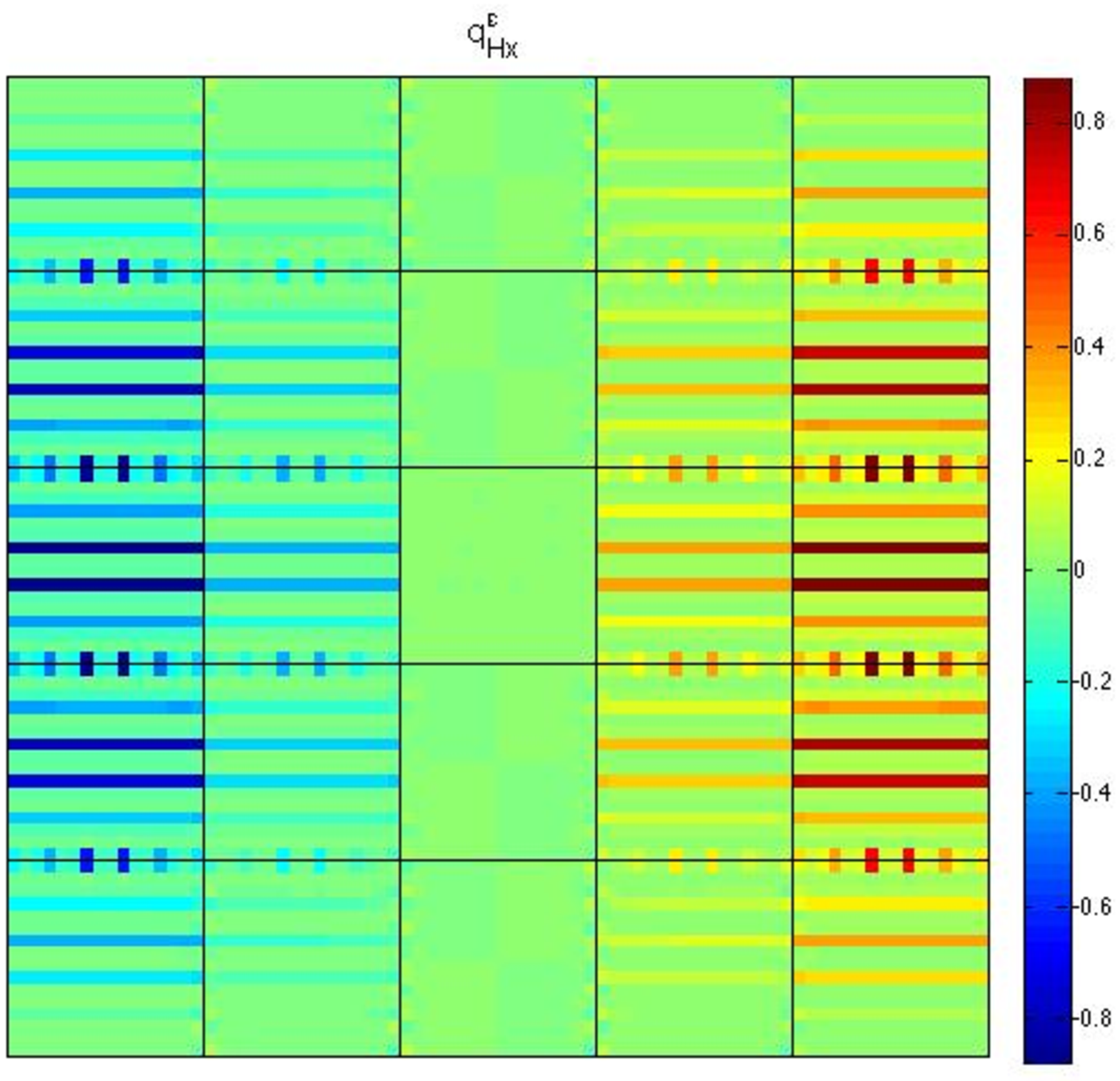}
\end{center}
\caption{2D periodic example: exact flux $\bq^\eps$ (left) and MsFEM flux $\bq^\eps_H$ (right) (we only show the first component of the flux; the figure for the second component is similar).}
\label{fig:solMsFEM2Dperflux}
\end{figure}

We consider an initial coarse mesh $\mT_H$ made of $5\times 5$ macro elements. Choosing $h_K=\eps/3$ for any $K$, and with no oversampling, the approximate MsFEM solution is represented in the Figures~\ref{fig:solMsFEM2Dperu} and~\ref{fig:solMsFEM2Dperflux} and compared to the exact solution (obtained by using a $500 \times 500$ fine mesh): solutions $u^\eps$ and $u^\eps_H$ are shown in Figure~\ref{fig:solMsFEM2Dperu} and fluxes $\bq^\eps$ and $\bq^\eps_H$ are shown in Figure~\ref{fig:solMsFEM2Dperflux}.

\begin{figure}[htbp]
\begin{center}
\includegraphics[width=75mm]{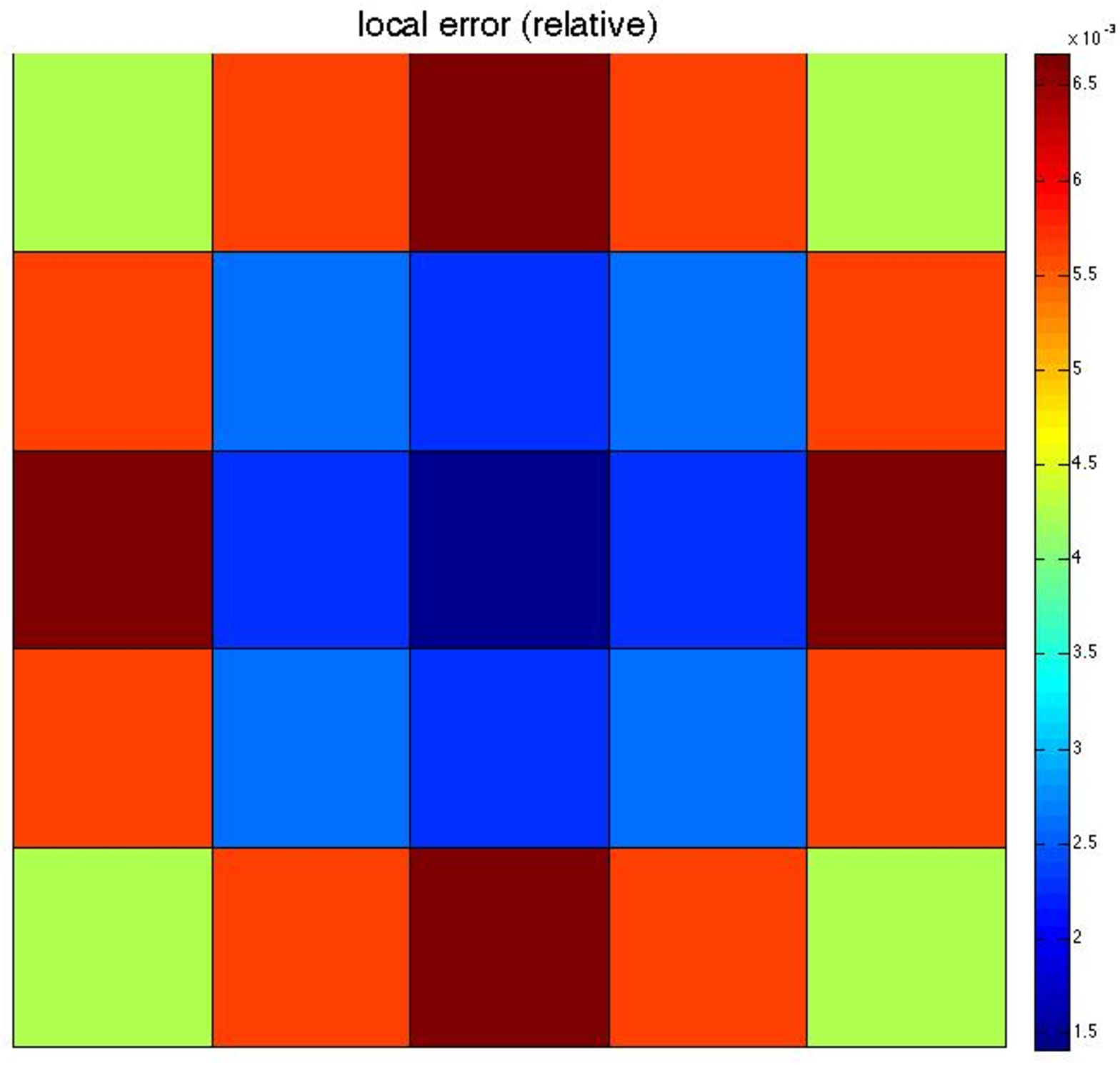} 
\includegraphics[width=75mm]{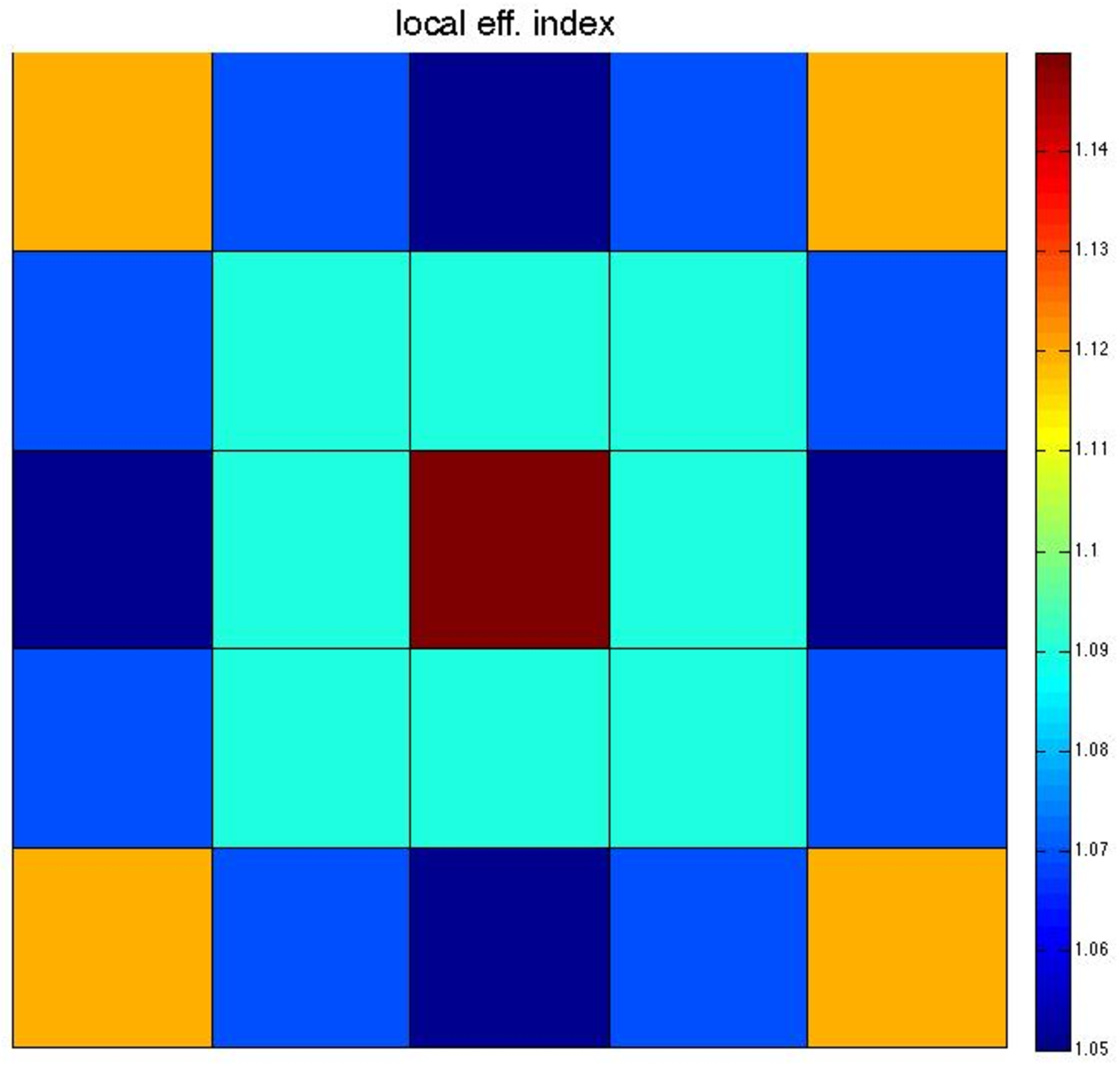} 
\end{center}
\caption{2D periodic example: values of $(\Delta^K_{MsFEM})^2$ for the specific computed MsFEM solution (left), and local effectivity indices (right).}
\label{fig:errorMsFEM2Dper}
\end{figure}

\begin{figure}[htbp]
\begin{center}
\includegraphics[width=75mm]{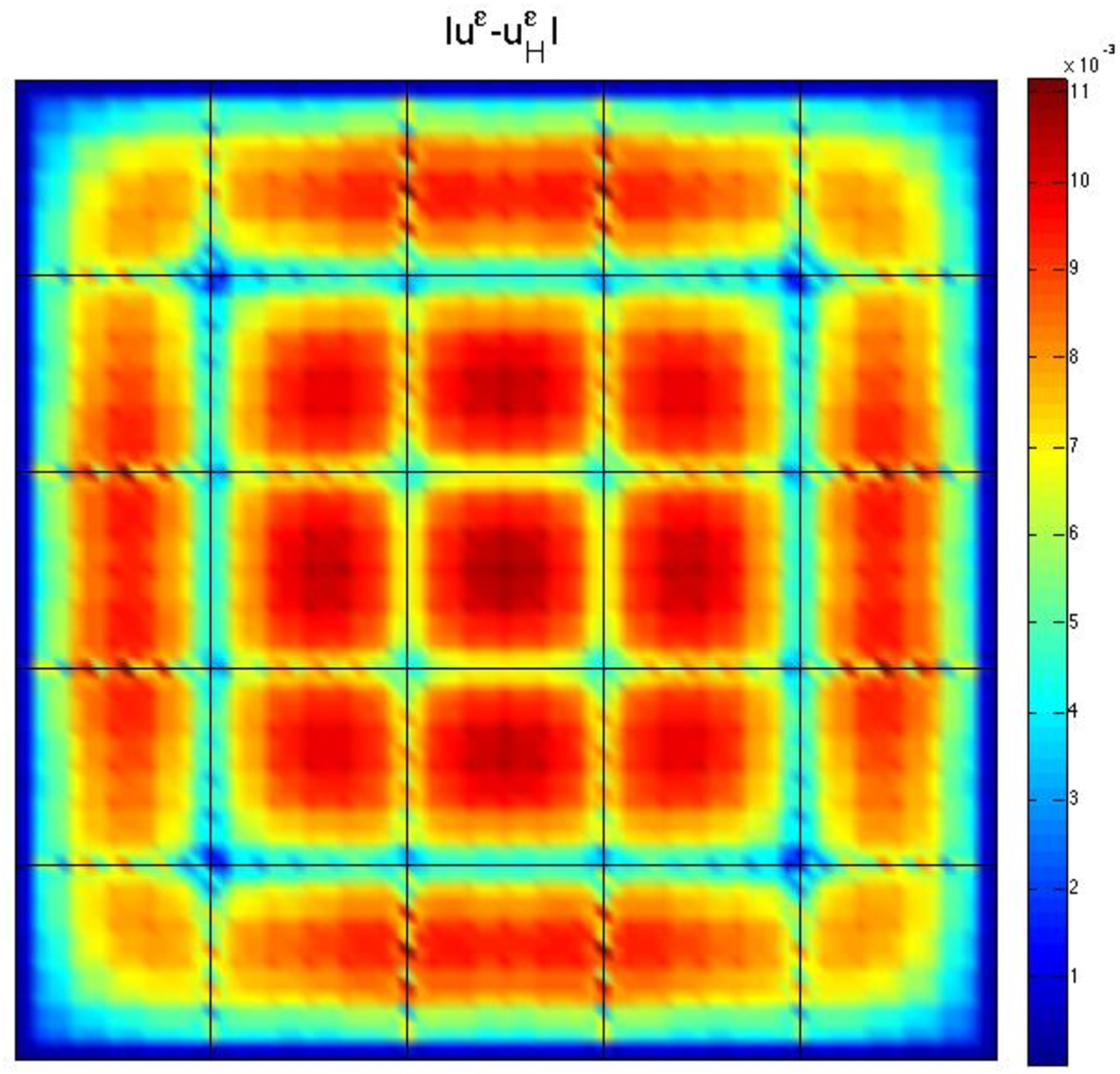} 
\includegraphics[width=75mm]{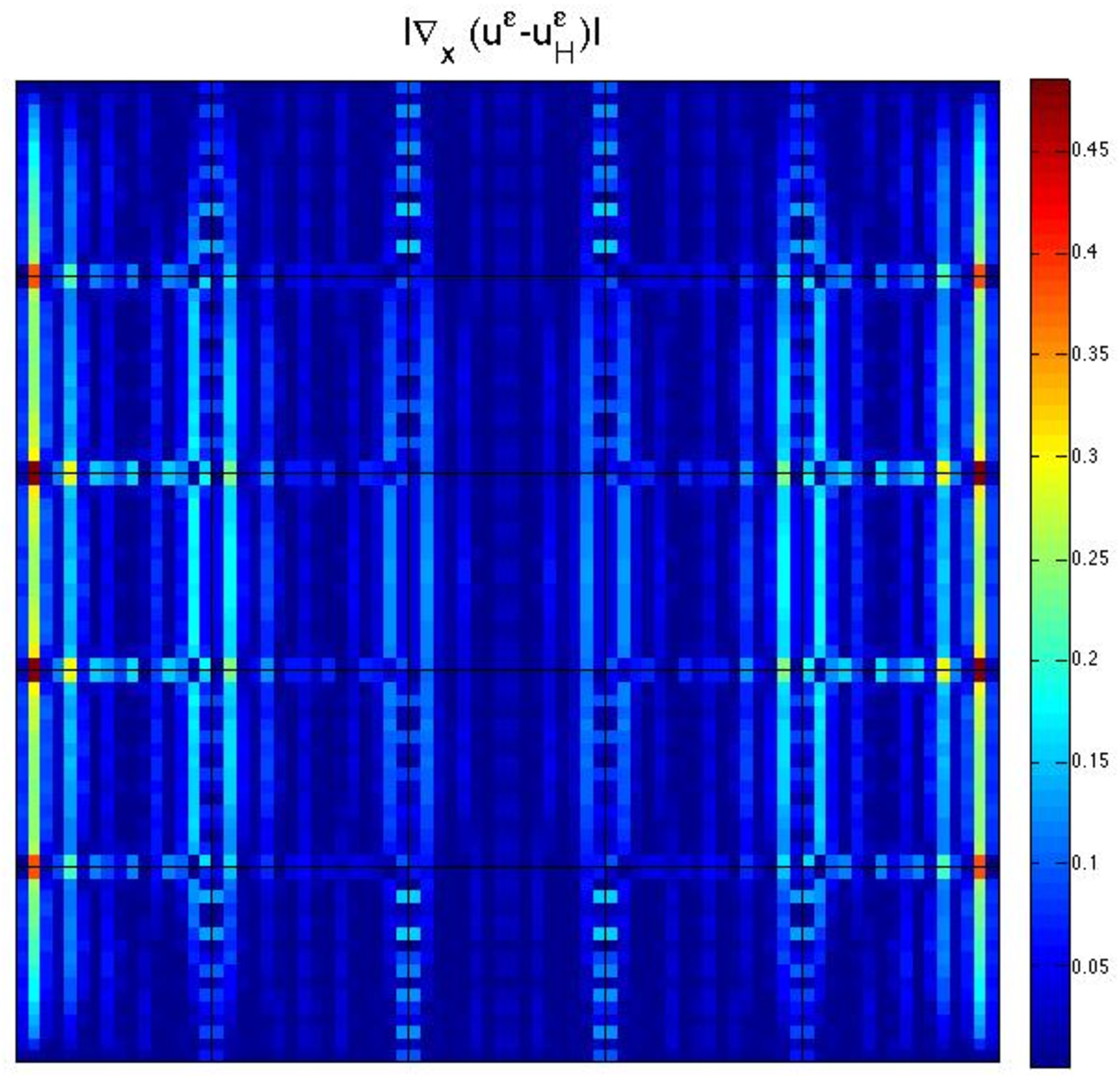} 
\end{center}
\caption{2D periodic example: plot of $|u^\eps-u_H^\eps|$ (left) and of $|\nabla(u^\eps-u_H^\eps) \cdot e_1|$ (right).}
\label{fig:distriberrorMsFEM2Dper}
\end{figure}

The value of the associated relative error estimate $\Delta_{MsFEM}/\vertiii{u^\eps_H}$ is about $33\%$, with an effectivity index of 1.08. Local contributions of $(\Delta_{MsFEM})^2 = \sum_K (\Delta^K_{MsFEM})^2$ over $\mT_H$, as well as local effectivity indices, are shown in Figure~\ref{fig:errorMsFEM2Dper}. The difference between the exact and the numerical solution, in terms of solution $u^\eps-u_H^\eps$ or gradient $|\nabla(u^\eps-u_H^\eps) \cdot e_1|$, is shown in Figure~\ref{fig:distriberrorMsFEM2Dper}. 

\begin{figure}[htbp]
\begin{center}
\includegraphics[width=150mm]{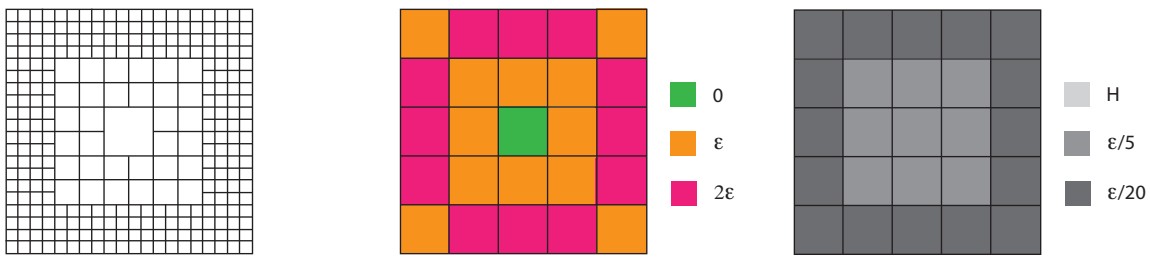}
\end{center}
\caption{2D periodic example: adapted coarse mesh (left), oversampling size (center) and fine mesh sizes $h_K$ (right).}
\label{fig:adapt2Dper}
\end{figure}

We now use the adaptive algorithm from this initial, non-accurate MsFEM solution (computed with $H_K = 0.2$, $h_K = H_K$ and $S_K=K$) and with a prescribed error tolerance of 5\%. This leads to the local MsFEM parameters shown in Figure~\ref{fig:adapt2Dper}. We observe that adapted discretization parameters (small coarse and fine mesh sizes, large oversampling) are required in regions located close to the boundary of $\Omega$.

Eventually, for the initial mesh $\mT_H$ made of $5\times 5$ macro elements and $h_K=\eps/5$, we analyze the effect of the value of the micro scale $\eps$ on the effectivity index. The results, given in Table~\ref{tab:effindex}, show that the quality of the error estimate is essentially independent of $\eps$. 

\begin{table}[htbp]
\centering
\begin{tabular}{| c | c c c | c c c | c c c | c c c |}
\hline
value of $\eps$ & & 0.04 & & & 0.01 & & & 0.004 & & & 0.001 & \\
\hline
effectivity index & & 1.08 & & & 1.12 & & & 1.10 & & & 1.15 & \\
\hline
\end{tabular}
\caption{2D periodic example: effectivity index for various values of $\eps$.}
\label{tab:effindex}
\end{table}

\subsection{2D example with a non-periodic coefficient}
\label{sec:num_2D_non-per}

We now consider a two-dimensional problem related to steady conduction through fiber composites. It is identical (up to the value of $\eps$) to the one discussed in~\cite{HOU97}. The problem is defined in the unit square domain $\Omega=(0,1)^2$ (unit cell of the composite material). Dirichlet boundary conditions $u(\bx)=(x_1-0.5)^2+(x_2-0.5)^2$ are applied on $\partial \Omega$, and a uniform load $f=-1$ is specified in $\Omega$. The conductivity properties of the medium are modeled by the tensor
\begin{equation*}
\Aa^\eps(\bx) = \Big[ 2 + P \cos\left(2\pi \tanh(w(r-0.3))/\eps\right) \Big] \, \II_2 = A^\eps(\bx) \, \II_2,
\end{equation*}
where $r=\sqrt{(x_1-0.5)^2+(x_2-0.5)^2}$ is the distance from the center of the fiber. The parameter $P$ controls the ratio between the conductivity of the fibers and that of the matrix, the parameter $w$ determines the total width of the reinforcement, and $\eps$ sets the wavelength of the unidirectional oscillations. In the following, we set $P=1.8$, $w=20$ and $\eps=0.2$, so that the shortest wavelength in the oscillations is about $\eps_0=0.01$. 

\begin{figure}[htbp]
\begin{center}
\includegraphics[width=75mm]{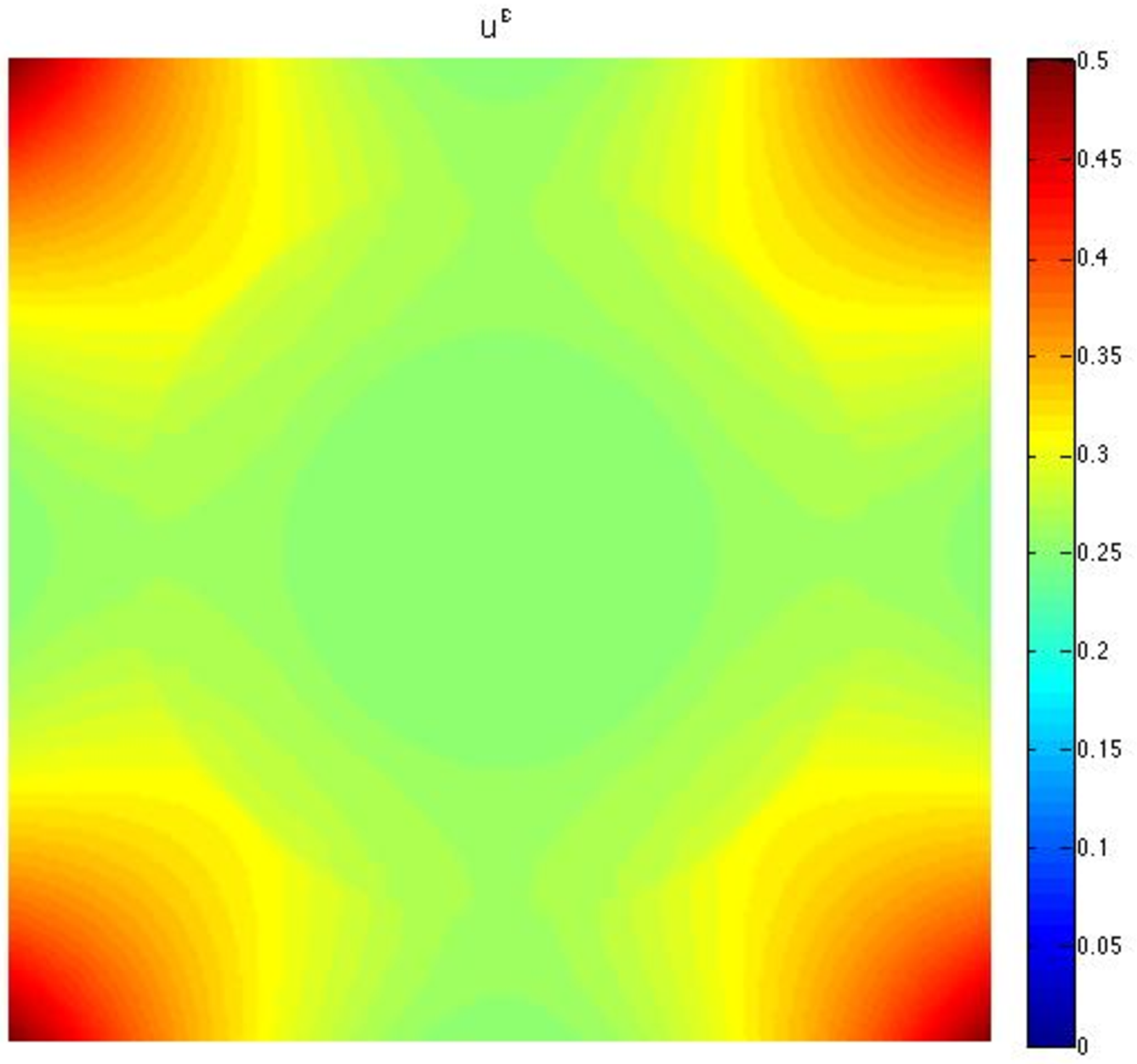}
\includegraphics[width=75mm]{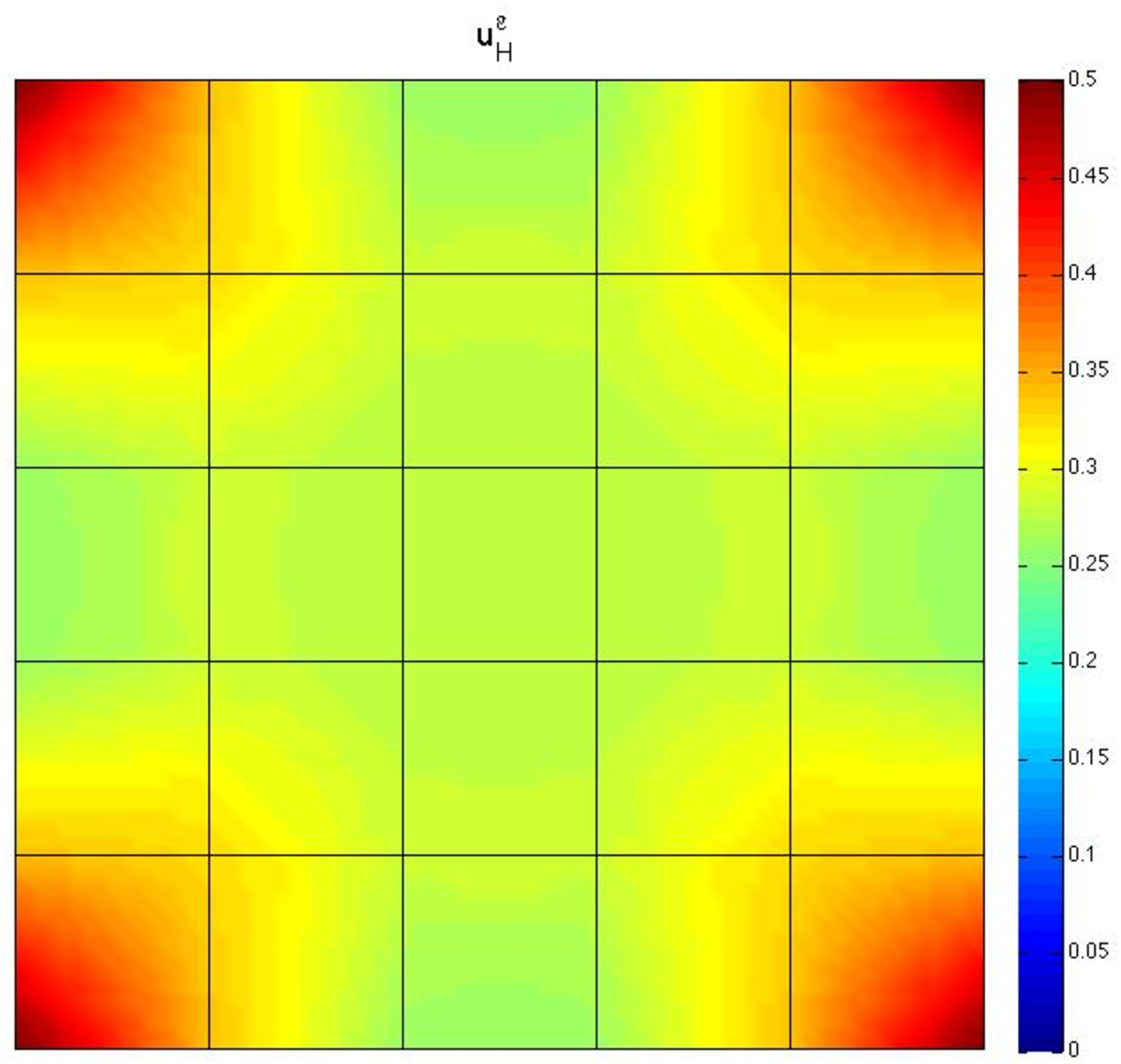}
\end{center}
\caption{2D non-periodic example: exact solution $u^\eps$ (left) and MsFEM solution $u^\eps_H$ (right).}
\label{fig:solMsFEM2Dfiberu}
\end{figure}
\begin{figure}[htbp]
\begin{center}
\includegraphics[width=75mm]{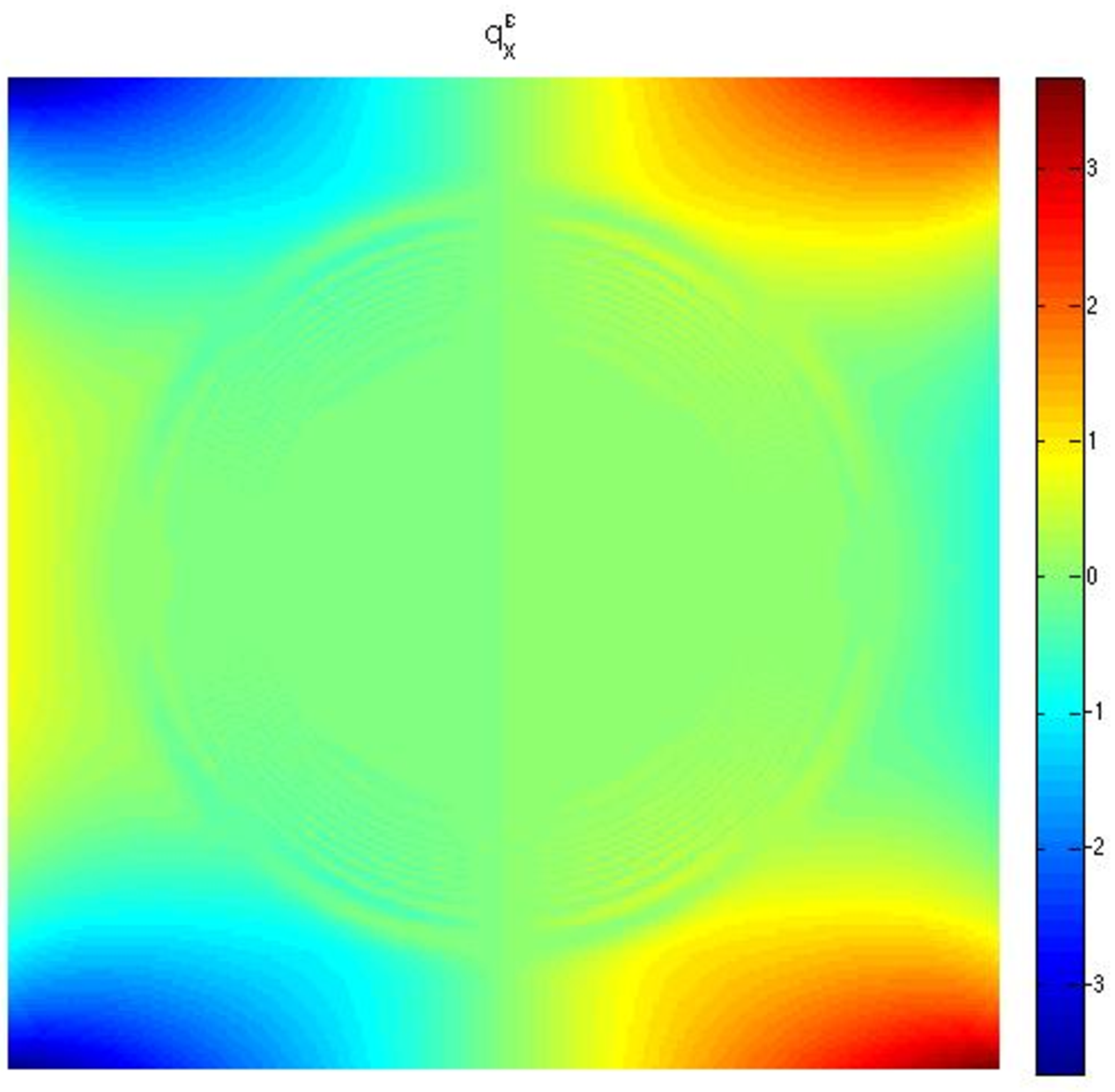}
\includegraphics[width=75mm]{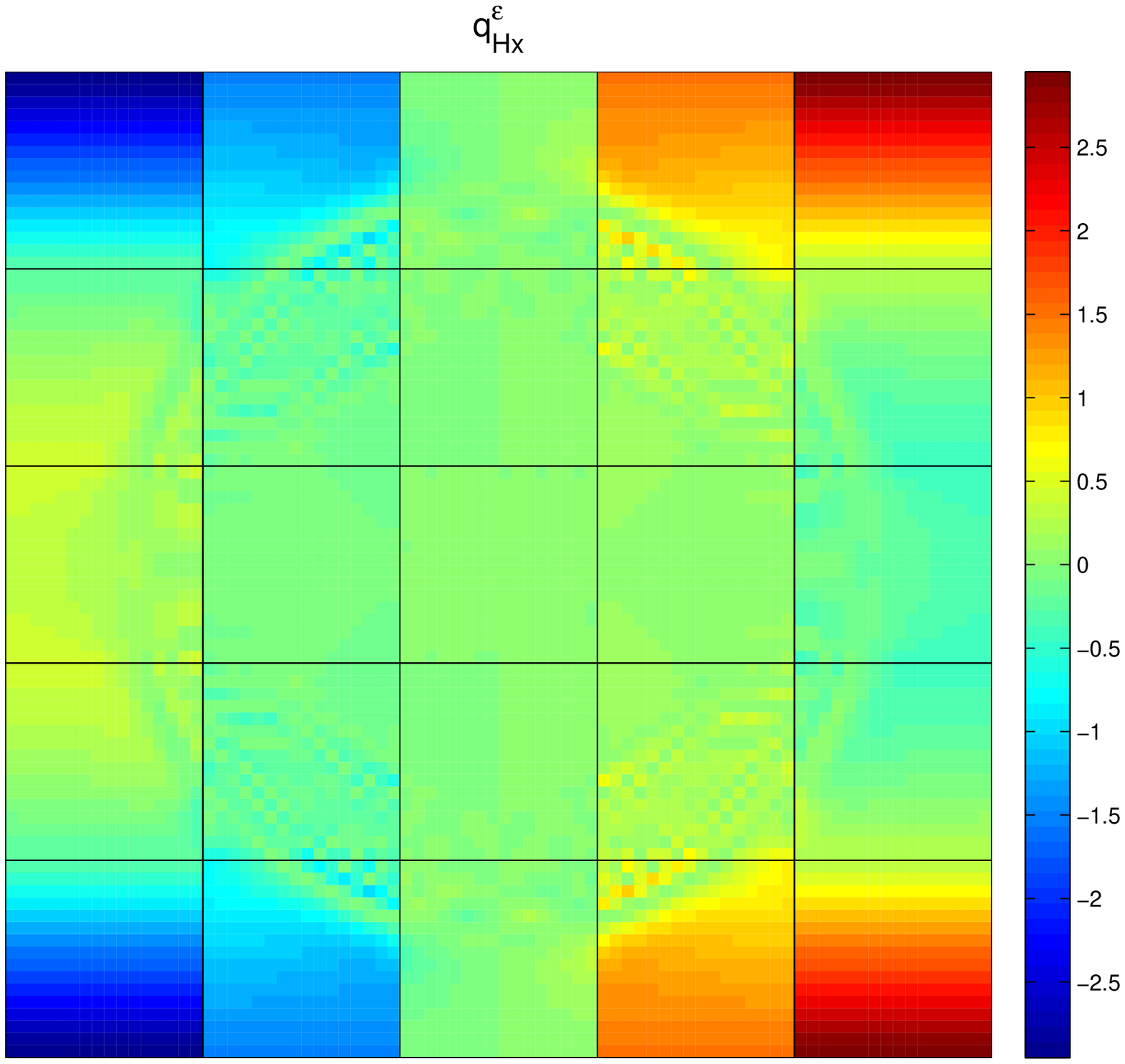}
\end{center}
\caption{2D non-periodic example: exact flux $\bq^\eps$ (left) and MsFEM flux $\bq^\eps_H$ (right) (we only show the first component of the flux; the figure for the second component is similar).}
\label{fig:solMsFEM2Dfiberflux}
\end{figure}

We take an initial coarse mesh $\mT_H$ made of $5\times 5$ macro elements. Choosing $h_K=\eps_0/5$ for any $K$, and with no oversampling, the approximate MsFEM solution is represented in the Figures~\ref{fig:solMsFEM2Dfiberu} and~\ref{fig:solMsFEM2Dfiberflux} and compared to the exact solution (computed using a $500 \times 500$ fine mesh): solutions $u^\eps$ and $u^\eps_H$ are shown in Figure~\ref{fig:solMsFEM2Dfiberu} and fluxes $\bq^\eps$ and $\bq^\eps_H$ are shown in Figure~\ref{fig:solMsFEM2Dfiberflux}.

\begin{figure}[htbp]
\begin{center}
\includegraphics[width=75mm]{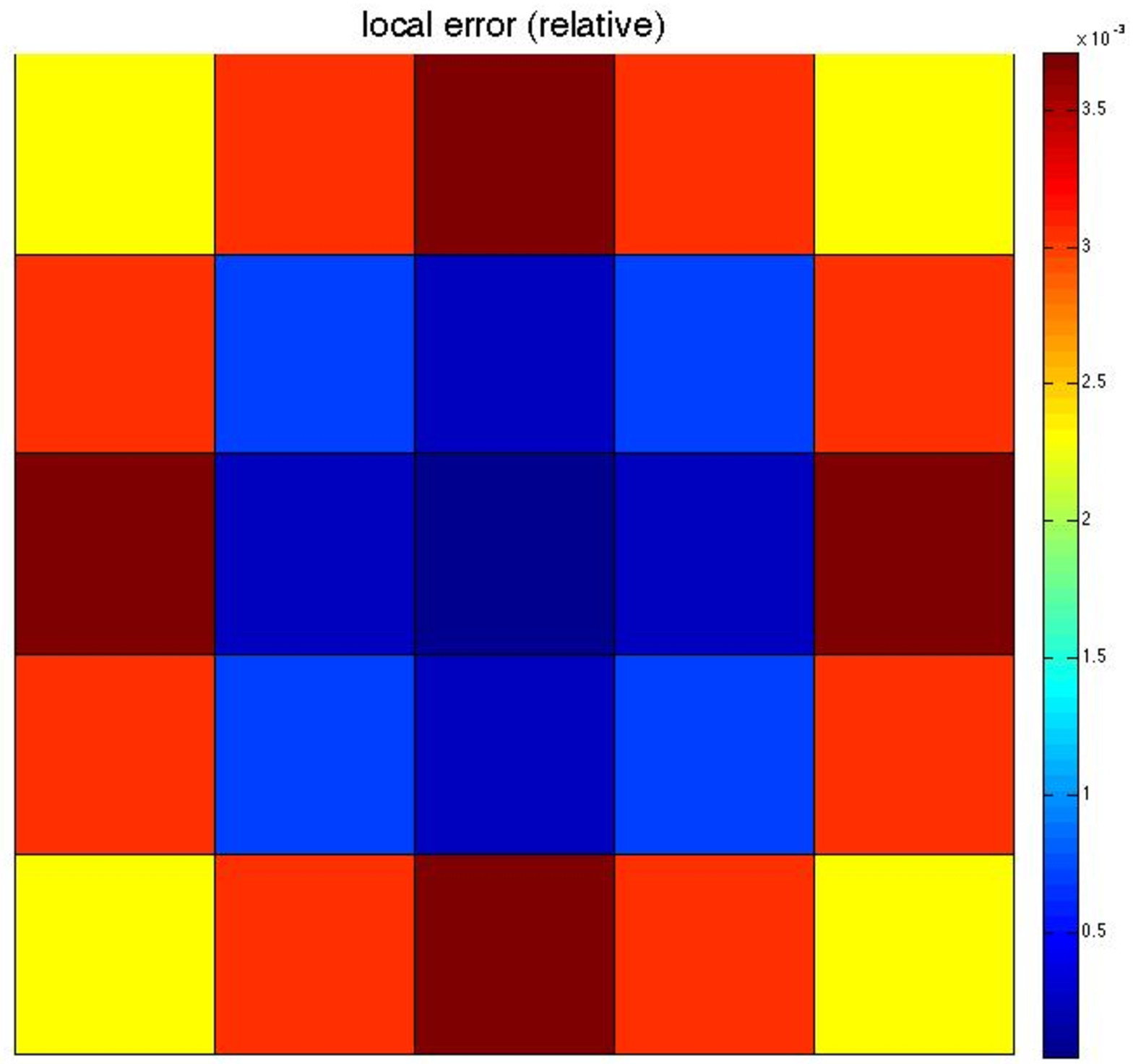} 
\includegraphics[width=75mm]{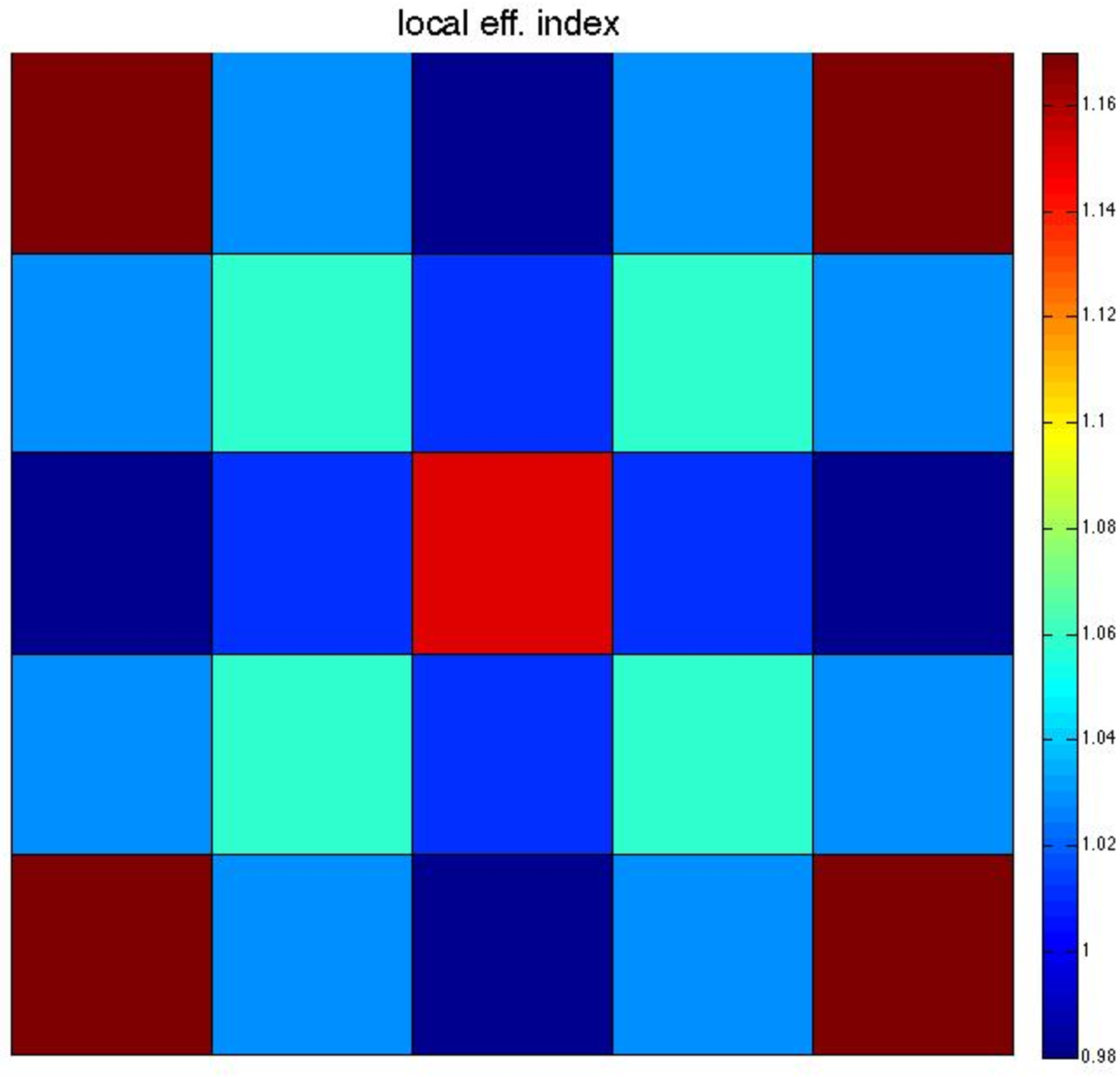} 
\end{center}
\caption{2D non-periodic example: values of $(\Delta^K_{MsFEM})^2$ for the specific computed MsFEM solution (left), and local effectivity indices (right). }
\label{fig:errorMsFEM2Dfiber}
\end{figure}

The value of the associated relative error estimate $\Delta_{MsFEM}/\vertiii{u^\eps_H}$ is about $36\%$, with an effectivity index of 1.09. Local contributions of $(\Delta_{MsFEM})^2$ over $\mT_H$, as well as local effectivity indices, are shown in Figure~\ref{fig:errorMsFEM2Dfiber}.

\begin{figure}[htbp]
\begin{center}
\includegraphics[width=150mm]{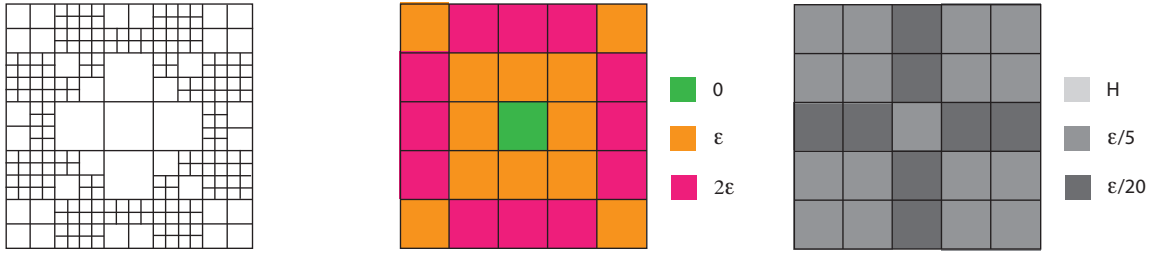}
\end{center}
\caption{2D non-periodic example: adapted coarse mesh (left), oversampling size (center) and fine mesh sizes $h_K$ (right).}
\label{fig:adapt2Dfiber}
\end{figure}

We now use the adaptive algorithm from this initial, non-accurate MsFEM solution ($H_K=0.2$, $h_K=H_K$ and $S_K=K$) and with a prescribed error tolerance of 5\%. This leads to the local MsFEM parameters shown in Figure~\ref{fig:adapt2Dfiber}, associated with an estimated error of 4.77\%. We observe that this error tolerance can be reached without resorting to accurate fine-scale computations in the center of the domain $\Omega$. Our adaptive algorithm leads to a discretization which is much more efficient than the one that would have been obtained using a uniform refinement. 

\subsection{A 2D structure with a crack}
\label{sec:num_2D_crack}

As a last example, we consider a problem over the cracked domain $\Omega = (0,1)^2 \setminus \{(x_1,0), \ \ 0.5 \le x_1 \le 1 \}$, with the loading $f=0$. It is identical (up to the value of $\eps$) to the one discussed in~\cite{ABD09}. On the boundary $\partial \Omega$ (including on the lips of the crack), we prescribe a Dirichlet boundary condition which matches the exact homogenized solution.


The diffusion tensor is chosen as
\begin{equation*}
\Aa^\eps(\bx)=\frac{64}{9 \sqrt{17}} \left[ \sin \left(\frac{2\pi (x_1-0.5)}{\eps}\right) + \frac{9}{8} \right] \left[ \cos \left( \frac{2\pi (x_2-0.5)}{\eps}\right) + \frac{9}{8} \right] \, \II_2 = A^\eps(\bx) \ \II_2,
\end{equation*}
with $\eps=0.04$. It is periodic, and the homogenized tensor $\Aa^0$ is equal to the identity matrix: $\Aa^0 = \II_2$. 


The homogenized solution $u^0$ is not in $\mH^2(\Omega)$. When using standard piecewise linear finite elements, the numerical error between $u^0$ and $u^0_H$ thus decreases with a sub-optimal convergence rate. One possibility to address this issue is to adaptively refine the mesh close to the crack tip. Within the MsFEM approach, we intend to have this same adaptive procedure for the macro  mesh $\mT_H$, all the more so as traditional scale separation assumptions no longer hold in the vicinity of the crack tip (high variation of the gradients).

\medskip

\begin{figure}[htbp]
\begin{center}
\includegraphics[width=75mm]{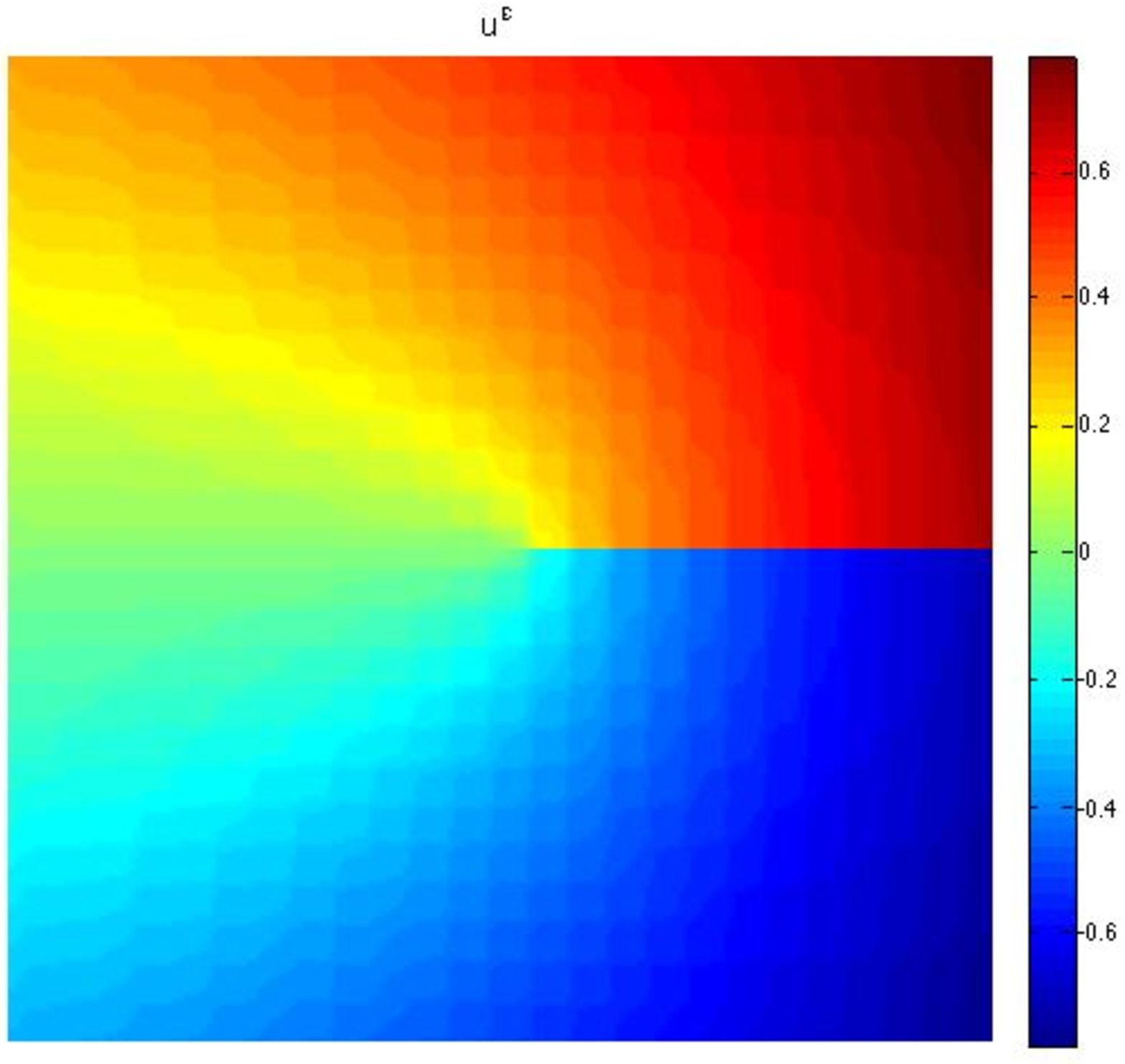}
\includegraphics[width=75mm]{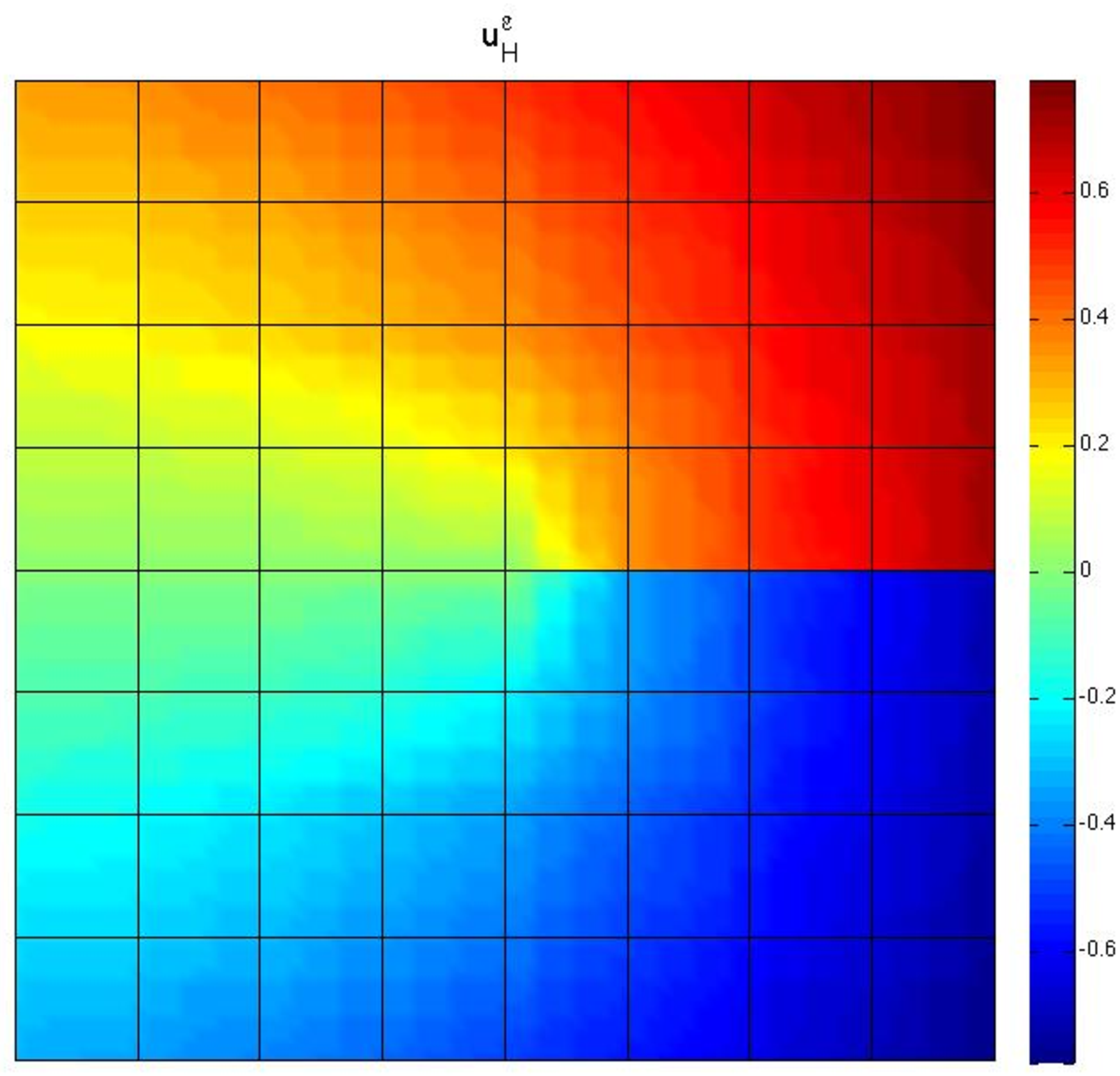}
\end{center}
\caption{2D example with a crack: exact solution $u^\eps$ (left) and MsFEM solution $u^\eps_H$ (right).}
\label{fig:solMsFEM2Dcracku}
\end{figure}

\begin{figure}[htbp]
\begin{center}
\includegraphics[width=75mm]{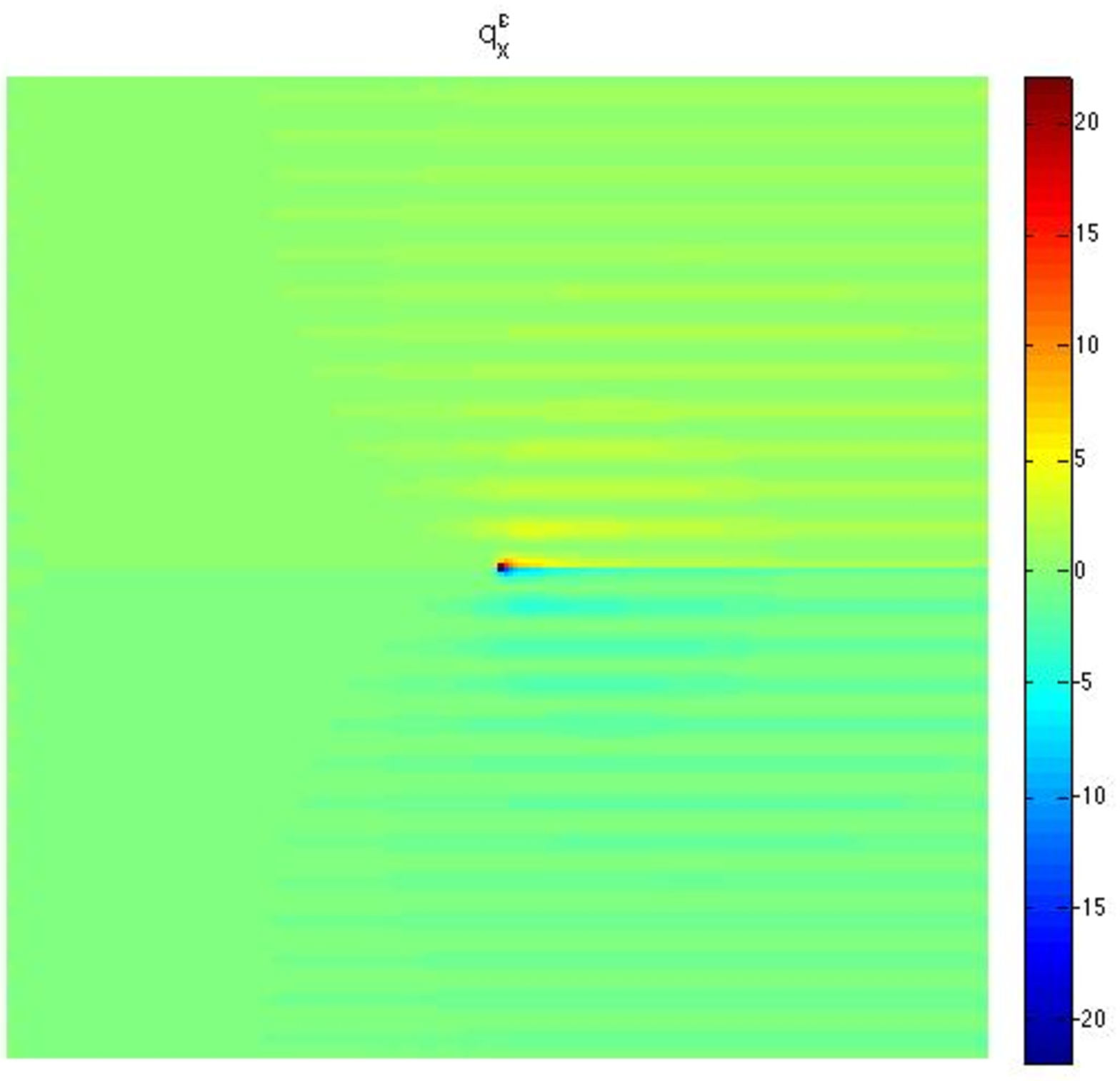}
\includegraphics[width=75mm]{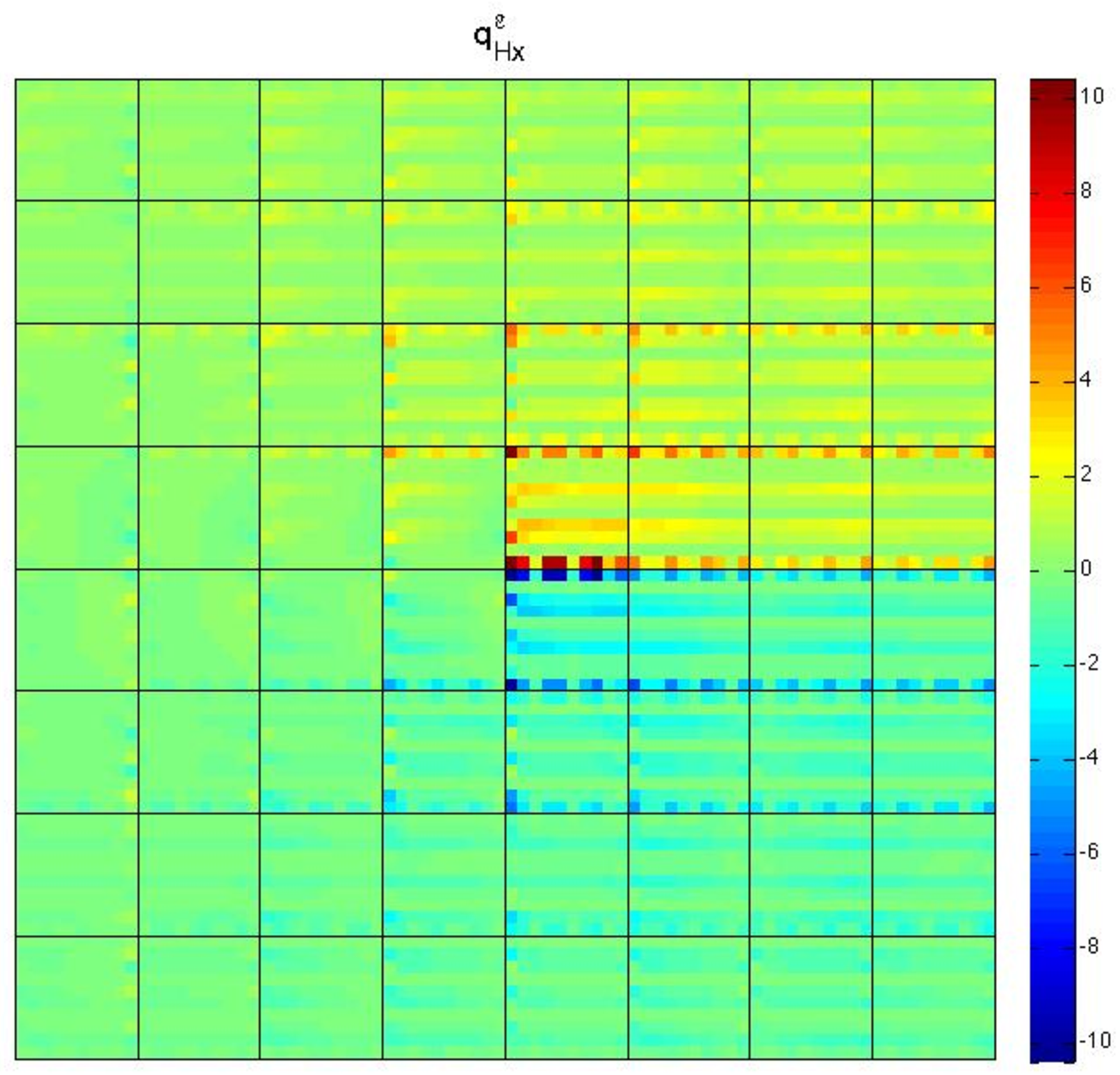}
\end{center}
\caption{2D example with a crack: exact flux $\bq^\eps$ (left) and MsFEM flux $\bq^\eps_H$ (right)  (we only show the first component of the flux; the figure for the second component is similar).}
\label{fig:solMsFEM2Dcrackflux}
\end{figure}

We consider an initial coarse mesh $\mT_H$ made of $8\times 8$ macro elements. Choosing $h_K=\eps/3$ for any $K$, and with no oversampling, the approximate MsFEM solution is represented on Figures~\ref{fig:solMsFEM2Dcracku} and~\ref{fig:solMsFEM2Dcrackflux} and compared to the exact solution (obtained using a $500 \times 500$ fine mesh): solutions $u^\eps$ and $u^\eps_H$ are shown in Figure~\ref{fig:solMsFEM2Dcracku}, and fluxes $\bq^\eps$ and $\bq^\eps_H$ are shown in Figure~\ref{fig:solMsFEM2Dcrackflux}.

\begin{figure}[htbp]
\begin{center}
\includegraphics[width=75mm]{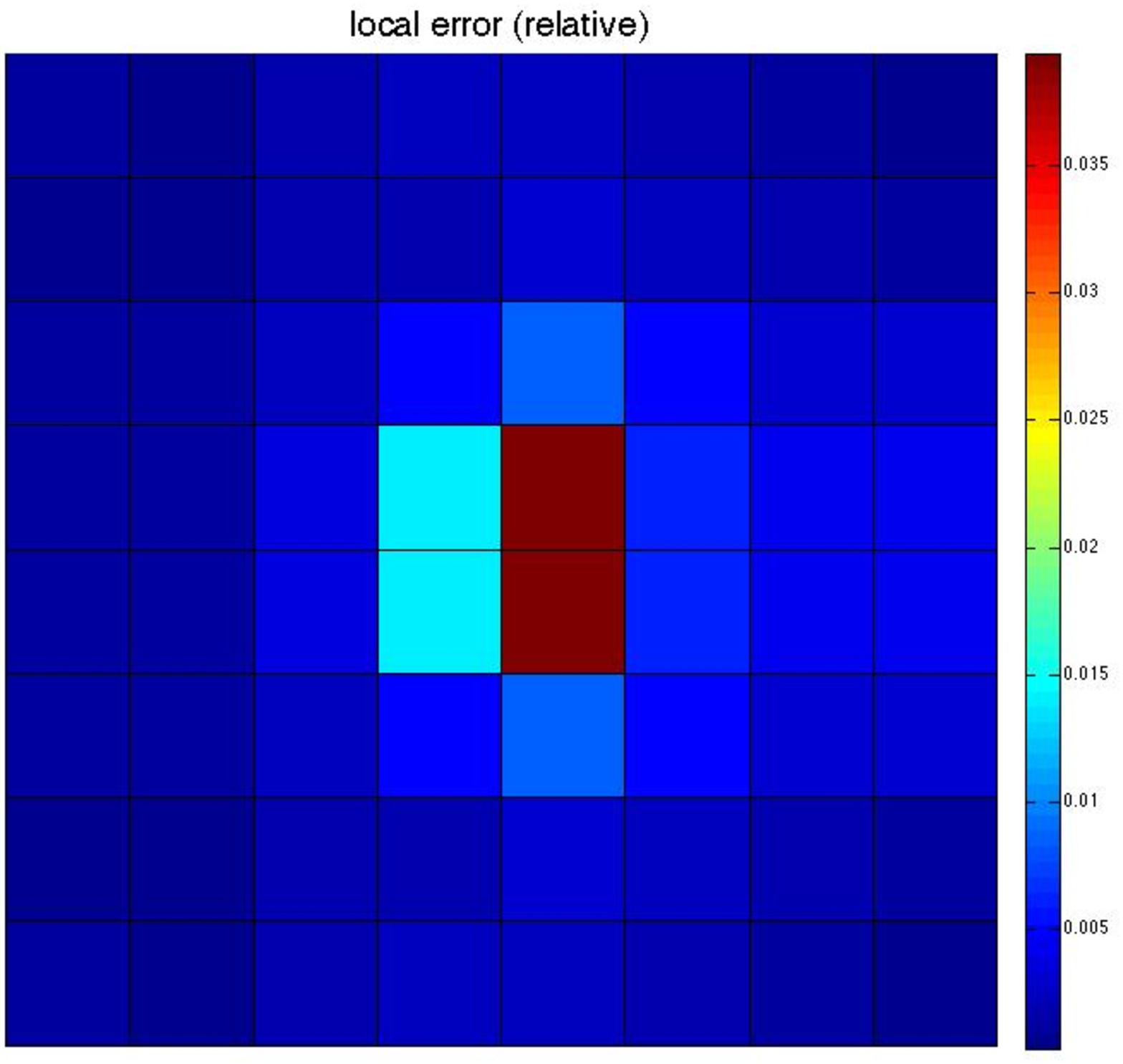} 
\includegraphics[width=75mm]{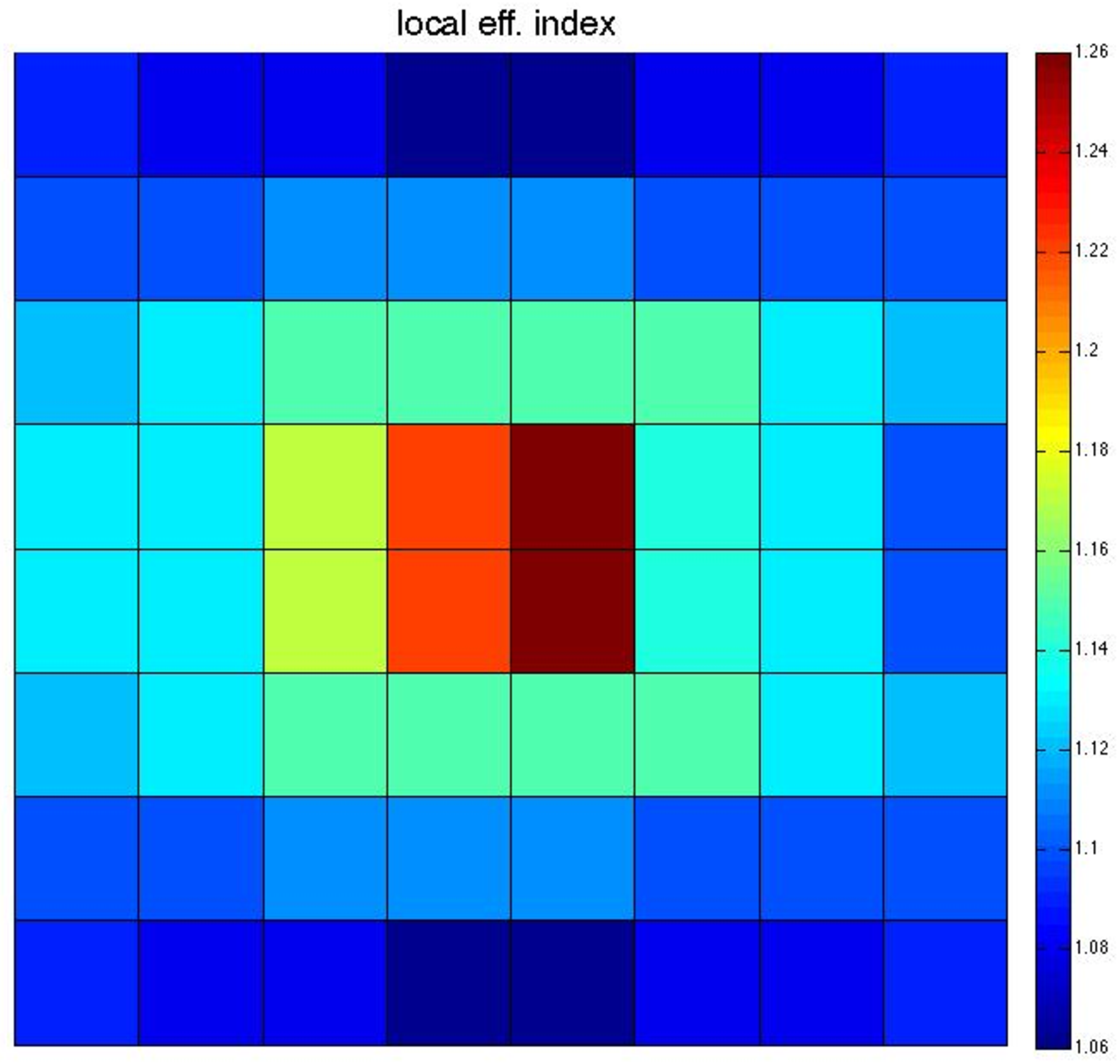} 
\end{center}
\caption{2D example with a crack: values of $(\Delta^K_{MsFEM})^2$ for the specific computed MsFEM solution (left), and local effectivity indices (right). }
\label{fig:errorMsFEM2Dcrack}
\end{figure}

The value of the associated relative error estimate $\Delta_{MsFEM}/\vertiii{u^\eps_H}$ is about $47\%$, with an effectivity index of 1.18. Local contributions of $(\Delta_{MsFEM})^2$ over $\mT_H$, as well as local effectivity indices, are shown in Figure~\ref{fig:errorMsFEM2Dcrack}. 

\begin{figure}[htbp]
\begin{center}
\includegraphics[width=150mm]{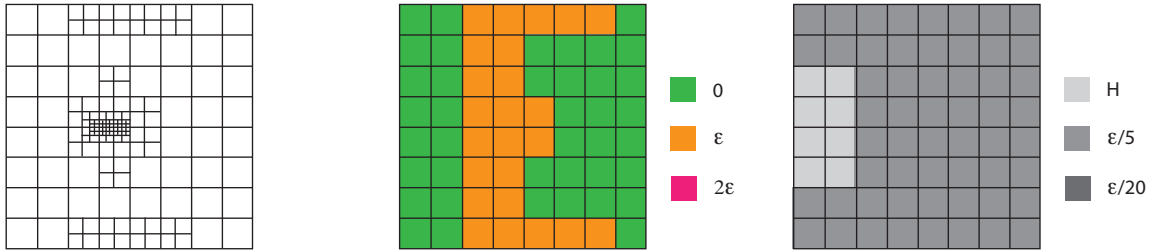}
\end{center}
\caption{2D example with a crack: adapted coarse mesh (left), oversampling size (center) and fine mesh sizes $h_K$ (right).}
\label{fig:adapt2Dcrack}
\end{figure}

Starting from the initial, non-accurate MsFEM approximation, we now use our adaptive algorithm, with a prescribed error tolerance of 5\%. This leads to the local MsFEM parameters shown in Figure~\ref{fig:adapt2Dcrack}. We observe that the coarse mesh needs to be highly refined in the vicinity of the crack tip. In contrast, the oversampling layer never has to be larger than $\eps$.

For this example, it would of course be interesting to use a MsFEM version of the XFEM method~\cite{MOE99}, which would consists in defining the MsFEM basis functions by means of local problems (posed on each coarse element $K$) complemented by {\em enriched boundary conditions}, in the spirit of the XFEM method. We did not pursue in that direction, which goes beyond the scope of this work. 

\section{Conclusions and prospects}\label{section:conclusions}

We have introduced numerical tools for \textit{a posteriori} error estimation and adaptive strategies in MsFEM computations. They lead to fully computable and guaranteed upper bounds on the error measured in the energy norm, and enable to perform adaptive multiscale computations by selecting the relevant MsFEM parameters to reach a prescribed error tolerance. An important feature of these numerical tools is that they are compatible with the {\em offline}/{\em online} and element-wise procedure of MsFEM. The additional computational costs are thus affordable. Using this procedure, the fine scale features of the solution are adaptively recovered during the online enrichment phase, if and where needed.

Further studies, investigated in future works as extensions to the present article, will be conducted towards the establishement of goal-oriented error estimators, allowing to control the error on specific outputs of the problem~\cite{Chamoin_Legoll_goal}. Possible extensions also include the use of the CRE concept for other multiscale problems solved with MsFEM, such as non-symmetric problems (advection-diffusion-reaction), hyperbolic problems (elastodynamics), stochastic problems~\cite{XU07,LEB14b} or some nonlinear problems~\cite{EFE04}.

\section*{Acknowledgements}

We thank C. Le Bris for stimulating discussions on the topics addressed in this article. LC thanks Inria for enabling his two-year leave (2014--2016) in the MATHERIALS project team. The work of FL is partially supported by ONR under grant N00014-15-1-2777 and by EOARD under grant FA9550-17-1-0294.


\begin{thebibliography}{00}

\bibitem{BAB94}
Babuska I., Caloz G., Osborn J.E.
\newblock Special finite element methods for a class of second order elliptic problems with rough coefficients.
\newblock {\em SIAM Journal on Numerical Analysis} 1994; \textbf{31}(4):945--981.

\bibitem{BAB83}
Babuska I., Osborn J.E.
\newblock Generalized finite element methods: their performance and their relation to mixed methods.
\newblock {\em SIAM Journal on Numerical Analysis} 1983; \textbf{20}:510--536. 

\bibitem{STR01}
Strouboulis T., Copps K., Babuska I.
\newblock The generalized finite element method.
\newblock {\em Computer Methods in Applied Mechanics and Engineering} 2001; \textbf{190}:4081--4193. 

\bibitem{FIS04} 
Fish J., Yuan Z. 
\newblock Multiscale enrichment based on partition of unity. 
\newblock {\em International Journal for Numerical Methods in Engineering} 2004; \textbf{70}:1341-1359.

\bibitem{HOU97}
Hou T., Wu X.H.
\newblock A multiscale finite element method for elliptic problems in composite materials and porous media.
\newblock {\em Journal of Computational Physics} 1997; \textbf{134}:169--189. 

\bibitem{HOU99} 
Hou T., Wu X.H., Cai Z.
\newblock Convergence of a multiscale finite element method for elliptic problems with rapidly oscillating coefficients.
\newblock {\em Mathematics of Computation} 1999; \textbf{68}(227):913--943. 

\bibitem{EFE09}
Efendiev Y., Hou T.
\newblock {\em Multiscale Finite Element Methods: Theory and Applications}.
\newblock Springer, New York 2009. 

\bibitem{HUG98}
Hughes T.J.R., Feijoo G.R., Mazzei L., Quincy J.B.
\newblock The variational multiscale method -- A paradigm for computational mechanics.
\newblock {\em Computer Methods in Applied Mechanics and Engineering} 1998; \textbf{166}(1-2):3--24. 

\bibitem{HOA05}
Hoang V.H., Schwab C.
\newblock High-dimensional finite elements for elliptic problems with multiple scales.
\newblock {\em SIAM Multiscale Modeling \& Simulation} 2005; \textbf{3}:168--194. 

\bibitem{HAC85}
Hackbusch W.
\newblock {\em Multi-Grid Methods and Applications}.
\newblock Springer Series in Computational Mathematics 1985. 

\bibitem{E03}
E W., Engquist B., Huang Z.
\newblock Heterogeneous multiscale method: a general methodology for multiscale modeling.
\newblock {\em Physical Review B} 2003; \textbf{67}(9):092101. 

\bibitem{ABD12} 
Abdulle A., E W., Engquist B., Vanden-Eijnden E.
\newblock The heterogeneous multiscale method.
\newblock {\em Acta Numerica} 2012; \textbf{21}:1--87. 

\bibitem{FEY03} 
Feyel F.
\newblock A multilevel finite element method (FE2) to describe the response of highly non-linear structures using generalized continua.
\newblock {\em Computer Methods in Applied Mechanics and Engineering} 2003; \textbf{192}:3233--3244. 

\bibitem{SUQ87} 
Suquet P.
\newblock Elements of homogenization theory for inelastic solid mechanics.
\newblock in {\em Homogenization techniques for composite media}, E. Sanchez-Palencia and A. Zaoui, eds., Lecture Notes in Physics, vol. 272, Springer 1987 (pp.~194--278). 

\bibitem{LEB14} 
Le Bris C., Legoll F., Lozinski A. 
\newblock An MsFEM type approach for perforated domains. 
\newblock {\em SIAM Multiscale Modeling \& Simulation} 2014; \textbf{12}(3):1046--1077.

\bibitem{VER96}
Verf\"urth R.
\newblock {\em A review of a posteriori error estimates and adaptive mesh-refinement techniques}.
\newblock Wiley-Teubner 1996. 

\bibitem{AIN00}
Ainsworth M., Oden J.T.
\newblock {\em A posteriori error estimation in finite element analysis}.
\newblock John Wiley \& Sons 2000. 

\bibitem{LAD04} 
Ladev\`eze P., Pelle J.-P. 
\newblock {\em Mastering Calculations in Linear and Nonlinear Mechanics}. 
\newblock Springer, 2004.

\bibitem{CHA15}
Chamoin L., D\'iez P. (ed).
\newblock {\em Verifying calculations, forty years on: an overview of classical verification techniques for FEM simulations}.
\newblock SpringerBriefs, 2015.

\bibitem{FRA72} 
Fraeijs de Veubeke B., Hugge M.A. 
\newblock Dual analysis for heat conduction problems by finite elements. 
\newblock {\em International Journal for Numerical Methods in Engineering} 1972; \textbf{5}(1):65-82.

\bibitem{ODE74} 
Oden J.T., Reddy J.N. 
\newblock On dual complementary variational principles in mathematical physics. 
\newblock {\em International Journal of Engineering Science} 1974; \textbf{12}:1-29.

\bibitem{FRA01} 
Fraeijs de Veubeke B.
\newblock Displacement and equilibrium models in the finite element method. 
\newblock {\em International Journal for Numerical Methods in Engineering, Classical Reprint Series} 2001; \textbf{52}:287-342.

\bibitem{PLE11} %
Pled F., Chamoin L., Ladev\`eze P. 
\newblock On the techniques for constructing admissible stress fields in model verification: performances on engineering examples. 
\newblock {\em International Journal for Numerical Methods in Engineering} 2011; \textbf{88}(5):409-441.

\bibitem{LAD10} %
Ladev\`eze P., Chamoin L. 
\newblock Calculation of strict error bounds for finite element approximations of nonlinear pointwise quantities of interest. 
\newblock {\em International Journal for Numerical Methods in Engineering} 2010; \textbf{84}:1638--1664.

\bibitem{EFE00}
Efendiev Y., Hou T., Wu X.H.
\newblock Convergence of a nonconforming multiscale finite element method.
\newblock {\em SIAM Journal of Numerical Analysis} 2000; \textbf{37}(3):888--910. 

\bibitem{HOU04}
Hou T., Wu X.H., Zhang Y.
\newblock Removing the cell resonance error in the multiscale finite element method via a Petrov-Galerkin formulation.
\newblock {\em Communications in Mathematical Sciences} 2004; \textbf{2}:185--205.

\bibitem{ALL06}
Allaire G., Brizzi R.
\newblock A multiscale finite element method for numerical homogenization.
\newblock {\em SIAM Multiscale Modeling \& Simulation} 2006; \textbf{4}(3):790--812.  

\bibitem{HES14}
Hesthaven J., Zhang S., Zhu X.
\newblock High-order multiscale finite element method for elliptic problems.
\newblock {\em SIAM Multiscale Modeling \& Simulation} 2014; \textbf{12}(2):650--666. 

\bibitem{STR06}
Strouboulis T., Zhang L., Wang D., Babuska I.
\newblock A posteriori error estimation for generalized finite element methods.
\newblock {\em Computer Methods in Applied Mechanics and Engineering} 2006; \textbf{195}:852--879. 

\bibitem{STR07} 
Strouboulis T., Zhang L., Babuska I.
\newblock Assessment of the cost and accuracy of the generalized FEM.
\newblock {\em International Journal for Numerical Methods in Engineering} 2007; \textbf{69}:250--283. 

\bibitem{LAR07}
Larson M.G., Malqvist A.
\newblock Adaptive variational multiscale methods based on a posteriori error estimation: energy norm estimates for elliptic problems.
\newblock {\em Computer Methods in Applied Mechanics and Engineering} 2007; \textbf{196}(21-24):2313--2324. 

\bibitem{NOL08} 
Nolen J., Papanicolaou G., Pironneau O.
\newblock A framework for adaptive multiscale methods for elliptic problems.
\newblock {\em SIAM Multiscale Modeling \& Simulation} 2008; \textbf{7}(1):171--196. 

\bibitem{ABD09} 
Abdulle A., Nonnenmacher A.
\newblock A posteriori error analysis of the heterogeneous multiscale method for homogenization problems.
\newblock {\em C.R. Acad. Sci. Paris, S\'erie I} 2009; \textbf{347}:1081--1086. 

\bibitem{ABD11a} 
Abdulle A., Nonnenmacher A.
\newblock Adaptive finite element heterogeneous multiscale method for homogenization problems.
\newblock {\em Computer Methods in Applied Mechanics and Engineering} 2011; \textbf{200}:2710--2726. 

\bibitem{LAR11}
Larsson F., Runesson K.
\newblock On two-scale adaptive FE analysis of micro-heterogeneous media with seamless scale-bridging.
\newblock {\em Computer Methods in Applied Mechanics and Engineering} 2011; \textbf{200}:2662--2674. 

\bibitem{PAL17}
Paladim D., Moitinho de Almeida J.P., Bordas S., Kerfriden P.
\newblock Guaranteed error bounds in homogenization: an optimum stochastic approach to preserve the numerical separation of scales.
\newblock {\em International Journal for Numerical Methods in Engineering} 2017; \textbf{110}(2):103--132. 

\bibitem{HEN14}
Henning P., Ohlberger M., Schweizer B.
\newblock An adaptive multiscale finite element method.
\newblock {\em SIAM Multiscale Modeling \& Simulation} 2014; \textbf{12}(3):1078--1107. 

\bibitem{CHU14}
Chung E.T., Efendiev Y., Li G.
\newblock An adaptive GMsFEM for high contrast flow problems.
\newblock {\em Journal of Computational Physics} 2014; \textbf{273}:54--76. 

\bibitem{ERN10}
Ern A., Stephansen A.F., Vohralik M.
\newblock Guaranteed and robust discontinuous Galerkin a posteriori error estimates for convection-diffusion-reaction problems.
\newblock {\em Journal of Computational and Applied Mathematics} 2010; \textbf{234}(1):114--130. 

\bibitem{CHA16} 
Chamoin L., Legoll F. 
\newblock A pedagogical view on a posteriori error estimation in finite element analysis.
\newblock {\em Submitted to Applied Mathematics Research Express} 2017.

\bibitem{BEN78}
Bensoussan A., Lions J.-L., Papanicolaou G.
\newblock {\em Asymptotic Analysis for Periodic Structures}.
\newblock Studies in Mathematics and its Applications, vol. 5, North-Holland, Amsterdam, New York 1978.

\bibitem{SAN80}
Sanchez-Palencia E.
\newblock {\em Non Homogeneous Media and Vibration Theory}.
\newblock Springer, Heidelberg 1980.

\bibitem{JIK94}
Jikov V.V., Kozlov S.M., Oleinik, O.A.
\newblock {\em Homogenization of differential operators and integral functionals}.
\newblock Springer-Verlag, Berlin Heidelberg 1994.

\bibitem{MUR97} 
Murat F., Tartar L. 
\newblock Calculus of variations and homogenization. 
\newblock in {\em Topics in the Mathematical Modelling of Composite Materials}, Progress in nonlinear differential equations and their applications, vol. 31, Birkhauser 1997 (pp.~139--173).

\bibitem{TAR10}
Tartar, L.
\newblock {\em The general theory of homogenization - A personalized introduction}.
\newblock Lecture Notes of the Unione Matematica Italiana, vol. 7, Springer-Verlag, Berlin Heidelberg 2010.

\bibitem{KAN09}
Kanout\'e P., Boso D.P., Chaboche J.-L., Schrefler B.A.
\newblock Multiscale methods for composites: a review.
\newblock {\em Archives of Computational Methods in Engineering} 2009; \textbf{16}(1):31--75. 

\bibitem{CIA78}
Ciarlet P.G.
\newblock {\em The Finite Element Method for Elliptic Problems}.
\newblock North Holland 1978. 

\bibitem{STR77}
Strang G., Fix G.J.
\newblock {\em An Analysis of the Finite Element Method}.
\newblock Prentice Hall 1977. 

\bibitem{HEN13}
Henning P., Peterseim D.
\newblock Oversampling for the multiscale finite element method.
\newblock {\em SIAM Multiscale Modeling \& Simulation} 2013; \textbf{11}(4):1149--1175. 

\bibitem{LAD83}
Ladev\`eze P., Leguillon D.
\newblock Error estimate procedure in the finite element method and applications. 
\newblock {\em SIAM Journal of Numerical Analysis} 1983; \textbf{20}(3):485--509.

\bibitem{DES99}
Destuynder P., M\'etivet B.
\newblock Explicit error bounds in a conforming finite element method.
\newblock {\em Mathematics of Computation} 1999; \textbf{68}(288):1379--1396. 

\bibitem{FRA65}
Fraejis de Veubeke B.
\newblock Displacement and equilibrium models in the finite element method.
\newblock in {\em Stress Analysis}, O.C. Zienkiewicz and G.S. Holister, eds., New York, Wiley 1965 (pp.~145--197). 

\bibitem{LAD96}
Ladev\`eze P., Maunder E.A.W.
\newblock A general method for recovering equilibrating element tractions.
\newblock {\em Computer Methods in Applied Mechanics and Engineering} 1996; \textbf{137}:111--151. 

\bibitem{FLO02}
Florentin E., Gallimard L., Pelle J.-P.
\newblock Evaluation of the local quality of stresses in 3d finite element analysis.
\newblock {\em Computer Methods in Applied Mechanics and Engineering} 2002; \textbf{191}:4441--4457.

\bibitem{PLE12}
Pled F., Chamoin L., Ladev\`eze P.
\newblock An enhanced method with local energy minimization for the robust a posteriori construction of equilibrated stress fields in finite element analyses.
\newblock {\em Computational Mechanics} 2012; \textbf{49}:357--378.

\bibitem{LAD10c}
Ladev\`eze P., Chamoin L., Florentin E.
\newblock A new non-intrusive technique for the construction of admissible stress fields in model verification.
\newblock {\em Computer Methods in Applied Mechanics and Engineering} 2010; \textbf{199}(9-12):766--777.

\bibitem{DES98b}
Destuynder P., M\'etivet B.
\newblock Explicit error bounds for a nonconforming finite element method.
\newblock {\em SIAM Journal of Numerical Analysis} 1998; \textbf{35}(5):2099--2115. 

\bibitem{ERN07}
Ern A., Nicaise S., Vohralik M.
\newblock An accurate H(div) flux reconstruction for discontinuous Galerkin approximations of elliptic problems.
\newblock {\em C.R. Acad. Sci. Paris, S\'erie I} 2007; \textbf{345}(12):709--712. 

\bibitem{LAD97}
Ladev\`eze P., Rougeot P.
\newblock New advances on a posteriori error on constitutive relation in finite element analysis.
\newblock {\em Computer Methods in Applied Mechanics and Engineering} 1997; \textbf{150}:239--249.

\bibitem{BAB94b}
Babu\u{s}ka I., Strouboulis T., Upadhyay C.S., Gangaraj S.K., Copps K.
\newblock Validation of a posteriori error estimators by numerical approach.
\newblock {\em International Journal for Numerical Methods in Engineering} 1994; \textbf{37}(7):1073--1123.

\bibitem{MOE99}
Mo\"es N., Dolbow J., Belytchko T. 
\newblock A finite element method for crack growth without remeshing. 
\newblock {\em International Journal for Numerical Methods in Engineering} 1999; \textbf{46}:131--150.

\bibitem{Chamoin_Legoll_goal} 
Chamoin L., Legoll F. 
\newblock Goal-oriented error estimation and adaptivity in MsFEM computations.
\newblock {\em in preparation}.

\bibitem{XU07} 
Xu F.X.
\newblock A multiscale stochastic finite element method on elliptic problems involving uncertainties.
\newblock {\em Computer Methods in Applied Mechanics and Engineering} 2007; \textbf{196}:2723--2736. 

\bibitem{LEB14b}
Le Bris C., Legoll F., Thomines F.
\newblock Multiscale finite element approach for ``weakly'' random problems and related issues.
\newblock {\em Mathematical Modelling and Numerical Analysis} 2014; \textbf{48}(3):815--858.

\bibitem{EFE04} 
Efendiev Y., Hou T., Ginting V.
\newblock Multiscale finite element methods for nonlinear problems and their applications.
\newblock {\em Communications in Mathematical Sciences} 2004; \textbf{2}(4):553--589. 

\end{thebibliography}
\end{document}